# Lectures on the linearized Kapustin-Witten equations for [0,∞) × Y

Clifford Henry Taubes[†]

Department of Mathematics
Harvard University
Cambridge, MA 02318

chtaubes@math.harvard.edu

ABSTRACT: This is a written version of lectures that I would have given myself about aspects of the differential operator that is obtained from the linearized Kapustin-Witten equations on the product of the half-line with a compact, oriented, Riemannian 3-manifold. These lectures concern for the most part certain instances of much more general theorems of R. Mazzeo and E. Witten. There is also a 'lecture series' about the asymptotics of solutions to the same Kapustin-Witten equations as the half-line parameter limits to ∞.

[†]Supported in part by National Science Foundation grants DMS 1708310 and DMS 2002771

# 0. Introduction

This paper consists of an introductory pair of lecture and then three more series of lectures that I would have given to myself on aspects of a linear operator that is canonically associated to a solution to the SU(2) Kapustin-Witten equations on the product of the half-line with a compact, oriented three manifold.

The introductory two lectures introduces the operator in question (it comes via a linearization of the Kapustin-Witten equations about a given solution) and then some technology and conventions.

After the introductory parts, the first of the lecture series is an account of the lectures I would have given to myself about an instance of theorems by Rafe Mazzeo and Edward Witten in [MW1] and [MW2] concerning the operator on $(0, \infty) \times \mathbb{R}^2 \times S^1$ when defined using certain model solutions to the Kapustin-Witten equations presented by Edward Witten in [W].

The next lecture series is a written account of lectures I would have given about an instance of theorems of Mazzeo-Witten in [MW1] and [MW2] concerning the analogous differential operator on $(0, \infty) \times Y$ with Y being a compact, oriented 3-manifold. (Siqi He in [H] also discusses instances of these theorems. See also [HM1] and [HM2].)

The last lecture series is a written account of lectures that I would have given to myself about the behavior of solutions to the Kapustin-Witten equation on $[1, \infty) \times Y$ as the parameter t from the $[1, \infty)$ factor gets ever larger in an unbounded sequence. This final series of lectures gives a proof of a Kapustin-Witten version of a fundamental proposition proved by John Morgan, Tom Mrowka and Danny Ruberman [MMR] concerning the t → ∞ limits of finite energy solutions to the anti-self dual Yang-Mills equations on $[1, \infty) \times Y$.

The introduction and the first and second series of lectures differ for the most part from a specific case of what is written by Mazzeo-Witten only by virtue of an attempt to prove almost all assertions from first principles rather than by invoking on occasion theorems in the literature. (And my notation is different.) My lectures in the last lecture series differ from the what is written in part of [MMR] only to the extent that SU(2) is not SL(2; $\mathbb{C}$) and to the extend that I give a proof of the relevant instance of a fundamental theorem of Leon Simon [S].

There is very little that will be new to Mazzeo and Witten (and Siqi He) in the first two series of lectures, and very little that is new to Morgan, Mrowka and Ruberman in the last series of lectures.

**Table of Contents**

**Acknowledgements**

I am tremendously grateful for conversations and insightful feedback from Greg Parker whose questions prodded me to write down my own account of the Mazzeo-Witten theorems and the proposition of Morgan-Mrowka-Ruberman.

I am also delighted to say thank-you again to Rafe Mazzeo for arranging a most pleasant visit this past winter at Stanford; and I thank the faculty, staff and students at Stanford's Mathematics Department for their generous hospitality.

# I. THE DEFINITION OF THE OPERATOR $\mathcal{D}$

The section that follows directly introduces and describes the differential operator $\mathcal{D}$. The subsequent section introduces basic technology and some conventions/notation.

## 1. The operator

Let Y denote an oriented 3-manifold with a geodesically complete Riemannian metric, and let $P \to Y$ denote a principal SU(2) bundle (the associated adjoint bundle is denoted by ad(P)). For the purposes of these notes, the Kapustin-Witten equations are equations on $(0, \infty) \times Y$ for a pair consisting of a connection on a P and 1-form on $(0, \infty) \times Y$ with values ad(P). Such a pair is denoted by $(A, \mathfrak{a})$. The Kapustin-Witten equations are obeyed by $(A, \mathfrak{a})$ when

- $F_A - \mathfrak{a} \wedge \mathfrak{a} = \hat{*} D_A \mathfrak{a}$,
- $D_A \hat{*} \mathfrak{a} = 0$,

(I.1.1)

where the notation has $\hat{*}$ signifying the Hodge star for the product metric on $(0, \infty) \times Y$; it has $F_A$ denoting A's curvature tensor; and it has $D_A$ denoting the covariant exterior derivative defined by A.

Of interest in these lectures is the differential operator that is defined from the linearization of these equations. The notation for the upcoming definition of this operator uses t to denote the Euclidean coordinate on the $(0, \infty)$ factor of $(0, \infty) \times Y$. Also: The exterior covariant derivative defined by a connection (call it A) on ad(P)-valued differential forms is written as $dt \nabla_t + d_A$ with $d_A$ denoting the covariant exterior derivative defined by A along the Y-factor of $(0, \infty) \times Y$.

Let $(A, \mathfrak{a})$ denote a pair of connection on P over $(0, \infty) \times Y$ and ad(P)-valued 1-form on this same space. Of particular interest in what follows are the cases where $\mathfrak{a}$ has no dt component (so it annihilates the vector field $\frac{\partial}{\partial t}$). This extra condition is assumed henceforth. The linearized operator that is defined by $(A, \mathfrak{a})$ acts on pairs of the form $(B = \mathfrak{b}_t dt + \mathfrak{b}, C = \mathfrak{c}_t dt + \mathfrak{c})$ with B giving the first order change in A and C giving the first order change in $\mathfrak{a}$ if $(A, \mathfrak{a})$ is a solution to (I.1.1). This linearized operator with an additional gauge slice constraint is denoted by $\mathcal{D}$. It sends any given pair (B,C) to a pair (P, Q) with $P = \mathfrak{p}_t dt + \mathfrak{p}$ and $Q = \mathfrak{q}_t dt + \mathfrak{q}$ defined as follows:

- $\mathfrak{p} \equiv \nabla_t \mathfrak{b} - d_A \mathfrak{b}_t - *d_A \mathfrak{c} - *(\mathfrak{b} \wedge \mathfrak{a} + \mathfrak{a} \wedge \mathfrak{b}) + [\mathfrak{a}, \mathfrak{c}_t]$.
- $\mathfrak{p}_t \equiv \nabla_t \mathfrak{b}_t + (\nabla_i \mathfrak{b})_i + [\mathfrak{a}_i, \mathfrak{c}_i]$ .
- $\mathfrak{q} \equiv \nabla_t \mathfrak{c} - d_A \mathfrak{c}_t - *d_A \mathfrak{b} + *(\mathfrak{c} \wedge \mathfrak{a} + \mathfrak{a} \wedge \mathfrak{c}) - [\mathfrak{a}, \mathfrak{b}_t]$.
- $\mathfrak{q}_t \equiv \nabla_t \mathfrak{c}_t + (\nabla_i \mathfrak{c})_i - [\mathfrak{a}_i, \mathfrak{b}_i]$ .

(I.1.2)

The notation uses $\nabla$ to denote the full covariant derivative defined by the connection A and the Levi-Civita connection on Y. What is denoted by $\nabla_t$ signifies A's covariant derivative along the $(0, \infty)$ factor of $(0, \infty) \times Y$; what is denoted by $\{\nabla_i\}_{i=1,2,3}$ signifies the A's covariant derivative in the direction of a chosen orthonormal basis for TY. The subscript on $\mathfrak{a}$, $\mathfrak{c}$ and $\mathfrak{b}$ in the second and fourth bullets denote the corresponding components in these directions. Also, repeated Latin indices are meant to be summed implicitly over the index set $\{1, 2, 3\}$.

### a) Depicting $\mathcal{D}$

There are various ways to depict the operator $\mathcal{D}$. A first depiction: This operator can be written as the sum of three parts:

$$\mathcal{D} = \nabla_t + D(\cdot) + \{\mathfrak{a}, \cdot\}$$

(I.1.3)

where D contains the A-covariant derivatives in (I.1.2) along the Y factor and where $\{\mathfrak{a}, \cdot\}$ denotes the zero'th order terms that involve commutators of components of $\mathfrak{a}$ with components of B and C. The formal $L^2$ adjoint is

$$\mathcal{D}^\dagger = -\nabla_t + D(\cdot) + \{\mathfrak{a}, \cdot\} .$$

(I.1.4)

For the second depiction: If the pair (B, C) is written as a column vector using a given orthornormal frame for T*Y, thus

$$\begin{pmatrix} b_1 \\ b_2 \\ b_3 \\ b_t \\ c_1 \\ c_2 \\ c_3 \\ c_t \end{pmatrix},$$

(I.1.5)

then there is a corresponding depiction of $\mathcal{D}$ as a matrix:

$$\mathcal{D} = \left( \begin{array}{cccc|cccc} \partial_t & \mathfrak{a}_3 & -\mathfrak{a}_2 & -\nabla_1 & 0 & \nabla_3 & -\nabla_2 & \mathfrak{a}_1 \\ -\mathfrak{a}_3 & \partial_t & \mathfrak{a}_1 & -\nabla_2 & -\nabla_3 & 0 & \nabla_1 & \mathfrak{a}_2 \\ \mathfrak{a}_2 & -\mathfrak{a}_1 & \partial_t & -\nabla_3 & \nabla_2 & -\nabla_1 & 0 & \mathfrak{a}_3 \\ \nabla_1 & \nabla_2 & \nabla_3 & \partial_t & \mathfrak{a}_1 & \mathfrak{a}_2 & \mathfrak{a}_3 & 0 \\ \hline 0 & \nabla_3 & -\nabla_2 & -\mathfrak{a}_1 & \partial_t & -\mathfrak{a}_3 & \mathfrak{a}_2 & -\nabla_1 \\ -\nabla_3 & 0 & \nabla_1 & -\mathfrak{a}_2 & \mathfrak{a}_3 & \partial_t & -\mathfrak{a}_1 & -\nabla_2 \\ \nabla_2 & -\nabla_1 & 0 & -\mathfrak{a}_3 & -\mathfrak{a}_2 & \mathfrak{a}_1 & \partial_t & -\nabla_3 \\ -\mathfrak{a}_1 & -\mathfrak{a}_2 & -\mathfrak{a}_3 & 0 & \nabla_1 & \nabla_2 & \nabla_3 & \partial_t \end{array} \right)$$

(I.1.6)

In this regard, the following short hand is used here and below: A matrix entry given by a Lie algebra valued 1-form (such as $\mathfrak{a}_1$ or $\mathfrak{a}_2$ or $\mathfrak{a}_3$ above) indicates the commutation action of that 1-form on the relevant entry of the column vector on which the matrix acts. For example, $\mathfrak{a}_3$ above in the top row and second column acts on $b_2$ in (1.4) as $[\mathfrak{a}_3, b_2]$.

The third depiction of $\mathcal{D}$ uses six generators of an 8-dimensional, real Clifford algebra which are written as two sets of three, $\{\gamma_1, \gamma_2, \gamma_3\}$ and $\{\rho_1, \rho_2, \rho_3\}$. These are $8 \times 8$, anti-symmetric, traceless matrices that obey

$$\gamma_i \gamma_j + \gamma_j \gamma_i = -2\delta_{ij} \quad \text{and} \quad \gamma_i \rho_j + \rho_j \gamma_i = 0 \quad \text{and} \quad \rho_i \rho_j + \rho_j \rho_i = -2\delta_{ij} \, .$$

(I.1.7)

These two sets of matrices (the set of $\gamma$'s and the set of $\rho$'s) define two Clifford module structures for Y that are compatible with the Levi-Civita covariant derivative. A specific realization of these matrices are depicted in Appendix A. (This last depiction of $\mathcal{D}$ comes from Sergey Cherkis' octonionic description of the Kapustin-Witten equations [C].)

The operator $\mathcal{D}$ can be written using these Clifford matrices as:

$$\mathcal{D} = \nabla_t + \gamma_i \nabla_i + \rho_i [\mathfrak{a}_i, \cdot ] .$$

(I.1.8)

(Remember that repeated indices are implicitly summed over the indexing set $\{1, 2, 3\}$.) By way of comparison with (I.1.3): What is denoted by D in (1.3) is $\gamma_i \nabla_i$ and what is denoted by $\{\mathfrak{a}, \cdot\}$ in (1.3) is $\rho_i[\mathfrak{a}_i, \cdot]$.

**b) Hilbert space domain and range**

The operator $\mathcal{D}$ is a bounded operator from a range Hilbert space to be denoted by $\mathbb{H}$ to a domain Hilbert space to be denoted by $\mathbb{L}$. To set the notation for these definitions, introduce $\mathfrak{su}(2)$ to denote the vector space of $2 \times 2$, traceless, anti-Hermitian matrices; the Lie algebra of the group SU(2). The norm on this vector space is the positive square root of the function that is defined by the rule $\sigma \to |\sigma|^2 = -\frac{1}{2}\operatorname{trace}(\sigma\sigma)$. This norm comes from the inner product on $\mathfrak{su}(2)$ that is defined by the rule whereby the inner product of elements $\sigma$ and $\sigma'$ is $-\frac{1}{2}\operatorname{trace}(\sigma\sigma')$. It induces a corresponding inner product on the bundle ad(P) and, with the Riemannian inner product, an inner product on bundles of ad(P)-valued tensors. The Riemannian inner product and all of the latter are denoted by $\langle\,,\,\rangle$ in what follows

Let $\mathbb{W}$ denote the vector bundle $\oplus_2 (\operatorname{ad}(P) \oplus (\operatorname{ad}(P) \otimes T^*Y))$. The domain Hilbert space is the completion of the space of compactly supported elements in $C^\infty((0,\infty) \times Y; \mathbb{W})$ using the norm whose square is

$$\psi \to \int_{(0,\infty) \times Y} (|\nabla \psi|^2 + |[\mathfrak{a},\psi]|^2) .$$

(I.1.9)

This norm is denoted by $\|\cdot\|_{\mathbb{H}}$. Meanwhile, the range Hilbert space is the completion of the space of compactly supported elements in $C^\infty((0,\infty) \times Y; \mathbb{W})$ using the norm whose square is given by the rule

$$\eta \to \int_{(0,\infty) \times Y} |\eta|^2 .$$

(I.1.10)

This domain norm is denoted by $\|\cdot\|_{\mathbb{L}}$ (it is the standard $L^2$-norm).

The operator $\mathcal{D}$ gives a bounded linear operator from $\mathbb{H}$ to $\mathbb{L}$ precisely because it is a linear combination of covariant derivatives and commutators with components of $\mathfrak{a}$.

## 2. Basic technology and then two conventions

This section constitutes a digression to introduce some of the basic technology and some conventions for the subsequent lectures.

### a) Rewriting the equations in (I.1.1)

It proves useful to introduce some notation to rewrite the equations in (I.1.1) when $(A, \mathfrak{a})$ has no dt component. To this end, first introduce by way of notation $*$ to denote the Hodge star for the metric on Y. Then write the A's curvature 2-form $F_A$ as $dt \wedge E_A + *B_A$ with $E_A$ and $B_A$ denoting ad(P)-valued sections of T*Y on $(0, \infty) \times Y$. The equations in (I.1.1) are equivalent to these:

- $E_A = *d_A \mathfrak{a}$ .
- $\nabla_t \mathfrak{a} = B_A - *(\mathfrak{a} \wedge \mathfrak{a})$ .
- $d_A * \mathfrak{a} = 0$.

(I.2.1)

This depiction of (I.2.1) is more in tune with the product structure of $(0, \infty) \times Y$.

A warning: The linearized version of (I.2.1) does *not* lead to $\mathcal{D}$ because (I.2.1) depicts (I.1.1) when $\mathfrak{a}$ has zero dt component; as a consequence (I.2.1) doesn't know what to do with a first order deformation of $\mathfrak{a}$ that has non-zero dt component (this component is $\mathfrak{c}_t$). The operator that is depicted in (I.2.1) observedly does allows for a non-zero dt component of the first order deformation of $\mathfrak{a}$.

### b) Hardy's inequality

Hardy's inequality plays a crucial role in subsequent proofs. Actually, multiple versions do; but the present version is the generic one. This version of Hardy's inequality asserts that the function

$$\psi \to \int_{(0,\infty) \times Y} \tfrac{1}{t^2} |\psi|^2$$

(I.2.2)

which is a priori defined on compactly supported elements in $C^\infty((0,\infty) \times Y; \mathbb{W})$ extends to $\mathbb{H}$ as a bounded, continuous function that obeys

$$\int_{(0,\infty) \times Y} \tfrac{1}{t^2} |\psi|^2 \leq 4 \|\psi\|_{\mathbb{H}}^2 .$$

(I.2.3)

Hardy's inequality in (I.2.3) follows directly from its 1-dimensional incarnation which says this: If $f$ is a continuous, piece-wise $C^1$ functions on the interval $(0, \infty)$ with compact support, then

$$\int_{(0,\infty)} \tfrac{1}{t^2} f^2 \, dt \leq 4 \int_{(0,\infty)} |\tfrac{df}{dt}|^2 \, dt .$$

(I.2.4)

This is proved by writing $\tfrac{1}{t^2} f^2 dt$ from the left hand integral as $-f^2 d\tfrac{1}{t}$ and then integrating by parts to identify the left hand integral with the $(0,\infty)$ integral of $\tfrac{2}{t} f \tfrac{df}{dt}$. Having made this identification, then (I.2.4) follows directly via the Cauchy-Schwarz inequality. (The inequality in (I.2.3) follows from the $f = |\psi|$ version of (I.2.3) because $|\tfrac{d|\psi|}{dt}| \leq |\nabla_t \psi|$.)

### c) A Bochner-Weitzenboch formula

The Bochner-Weitzenboch formula for $\mathcal{D}^\dagger \mathcal{D}$ has the schematic form

$$\mathcal{D}^\dagger \mathcal{D} = \nabla^\dagger \nabla + [\mathfrak{a}_i [\,\cdot\,, \mathfrak{a}_i]] + \mathbb{X}$$

(I.2.5)

with $\nabla^\dagger$ denoting the formal adjoint of $\nabla$ (defined using $\mathbb{L}$'s inner product), and with $\mathbb{X}$ denoting a symmetric endomorphism of $\mathbb{W}$. The notation in (I.2.5) is such that the repeated indices in $[\mathfrak{a}_i, [\,\cdot\,, \mathfrak{a}_i]]$ are summed over the set $\{1, 2, 3\}$.

For the present purpose, there are two key observations to be made about $\mathbb{X}$, the first being that

$$|\mathbb{X}| \leq c_1 (|F_A| + |\nabla_A \mathfrak{a}|) + c_2 |\mathcal{R}_Y|$$

(I.2.6)

with $\mathcal{R}_Y$ denoting the Riemann curvature tensor of Y and with $c_1$ and $c_Y$ denoting positive numbers that are independent of $(A, \mathfrak{a})$ and the metric on Y.

The second key observation about $\mathbb{X}$ is this: This endomorphism annihilates elements of $\mathbb{W}$ with only the $\mathfrak{c}_t$ component in (I.1.6) non-zero.

An integral version of the the Bochner-Weitzenboch formula from (I.2.5) holds for elements in $\mathbb{H}$ if the norm of $\mathbb{X}$ is bounded by a constant multiple of $\tfrac{1}{t^2}$, thus if

$$|\mathbb{X}| \leq c_0 \tfrac{1}{t^2}$$

(I.2.7)

with $c_0$ being constant. The following is true under this assumption: If $\psi$ and $\xi$ are any two elements in $\mathbb{H}$, then

$$\int_{(0,\infty) \times \mathbb{R}^2 \times S^1} \langle \mathcal{D}\xi, \mathcal{D}\psi \rangle = \int_{(0,\infty) \times \mathbb{R}^2 \times S^1} (\langle \nabla \xi, \nabla \psi \rangle + \langle [\mathfrak{a}, \xi], [\mathfrak{a}, \psi] \rangle) + \int_{(0,\infty) \times \mathbb{R}^2 \times S^1} \langle \xi, \mathbb{X}\psi \rangle$$

(I.2.8)

Remember that notation used here has $\langle \cdot , \cdot \rangle$ denoting the pointwise inner product on the relevant product vector bundles over $(0, \infty) \times Y$ that is induced from the inner product $\langle \cdot \, \cdot \rangle$ on $\mathfrak{su}(2)$ (which is $-\frac{1}{2}$ trace($\cdot \, \cdot$)).

With regards to (I.2.7): The left most integral on the right in (I.2.7) is the inner product in $\mathbb{H}$ between $\psi$ and $\xi$. Meanwhile, the right most integral on the right side of (I.2.7) is a priori bounded by a $\psi$ and $\xi$ independent multiple of $\|\psi\|_{\mathbb{H}} \|\xi\|_{\mathbb{H}}$ by virtue of the Hardy's inequality (I.2.3).

### d) Notation and conventions

The first convention is with regards to notation: The notation has $c$ denoting a number that is greater than 1 and independent of what ever is relevant to a given inequality. For example, if an inequality concerns a value or values of t from $[0, \infty)$ and/or points in Y, then $c$ will be independent of t and the points in y. In general, it should be clear from the context what does and doesn't determine an upper bound for $c$. It is always the case that $c$ increases between successive appearances.

The second convention is with regards to what are called 'cut-off functions'. All such functions will be constructed from a basic model function on $\mathbb{R}$ to be denoted by $\chi$. This is a smooth, non-increasing function that equals 1 on $(-\infty, \frac{1}{4}]$ and equals zero on $[\frac{3}{4}, \infty)$. Cut-off functions in other contexts on other manifolds (such as on $(0, \infty) \times Y$) can be obtained from $\chi$ by composing the latter with a suitably chosen map to $\mathbb{R}$. (The advantage of such a construction is this: The norms of the derivatives of these other cut off functions have a priori bounds given bounds for the norms of the derivatives of the map to $\mathbb{R}$.)

## II. THE OPERATOR $\mathcal{D}$ ON $(0, \infty) \times \mathbb{R}^2 \times S^1$

The sections in this 'lecture series' concern the version of $\mathcal{D}$ on $(0, \infty) \times \mathbb{R}^2 \times S^1$ that is defined by certain model solutions to (I.2.1) on this domain. The goal of the sections is to first state and then prove a fundamental theorem (the upcoming Theorem II.1) to the effect that each such version of $\mathcal{D}$, if viewed as map from the corresponding version of $\mathbb{H}$ to $\mathbb{L}$, defines a Fredholm operator with trivial kernel and cokernel.

### 1. The operator $\mathcal{D}$ for the model solutions

This section presents Witten's model solutions and then the fundamental theorem regarding the associated version of $\mathcal{D}$.

### a) Witten's model solutions

Witten's model solutions to (I.2.1) on $(0, \infty) \times \mathbb{R}^2 \times S^1$ are depicted momentarily. To set the stage the depiction: The metric used here for $\mathbb{R}^2 \times S^1$ (which is Y in this case) is the product of the Euclidean metric on $\mathbb{R}^2$ with a metric on $S^1$. The depiction of the model solutions uses an orthonormal basis $\{\sigma_1, \sigma_2, \sigma_3\}$ for the Lie algebra of SU(2) (which is denoted by $\mathfrak{su}(2)$) that obey:

- $\sigma_1^2 = \sigma_2^2 = \sigma_3^2 = -1$
- $\sigma_1 \sigma_2 = -\sigma_3$ *and* $\sigma_2 \sigma_3 = -\sigma_1$ *and* $\sigma_3 \sigma_1 = -\sigma_2$ .

(II.1.1)

Note the - sign in the second bullet. This is not the convention used by other people.

The formulas that follow use $\theta_0$ to denote the product connection on the product principal SU(2) bundle over $(0, \infty) \times \mathbb{R}^2 \times S^1$. The basis $\{\sigma_1, \sigma_2, \sigma_3\}$ for $\mathfrak{su}(2)$ is viewed as a $\theta_0$-covariantly constant basis for the associated Lie algebra bundle of the product principal SU(2) bundle.

The notation used below implicitly refers to Euclidean coordinates $(z_1, z_2)$ for the $\mathbb{R}^2$ factor of $\mathbb{R}^2 \times S^1$. It uses these to identify $\mathbb{R}^2$ with $\mathbb{C}$ using the $\mathbb{C}$-valued coordinate $z = z_1 + i z_2$. It also uses the coordinate differentials $dz_1$ and $dz_2$ with the differential of the Euclidan coordinate $x_3$ for the $S^1$ factor as an orthonormal frame for $T^*(\mathbb{R}^2 \times S^1)$. (The coordinate $x_3$ is $\mathbb{R}/(\ell \mathbb{Z})$-valued with $\ell$ denoting the length of the $S^1$ factor.)

The notation also introduces the real valued function $\Theta$ on $(0, \infty) \times (\mathbb{C} - 0)$ that is defined by setting $\sinh(\Theta) = \frac{t}{|z|}$. Meanwhile, a function on the $(0, \infty) \times \mathbb{R}^2$ factor of $(0, \infty) \times \mathbb{R}^2 \times S^1$ that is denoted by $x$ is $(t^2 + |z|^2)^{1/2}$.

Witten's model solutions (from [W] are indexed by a non-negative integer which is denoted by $m$ in what follows. The integer $m$ version is defined as follows:

- *Higgs field*:
    a) $\mathfrak{a}_3 = -\frac{1}{2t} \frac{(m+1)\sinh(\Theta)}{\sinh((m+1)\Theta)} \frac{\cosh((m+1)\Theta)}{\cosh(\Theta)} \sigma_3$.
    b) $\mathfrak{a}_1 - i\mathfrak{a}_2 = -\frac{1}{2t} \frac{(m+1)\sinh(\Theta)}{\sinh((m+1)\Theta)} (\frac{z}{|z|})^m (\sigma_1 - i\sigma_2)$ .
- *Connection and curvature*:
    a) $A = \theta_0 + \frac{(m+1)}{2}(1 - \frac{\sinh(\Theta)}{\cosh(\Theta)} \frac{\cosh((m+1)\Theta)}{\sinh((m+1)\Theta)}) \frac{1}{|z|^2} (z_1 dz_2 - z_2 dz_1) \sigma_3$ .
    b) $B = \frac{(m+1)}{2x^2} \frac{\sinh(\Theta)}{\cosh(\Theta)} \frac{\cosh((m+1)\Theta)}{\sinh((m+1)\Theta)} (1 - \frac{(m+1)\sinh(\Theta)\cosh(\Theta)}{\sinh((m+1)\Theta)\cosh((m+1)\Theta)}) \sigma_3 \, dx_3$
    c) $E = -\frac{(m+1)}{2x^2} \frac{\cosh((m+1)\Theta)}{\sinh((m+1)\Theta)} (1 - \frac{(m+1)\sinh(\Theta)\cosh(\Theta)}{\sinh((m+1)\Theta)\cosh((m+1)\Theta)}) \sigma_3 \frac{1}{x} (z_1 dz_2 - z_2 dz_1)$ .

(II.1.2)

These solutions have certain special properties that play roles in what is to come. The salient ones are listed below in (II.1.3)–(II.1.4).

- *The Lie-algebra element $\sigma_3$ is A-covariantly constant; so the connection A is Abelian.*
- *The component $\mathfrak{a}_3$ is proportional to $\sigma_3$; it will be written as $\alpha\sigma_3$.*
    a) *$\alpha$ is negative with values between $-\frac{1}{2t}$ and $-\frac{(m+1)}{2t}$.*
    b) *$\frac{\partial}{\partial t}\alpha$ is positive and greater than a positive, constant multiple of $\frac{1}{t^2}$.*
- *The components $\mathfrak{a}_1$ and $\mathfrak{a}_2$ are pointwise orthogonal to $\sigma_3$. Moreover*
    a) *They are pointwise orthogonal to each other and they have the same norm. Thus, the square of $\varphi = \mathfrak{a}_1 - i\mathfrak{a}_2$ has trace zero*
    b) *$[\frac{i}{2}\sigma_3, \varphi] = \varphi$.*
    c) *$|\varphi| \leq \frac{1}{\sqrt{2}t}$ with equality only in the case when $m = 0$.*
- *The $\mathfrak{su}(2)$ valued 1-forms $B \equiv B_A$ and $E \equiv E_A$ and the components of $\nabla\mathfrak{a}$ obey*
    a) *$B_1 \equiv B_2 \equiv E_3 \equiv \nabla_3\mathfrak{a} \equiv 0$.*
    b) *The norms of $B_3$, $E_1$, $E_2$ and $\nabla_1\mathfrak{a}$, $\nabla_2\mathfrak{a}$ are bounded by a constant multiple of $\frac{t}{x^3}$.*
    c) *The norm of $\nabla_t\mathfrak{a}$ is bounded by a constant multiple of $\frac{1}{t^2}$.*

(II.1.3)

The next property concerns the family of coordinate rescaling diffeomorphisms of $(0,\infty) \times \mathbb{R}^2$ parametrized by $(0,\infty)$ whose action is defined by the rule whereby the diffeomorphism that is parametrized by any given $\lambda > 0$ sends $(t, z_1, z_2) \to (\lambda t, \lambda z_1, \lambda z_2)$.

*The model solutions from (II.1.2) are fixed by every coordinate rescaling diffeomorphism from this family: Each model solution is the same as its pull-back by any of them.*

(II.1.4)

A final property: Write $\mathfrak{a}_3$ as $\alpha\sigma_3$ with $\alpha \equiv -\frac{1}{2t}\frac{(m+1)\sinh(\Theta)}{\sinh((m+1)\Theta)}\frac{\cosh((m+1)\Theta)}{\cosh(\Theta)}$. Then the equations in (I.2.1) when written using $(\alpha, \varphi)$ and $(E_1, E_2, B_3)$ say this:

- $\nabla_t\varphi - 2\alpha\varphi = 0$ *and* $(\nabla_1 + i\nabla_2)\varphi = 0$.
- $E_1 = \frac{\partial\alpha}{\partial z_2}\sigma_3$ *and* $E_2 = -\frac{\partial\alpha}{\partial z_1}\sigma_3$.
- $B_3 = (\frac{\partial\alpha}{\partial t} - |\varphi|^2)\sigma_3$.

(II.1.5)

By way of terminology: The $m = 0$ version of (II.1.2) is called the Nahm pole solution; it is the solution with $A = \theta_0$ and $\mathfrak{a}_i = -\frac{1}{2t}\sigma_i$.

### b) The operator $\mathcal{D}$ for Witten's model solutions

The central concern of this lecture series is the following theorem which is an instance of a far more general theorem from the paper [MW2] by Rafe Mazzeo and Edward Witten:

**Theorem II.1**: *The operator $\mathcal{D}$ when defined using any of Witten's model solutions is a bounded, Fredholm map from $\mathbb{H}$ to $\mathbb{L}$ with trivial kernel and trivial cokernel.*

The proof of Theorem II.1 occupies Sections II.2-II.6 of these notes.

With regards to the theorem: A bounded operator between Banach spaces is Fredholm when it has the following properties:

- *The kernel is finite dimensional.*
- *The range is closed.*
- *The cokernel is finite dimensional.*

(II.1.6)

A proof that $\mathcal{D}$ is Fredholm as a map from $\mathbb{H}$ and $\mathbb{L}$ must demonstrate that all three of these conditions are met. That is what these notes will do when $\mathcal{D}$ is defined by one of Witten's model solutions.

The arguments presented here for the three conditions in (II.1.6) amount to little more than a distillation of various observations in [MW1] and [MW2] that are used to prove [MW2]'s more general theorem. Even so, these notes aspire to be self-contained to the extent possible. In particular, most of the analytic observations are proved directly, even instances of more general theorems in the literature.

By way of a look ahead, the only truly subtle part of the argument concerns the proof that $\mathcal{D}$ has closed range. In this regard, the assertion that $\mathcal{D}$ has closed range follows tautologically from an inequality of the form

$$\|\cdot\|_{\mathbb{H}} \leq \kappa \|\mathcal{D}(\cdot)\|_{\mathbb{L}}$$

(II.1.7)

with $\kappa$ being a real number.

If (II.1.7) holds, then the kernel of $\mathcal{D}$ in $\mathbb{H}$ is also trivial. The assertion that the cokernel is trivial (assuming (II.1.7)) is (almost) a direct consequences of the coordinate rescaling property in (II.1.4). First of all, if the range is closed (implied by (II.1.7)), then the cokernel of $\mathcal{D}$ is the kernel in $\mathbb{L}$ of its formal $\mathbb{L}$ adjoint, $\mathcal{D}^\dagger$ (see (I.1.4).) Assuming that $\mathcal{D}$ has closed range, and granted that (II.1.4) holds, then the respective kernels in $\mathbb{L}$ of $\mathcal{D}^\dagger$ (and also $\mathcal{D}$) are spanned by elements that have the form $x^\lambda u$ with $\lambda$ being a real number and $u$ being independent of $x$. Since no element of this sort is in $\mathbb{L}$, both $\mathcal{D}^\dagger$ and $\mathcal{D}$ must have trivial kernel in $\mathbb{L}$.

## 2. Basic technology for the case of the model solutions

This section introduces some basic technology that is specific to the case when Y is $(0,\infty) \times \mathbb{R}^2 \times S^1$ and when $(A, \mathfrak{a})$ is one of Witten's model solutions.

### a) Versions of Hardy's inequality

Hardy's inequality as depicted in (I.2.3) plays a role in the cases at hand, but so do two somewhat more subtle versions of Hardy's inequality. To set the stage for the first of these other versions, reintroduce Euclidean coordinates $(z_1, z_2)$ for the $\mathbb{R}^2$ factor of $\mathbb{R}^2 \times S^1$ and then the function $x = (t^2 + z_1^2 + z_2^2)^{1/2}$. This second version of Hardy's inequality says that the function

$$\psi \to \int_{(0,\infty)\times\mathbb{R}^2\times S^1} \tfrac{1}{x^2} |\psi|^2$$

(II.2.1)

which is a priori defined on compactly supported elements in $C^\infty((0,\infty)\times\mathbb{R}^2\times S^1; \mathbb{W})$ extends to $\mathbb{H}$ as a bounded, continuous function obeying

$$\int_{(0,\infty)\times\mathbb{R}^2\times S^1} \tfrac{1}{x^2} |\psi|^2 \leq \tfrac{4}{9} \|\psi\|_{\mathbb{H}}^2 .$$

(II.2.2)

For the purpose of proving this version of Hardy's inequality, suppose that $f$ is a continuous, piece-wise $C^1$ function on $(0,\infty)$ with compact support. Then an integration by parts and the Cauchy-Schwarz inequality much like the one used in proving (I.2.4) leads to this:

$$\int_{(0,\infty)} f^2 \, dx \leq 4 \int_{(0,\infty)} (\tfrac{\partial f}{\partial x})^2 x^2 \, dx .$$

(II.2.3)

Keeping (II.2.3) handy, suppose instead that $f$ is a continuous, piece-wise $C^1$ function on $[-\pi, \pi]$ that vanishes at both endpoints. There is in this instance the inequality

$$\int_{[\pi,-\pi]} f^2 \sin\theta \, d\theta \leq \tfrac{1}{2} \int_{[\pi,-\pi]} (\tfrac{df}{d\theta})^2 \sin\theta \, d\theta$$

(II.2.4)

which follows by virtue of the fact that 2 is the smallest Dirichelet eigenvalue of the unbounded, postive, self adjoint operator $-\tfrac{1}{\sin\theta}\tfrac{d}{d\theta}(\sin\theta \tfrac{d}{d\theta}(\cdot))$ acting on the space of square integrable function on the interval $[-\pi, \pi]$. (The corresponding eigenfunction is the function $\theta \to \cos\theta$.) The sum of the respective $f = |\psi|$ versions of these last two inequalities lead directly to (II.2.2) because $|d|\psi|| \leq |\nabla\psi|$.

To set the stage for the next version of Hardy's inequality: Let $S^+$ denote the hemisphere in $(0, \infty) \times \mathbb{R}^2$ where $x = 1$. Let A denote a connection on the restriction of P to $S^+$ and let $\nabla^S$ denote the corresponding covariant derivative, a map from sections of

ad(P) over $S^+$ to sections of $ad(P) \otimes T^*S^+$. Now let $\eta$ denote a section of ad(P)'s restriction to $S^+$ with compact support. Version three of Hardy's inequality says that

$$\int_{S^+} \tfrac{1}{t^2}|\eta|^2 \leq 4 \int_{S^+} |\nabla^S \eta|^2 \ .$$

(II.2.5)

This inequality can be derived from (II.2.4) and (I.2.3). Here is another way: Write t using the functions $x$ and $\Theta$ as $x \tfrac{\sinh\Theta}{\cosh\Theta}$. Then (II.2.5) follows from the following inequality for bounded functions on $(0, \infty)$ that vanish near 0:

$$\int_{(0,\Theta)} \tfrac{1}{\sinh^2\Theta} f^2 \, d\Theta \leq 4 \int_{(0,\infty)} \tfrac{1}{\cosh^2\Theta} |\tfrac{\partial f}{\partial \Theta}|^2 \, d\Theta \ .$$

(II.2.6)

(In this regard: The round metric on $S^+$ when written using the function $\Theta$ with the usual longitudinal angle as coordinates is conformally flat with conformal factor $\tfrac{1}{\cosh^2\Theta}$.) To prove (II.2.6), write $\tfrac{1}{\sinh^2\Theta}$ on the left hand side as $-\tfrac{1}{\cosh\Theta} d(\tfrac{1}{\sinh\Theta})$, integrate by parts and then use the Cauchy-Schwarz inequality. (The integration by parts term with the derivative of $\tfrac{1}{\cosh\Theta}$ can be discarded because of its sign.)

**b) The Bochner-Weitzenboch formula again**

The upcoming equation (II.2.7) depicts the endomorphism $\mathbb{X}$ that appears in (I.2.5) when Y is $\mathbb{R}^2 \times S^1$ and when $(A, \mathfrak{a})$ is one of Witten's model solutions.

With regards to notation: The notation has $x_3$ denoting the Euclidean coordinate along the $S^1$ factor. Meanwhile, Euclidean coordinates $(z_1, z_2)$ are taken for the $\mathbb{R}^2$ factor. Having specified a connectoin A on the principal bundle P, then $\nabla_3$ denotes A's directional covariant derivative along the unit length, tangent vector to the $S^1$ factor of $S^1$-factor of $\mathbb{R}^2 \times S^1$. And, $\nabla_1$ and $\nabla_2$ are the A-covariant directional derivatives as defined by the Euclidean coordinates $(z_1, z_2)$ for the $\mathbb{R}^2$ factor. Supposing that $\mathfrak{a}$ is a 1-form on $\mathbb{R}^2 \times S^1$ with coefficients in ad(P), then $\mathfrak{a}_1$, $\mathfrak{a}_2$ and $\mathfrak{a}_3$ denote its respective $dz_1$, $dz_2$ and $dx_3$ components.

The upcoming Bochner-Weitzenboch formula holds when $(A, \mathfrak{a})$ is an $S^1$-invariant solutions to the Kapustin-Witten equations on $\mathbb{R}^2 \times S^1$. This means in particular that $\nabla_3 \mathfrak{a} = 0$ and that $\nabla_3$ commutes with $\nabla_t$ and $\nabla_1$ and $\nabla_2$. This is to say that $E_3$ and $B_1$ and $B_2$ vanish identically. These assumptions with regards to $\nabla_3$ hold for Witten's model solutions.

The following Bochner-Weitzenboch depiction of $\mathbb{X}$ can be derived using any of the three depictions of $\mathcal{D}$ given in Section I.1a:

$$\begin{pmatrix} 0 & -2B_3 & 0 & 2E_1 & -\mathfrak{A}_{11} & -\mathfrak{A}_{12} & 2E_2 & 0 \\ 2B_3 & 0 & 0 & 2E_2 & -\mathfrak{A}_{21} & -\mathfrak{A}_{22} & -2E_1 & 0 \\ 0 & 0 & 0 & 0 & 0 & 0 & 0 & 0 \\ -2E_1 & -2E_2 & 0 & 0 & 2[\mathfrak{a}_2,\mathfrak{a}_3] & 2[\mathfrak{a}_3,\mathfrak{a}_1] & 2[\mathfrak{a}_1,\mathfrak{a}_2]-2B_3 & 0 \\ \mathfrak{A}_{11} & \mathfrak{A}_{12} & 0 & 2[\mathfrak{a}_3,\mathfrak{a}_2] & 0 & 2[\mathfrak{a}_2,\mathfrak{a}_1] & 2[\mathfrak{a}_3,\mathfrak{a}_1] & 0 \\ \mathfrak{A}_{21} & \mathfrak{A}_{22} & 0 & 2[\mathfrak{a}_1,\mathfrak{a}_3] & 2[\mathfrak{a}_1,\mathfrak{a}_2] & 0 & 2[\mathfrak{a}_3,\mathfrak{a}_2] & 0 \\ -2E_2 & 2E_1 & 0 & 2[\mathfrak{a}_2,\mathfrak{a}_1]+2B_3 & 2[\mathfrak{a}_1,\mathfrak{a}_3] & 2[\mathfrak{a}_2,\mathfrak{a}_3] & 0 & 0 \\ 0 & 0 & 0 & 0 & 0 & 0 & 0 & 0 \end{pmatrix}$$

(II.2.7)

By way of notation, each entry in this matrix acts on the relevant entry in (I.1.5) by commutation. For example, when $\mathbb{X}$ acts on the column vector depicted in (I.1.5) to give a new column vector, the top entry of this new column vector is the sum of $-2[B_3, \mathfrak{b}_2]$ and $2[E_1, \mathfrak{b}_t]$ and $-[\mathfrak{A}_{11}, \mathfrak{c}_1]$ and $-[\mathfrak{A}_{12}, \mathfrak{c}_2]$ and $2[E_2, \mathfrak{c}_3]$. In this regard, the Lie algebra valued functions $E_1$ and $E_2$ and $B_3$ are given by the rules $E_1 = [\nabla_t, \nabla_1]$ and $E_2 = [\nabla_r, \nabla_2]$; and $B_3 = [\nabla_1, \nabla_2]$; they are the non-zero components of A's curvature tensor. Also: What is denoted by $\mathfrak{A}_{11}$ is short-hand for $\nabla_1 \mathfrak{a}_1$, what is denoted by $\mathfrak{A}_{22}$ is $\nabla_2 \mathfrak{a}_2$, and what are denoted by $\mathfrak{A}_{12}$ and $\mathfrak{A}_{21}$ are $\nabla_1 \mathfrak{a}_2$ and $\nabla_2 \mathfrak{a}_1$. Note that $\nabla_1 \mathfrak{a}_1 = -\nabla_2 \mathfrak{a}_2$ and $\nabla_1 \mathfrak{a}_2 = \nabla_2 \mathfrak{a}_1$ when $(A, \mathfrak{a})$ is an $S^1$ invariant solution. (Since the metric on $\mathbb{R}^2 \times S^1$ is flat, there are no Riemann curvature terms contributing to $\mathbb{X}$.)

Concerning Witten's model solutions and the integral version of the Bochner-Weitzenboch formula in (I.2.8): This formula can be invoked using (II.2.7)'s version of $(A, \mathfrak{a})$ because of what is said by the fourth bullet in (II.1.3). (The norm of the corresponding version of $\mathbb{X}$ is bounded by a constant multiple of $\frac{1}{t^2}$.)

**c) Second order operators on the half-sphere**

Let $S^+$ denote the $t > 0$ part of the sphere in $(0, \infty) \times \mathbb{R}^2$ where $x = 1$ (which is the locus where $t^2 + |z|^2 = 1$). Suppose in what follows that $\mathbb{V} \to S^+$ is a given vector bundle with fiber metric. Fix a metric compatible connection and let $\nabla^S$ denote the induced covariant derivative on the space $C^\infty(S^+; \mathbb{V})$. (It maps a section of $\mathbb{V}$ to a section of $\mathbb{V} \otimes T^*S^2$.) Let W denote a symmetric endomorphism of $\mathbb{V}$. The analysis in this subsection concerns the bilinear form that is defined initially on the space of compactly supported sections of $\mathbb{V}$ given below in (II.2.8). This form is denoted by $\mathcal{E}$. Supposing that $\psi$ and $\eta$ are two compactly supported sections, then

$$\mathcal{E}(\psi, \eta) = \int_{S^+} (\langle \nabla^S \psi, \nabla^S \eta \rangle + \langle \psi, W\eta \rangle)$$

(II.2.8)

The notation here has $\langle\,,\,\rangle$ denoting both the fiber inner product on $\mathbb{V}$ and the fiber inner product on $\mathbb{V}\otimes T^*S^2$ that comes from the one on $\mathbb{V}$ and the round metric's inner product on $T^*S^2$. Let $\mathbb{H}_S$ denote the completion of the space of compactly supported sections of $\mathbb{V}$ using the norm whose square is given by the rule

$$\psi \to \int_{S^+} \langle \nabla^S\psi, \nabla^S\psi \rangle \;.$$

(II.2.9)

(This $\mathbb{H}_S$ is a separable Hilbert space whose inner product is the polarization of the square of the norm.) The upcoming Proposition II.2.1 says in part that if W is bounded from below and if its norm is such that $t^2|W|$ is bounded, then $\mathcal{E}$ extends from its dense domain to define a bounded, bilinear form on $\mathbb{H}_S$.

Use $\mathbb{L}_S$ to denote the Hilbert space completion of the space of compactly supported sections of $\mathbb{V}$ using the norm whose square is given by the rule

$$\psi \to \int_{S^+} |\psi|^2 \;.$$

(II.2.10)

This is the classical $L^2$ norm and $\mathbb{L}_S$ is the classical $L^2$–Hilbert space. The forgetful map from $\mathbb{H}_S$ to $\mathbb{L}_S$ is bounded (this follows from the $f = |\psi|$ and $\psi \in \mathbb{H}_S$ version of (II.2.4)).

A section of $\psi$ in $\mathbb{H}_S$ is said to be an eigensection of $\mathcal{E}$ if

$$\mathcal{E}(\psi,\eta) = \lambda \int_{S^+} \langle \psi,\eta \rangle$$

(II.2.11)

for all $\eta \in \mathbb{H}_S$. The real number $\lambda$ is the corresponding eigenvalue.

**Proposition II.2.1**: *Assume that $t^2|W|$ is bounded on $S^+$ and that W is bounded from below. Then $\mathcal{E}$ extends to $\mathbb{H}_S$ as a bounded, bilinear form. Moreover,*
- *There is an orthonormal basis of $\mathbb{L}_S$ from $\mathbb{H}_S$ that are eigensections for $\mathcal{E}$.*
- *The corresponding set of eigenvalues is bounded from below, has no accumulation points and each eigenvalue has finite multiplicity.*

*Proof of Proposition II.2.1*: To see that $\mathcal{E}$ extends to $\mathbb{H}_S$: There are two key points to proving this: The first is that

$$\int_{S^+} \langle \nabla^S \psi, \nabla^S \psi \rangle \geq 2 \int_{S^+} |\psi|^2$$

(II.2.12)

when ψ is a compactly supported section of $\mathbb{V}$. This follows from (II.2.4) by taking $f = |\psi|$. The second point is that

$$\int_{S^+} \langle \nabla^S \psi, \nabla^S \psi \rangle \geq \tfrac{1}{4} \int_{S^+} \tfrac{1}{t^2} |\psi|^2$$

(II.2.13)

which holds when ψ is again a compactly supported section of $\mathbb{V}$. This follows from the third version of Hardy's inequality, the one in (II.2.5).

Given the inequality in (II.2.13) and given that W is bounded from below, then the existence of a complete orthonormal basis of eigensections for $\mathcal{E}$ with the asserted properties of the eigenvalue set is proved by mimicking the proof of the analogous assertion for the Dirichelet eigenfunctions for the Laplacian on $S^+$ via the bilinear form

$$(f, f\prime) \to \int_{S^+} \langle df, df\prime \rangle \ .$$

(II.2.14)

Indeed, the linear algebra for the latter proof and for the proof of the bullets in Propostion II.2.1 is the same as that used for proving that a symmetric bilinear form on a finite dimensional vector space has a complete set of orthonormal eigenvectors. The only novel input for these infinite dimensional analogs is the Rellich theorem (see for example [F], page 305) which asserts in this case that the tautological map from the completion of the space of compactly supported functions on $S^+$ with finite $|d(\cdot)|^2$ integral to the space of square integrable functions on $S^+$ takes sequences with bounded $|d(\cdot)|^2$ integral to sequences with convergent subsequences with respect to the $L^2$ inner product. This implies that the limit of any weakly convergent sequence in $\mathbb{H}_S$ with the integral in (II.2.10) equal to 1 for each member will be an element in $\mathbb{H}_S$ with the integral in (II.2.10) also equal to 1.

Some relevant consequences of the Rellich lemma are stated below. These use $\mathbb{H}_{S1}$ to denote the subspace in $\mathbb{H}_S$ where the integral in (II.2.10) is equal to 1.

- Because W is bounded from below, the function $\psi \to \mathcal{E}(\psi, \psi)$ is bounded from below on $\mathbb{H}_{S1}$.
- Since $\mathbb{H}_S$ is a separable Hilbert space, a bounded sequence has a weakly convergent subsequence. As a consequence of the Rellich lemma, the limit is in $\mathbb{H}_{S1}$ if the sequence is in $\mathbb{H}_{S1}$.

- Because W is bounded from below, a minimizing sequence for $\mathcal{E}$ (or a minimizing sequence for the restriction of $\mathcal{E}$ to a closed subspace of $\mathbb{H}_{S1}$) has bounded $\mathbb{H}_S$-norm.
- As a consequence of weak convergence (and strong convergence of the integrals in (II.2.10)), the value of $\mathcal{E}(\cdot, \cdot)$ on the weak limit of a bounded sequence in $\mathbb{H}_{S1}$ is no greater than the lim-inf of its values on the subsequence.

**d) Elliptic regularity**

This subsection offers a technical lemma regarding elements in $\mathbb{L}$ that are annihilated by $\mathcal{D}^\dagger$ (or $\mathcal{D}$). But first a preliminary comment: These elements from $\mathbb{L}$ are a priori smooth on $(0, \infty) \times \mathbb{R}^2 \times S^1$. The promised lemma gives something by way of an a priori integral bound on the first derivatives:

**Lemma II.2.2**: *Supposing that $\mathcal{D}$ is defined by one of Witten's model solutions, then there exists $\kappa > 1$ with the following significance: If $\eta$ is from the $\mathbb{L}$-kernel of $\mathcal{D}^\dagger$ (or $\mathcal{D}$), and if $t \in (0, \infty)$, then*

$$\int_{[t, 2t] \times \mathbb{R}^2 \times S^1} (|\nabla \eta|^2 + \|[\mathfrak{a}, \eta]\|^2) \leq \frac{\kappa}{t^2} \int_{[\frac{1}{2}t, 4t] \times \mathbb{R}^2 \times S^1} |\eta|^2 \ .$$

*Proof of Lemma II.2.2*: Remember that a smooth, non-increasing function on $\mathbb{R}$ has been fixed (denoted by $\chi$) which is equal to 1 on $(-\infty, \frac{1}{4}]$ and equal to 0 on $[\frac{3}{4}, \infty)$. Given $R > 1$, use $\chi_R$ to denote the function on $\mathbb{R}^2$ that is given by the rule $z \to \chi(\frac{|z|}{R} - 1)$. This function is equal to 1 where $|z|$ is less than $R$ and equal to 0 where $|z|$ is greater than $2R$. A cut-off function depending on t is also needed: For this purpose, fix for the moment $\tau \in (0, \infty)$, let $X_\tau$ denote the function on $(0, \infty)$ that is defined by $t \to \chi(\frac{t}{2\tau} - 1)\chi(2(1 - \frac{t}{\tau}))$. This function has support where $t \in [\frac{1}{2}\tau, 4\tau]$ and it equals one when $t \in [\tau, 2\tau]$.

The operator $\mathcal{D}^\dagger$ has a Bochner-Weitzenboch formula which has the same general form as the one depicted in (I.2.5). For the present purposes, it is enough to know that

$$\mathcal{D}\mathcal{D}^\dagger = \nabla^\dagger \nabla + [\mathfrak{a}_i, [\cdot, \mathfrak{a}_i]] + \mathbb{X}'$$

(II.2.15)

with $\mathbb{X}'$ being an endomorphism whose norm is bounded by a contant multiple of $\frac{1}{t^2}$. Supposing that $\eta$ is in the $\mathbb{L}$-kernel of $\mathcal{D}^\dagger$, use this formula for $\mathcal{D}\mathcal{D}^\dagger \eta$. Since the left hand side is zero, so is the right. In any event, take the inner product of the right hand side with $\chi_R^2 X_\tau^2 \eta$ and then integrate the result. Then integrate by parts and use the Cauchy-Schwarz inequality to see that

$$\int_{[\frac{1}{2}\tau,4\tau]\times\mathbb{R}^2\times S^1} \chi_R^2 x_\tau^2 (|\nabla\eta|^2 + \|[\mathfrak{a},\eta]\|^2) \leq c(\tfrac{1}{\tau^2} + \tfrac{1}{R^2}) \int_{[\frac{1}{2}\tau,4\tau]\times\mathbb{R}^2\times S^1} |\eta|^2 \ .$$

(II.2.16)

Here, $c$ is a number that is independent of $\tau$, R and $\eta$. Taking the $R \to \infty$ limit of both sides of this inequality gives the inequality in the lemma after changing the notation by replacing $\tau$ with t.

### 3. The case of the Nahm pole solution

The Nahm pole solution is the m = 0 version of (II.1.2). The connection in this case is the product connection $\theta_0$ and $\mathfrak{a}$ in this case is such that $\{\mathfrak{a}_i = -\tfrac{1}{2t}\sigma_i\}_{i=1,2,3}$. This section serves as a warm-up to prove that the Nahm pole version of $\mathcal{D}$ is a Fredholm map with trivial kernel and cokernel from the Nahm pole version of $\mathbb{H}$ to the space $\mathbb{L}$ (this is an instance of a more general theorem in [MW1].) For the record, the square of the norm for the Nahm pole version of $\mathbb{H}$ is given by the rule

$$\psi \to \int_{(0,\infty)\times\mathbb{R}^2\times S^1} (|\nabla\psi|^2 + \tfrac{2}{t^2}|\psi|^2) \ .$$

(II.3.1)

By way of a reminder, the norm on $\mathbb{L}$ is the $(0,\infty)\times\mathbb{R}^2\times S^1$ integral of $|\cdot|^2$.

### a) Preliminary observations for the Nahm pole case

Important points are that E = B = 0 (the connection is flat) and that $\nabla_i \mathfrak{a} \equiv 0$ for $i \in \{1, 2, 3\}$. Meanwhile $\nabla_t \mathfrak{a} = -\tfrac{1}{t}\mathfrak{a}$.

These observations have the following implications: Write $\mathcal{D}$ as in (I.1.3) so as to define D and $\{\mathfrak{a}, \cdot\}$. (Using the Clifford algebra notation of (I.1.8), then D is $\gamma_i \nabla_i$ and $\{\mathfrak{a}, \cdot\}$ is $\rho_i[\mathfrak{a}_i, \cdot]$.)

- $[\nabla_t, D] = 0$.
- $D(\{\mathfrak{a},\cdot\}) + \{\mathfrak{a},D(\cdot)\} = 0$.
- $[\nabla_t, \{\mathfrak{a},\cdot\}] = -\tfrac{1}{t}\{\mathfrak{a},\cdot\}$.
- $D^2 = -\nabla_i \nabla_i$

(II.3.2)

There is one other crucial fact that enters the arguments:

*The eigenvalues of the endomorphism $\{\mathfrak{a},\cdot\}$ are $\{\pm\tfrac{1}{t}, \pm 2\tfrac{1}{t}\}$ .*

(II.3.3)

This is left as an exercise to prove. (Here is a hint: Use endomorphisms constructed from $\rho_1, \rho_2, \rho_3$ and $\gamma_1, \gamma_2$ and $\gamma_3$ and $[\sigma, \cdot]$ to find a pairwise commuting set that commutes with $\{\mathfrak{a}, \cdot\}$.)

**b) Closed range and trivial cokernel**

Supposing that $\psi$ is a section from $\mathbb{H}$, then the top three bullets in (II.3.2) can be used to derive the following partial Bochner-Weitzenboch formula (via integration by parts):

$$\int_{(0,\infty)\times\mathbb{R}^2\times S^1} |\mathcal{D}\psi|^2 = \int_{(0,\infty)\times\mathbb{R}^2\times S^1} (|\nabla_t\psi|^2 + |D\psi|^2 + \tfrac{1}{t}\langle\psi,\{\mathfrak{a},\psi\}\rangle + |\{\mathfrak{a},\psi\}|^2) \ .$$

(II.3.4)

Then, by virtue of (II.3.3), this inequality implies in turn that

$$\int_{(0,\infty)\times\mathbb{R}^2\times S^1} |\mathcal{D}\psi|^2 \geq \int_{(0,\infty)\times\mathbb{R}^2\times S^1} (|\nabla_t\psi|^2 + |D\psi|^2)$$

(II.3.5)

whose right hand side can be written using the fourth bullet of (II.3.2) as

$$\int_{(0,\infty)\times\mathbb{R}^2\times S^1} |\nabla\psi|^2$$

(II.3.6)

Thus, $\|\mathcal{D}\psi\|_\mathbb{L} \geq \|\nabla\psi\|_\mathbb{L}$. This implies that $\|\mathcal{D}\psi\|_\mathbb{L} \geq \tfrac{1}{2}\|\tfrac{1}{t}\psi\|_\mathbb{L}$ because of (I.2.4). These last two inequalities and (II.3.1) lead directly to this: $9\|\mathcal{D}\psi\|_\mathbb{L} \geq \|\psi\|_\mathbb{H}$. As observed in Section 1d, an inequality of this sort implies that $\mathcal{D}$ has closed range and trivial kernel.

**c) The cokernel of $\mathcal{D}$**

Because $\mathcal{D}$ has closed range, its cokernel is the kernel $\mathcal{D}^\dagger$ in $\mathbb{L}$. The latter operator is depicted in (I.1.4). Since $\mathcal{D}^\dagger$ is elliptic and its coefficients are smooth, elements in its kernel are smooth also. (Even so, they can't be in $\mathbb{H}$ because the same argument but that leads to (II.3.5) and the bound $9\|\mathcal{D}\psi\|_\mathbb{L} \geq \|\psi\|_\mathbb{H}$ leads to the bound $9\|\mathcal{D}^\dagger\psi\|_\mathbb{L} \geq \|\psi\|_\mathbb{H}$.)

To start the proof that $\mathcal{D}^\dagger$ has no $\mathbb{L}$-kernel, remember that the eigenvalues of $\{\mathfrak{a},\cdot\}$ on any given constant t slice of $(0,\infty)\times\mathbb{R}^2\times S^1$ are $\pm\tfrac{1}{t}$ and $\pm\tfrac{2}{t}$. Because of the second bullet of (II.3.2), the operator D preserves the norm of the eigenvalue but switches the sign. Therefore (and by virtue of the third bullet of (II.3.2)), it is sufficient to study those elements in the kernel of $\mathcal{D}^\dagger$ that have the form $\eta = \eta^+ + \eta^-$ with $\{\mathfrak{a},\eta^+\}$ is equal to either $\frac{1}{t}\eta^+$ or $\frac{2}{t}\eta^+$ (but not both) and with $\{\mathfrak{a}, \eta^-\}$ being the corresponding $-\frac{1}{t}\eta^-$ or $-\frac{2}{t}\eta^-$ as the case may be. Then, for $\eta$ as just described, $\mathcal{D}^\dagger \eta$ can be projected to the relevant $\pm\lambda$ eigenspace of $\{\mathfrak{a}, \cdot\}$ to give the two equations

- $-\nabla_t \eta^+ + \frac{\lambda}{t}\eta^+ + D\eta^- = 0$,
- $-\nabla_t \eta^- - \frac{\lambda}{t}\eta^- + D\eta^+ = 0$.

(II.3.7)

Take the inner product of the top equation with $\eta^+$ and the lower one with $\eta^-$ and subtract the result of doing that to the lower one from the result of doing that to the top one. Then integrate over a given slice $\{t\} \times \mathbb{R}^2 \times S^1$. Integration over all but a measure zero set of slices is finite, and when it is, then integration by parts is allowed (this follows from what is said by Lemma II.2.2). Because D is symmetric on each slice, integration by parts eliminates the terms with D and gives the following identity:

$$\frac{d}{dt} \int_{\{t\}\times\mathbb{R}^2\times S^1} (|\eta^+|^2 - |\eta^-|^2) = \frac{2\lambda}{t} \int_{\{t\}\times\mathbb{R}^2\times S^1} (|\eta^+|^2 + |\eta^-|^2)$$

(II.3.8)

Lemma II.2.2 implies that both integrals that appear in (II.3.8) when viewed as functions on $(0, \infty)$ have $L^1$ derivatives on compact subsets of $(0, \infty)$. This identity implies (among other things) that the left hand integral has a $C^0$ derivative.

Now suppose that there exists some $t_0 > 0$ such that the $\{t_0\} \times \mathbb{R}^2 \times S^1$ integral of $|\eta^+|^2$ is greater than the $\{t_0\} \times \mathbb{R}^2 \times S^1$ integral of $|\eta^-|^2$. Because (II.3.8) leads to this:

$$\frac{d}{dt} \int_{\{t\}\times\mathbb{R}^2\times S^1} (|\eta^+|^2 - |\eta^-|^2) \geq \frac{2\lambda}{t} \int_{\{t\}\times\mathbb{R}^2\times S^1} (|\eta^+|^2 - |\eta^-|^2),$$

(II.3.9)

it follows that

$$\int_{\{t\}\times\mathbb{R}^2\times S^1} (|\eta^+|^2 - |\eta^-|^2) \geq (\tfrac{t}{t_0})^{2\lambda} \int_{\{t_0\}\times\mathbb{R}^2\times S^1} (|\eta^+|^2 - |\eta^-|^2)$$

(II.3.10)

for $t \geq t_0$. Since $\lambda \geq 1$, and since the right hand side is positive (by assumption), this growth is not compatible with $\eta$ being in $\mathbb{L}$.

Suppose on the other hand that there exists some $t_0 > 0$ such that the $\{t_0\} \times \mathbb{R}^2 \times S^1$ integral of $|\eta^-|^2$ is greater than the $\{t_0\} \times \mathbb{R}^2 \times S^1$ integral of $|\eta^+|^2$. Because (II.3.8) leads to the inequality

$$\frac{d}{dt} \int_{\{t\} \times \mathbb{R}^2 \times S^1} (|\eta^-|^2 - |\eta^+|^2) \leq -\frac{2\lambda}{t} \int_{\{t\} \times \mathbb{R}^2 \times S^1} (|\eta^-|^2 - |\eta^+|^2) \,,$$

(II.3.11)

it follows that

$$\int_{\{t\} \times \mathbb{R}^2 \times S^1} (|\eta^-|^2 - |\eta^+|^2) \geq \left(\frac{t_0}{t}\right)^{2\lambda} \int_{\{t_0\} \times \mathbb{R}^2 \times S^1} (|\eta^-|^2 - |\eta^+|^2)$$

(II.3.12)

for $t \leq t_0$ which is likewise incompatible with $\eta$ being in $\mathbb{L}$.

Thus, if $\eta \in \mathbb{L}$, then the $\{t\} \times \mathbb{R}^2 \times S^1$ integrals of $|\eta^+|^2$ and $|\eta^-|^2$ must be equal for all t. But that runs afoul of (II.3.8) unless these integrals are identically zero.

## 4. Symmetries of $\mathcal{D}$

This section describes certain invariance properties of $\mathcal{D}$ which are used to prove that the Witten's model solution versions are Fredholm with trivial kernel and cokernel. The first is $\mathcal{D}$'s invariance with respect to the constant translations of the $S^1$ coordinate (this is the symmetry $x_3 \to x_3$ + constant). The second is $\mathcal{D}$'s homogeniety with respect to the coordinate rescalings of the (t, z) coordinates on the $x_3$-invariant subspace of its domain. (These coordinate rescaling diffeomorphism appear in (II.1.4)). These two properties are used to prove that $\mathcal{D}$ has trivial kernel (see Proposition II.4.1). There are also two important algebraic symmetries that play a role in Section 5's proof that $\mathcal{D}$ has closed range and in Section 6's proof that $\mathcal{D}$ has trivial cokernel.

The last subsection in this section describes another continuous symmetry of $\mathcal{D}$ which is its equivariance with respect to the constant rotations of the $\mathbb{R}^2$ factor of $(0,\infty) \times \mathbb{R}^2 \times S^1$. This $\mathbb{R}^2$-rotation symmetry plays no role in subsquent arguments.

### a) Symmetry with respect $x_3$

The action of constant rotations of the $S^1$ factor in $\mathbb{R}^2 \times S^1$ induces isometric actions of $S^1$ on $\mathbb{H}$ and on $\mathbb{L}$ via pull-back which are intertwined by the operator $\mathcal{D}$. (This is because both $[\nabla_3, \nabla] = 0$ and $\nabla_3 \mathfrak{a} = 0$.) This $S^1$ translation symmetry has the following implications: Distinguish the $x_3$-derivative part of $\mathcal{D}$ by writing $\mathcal{D}$ as $\Xi + \gamma_3 \nabla_3$ with $\Xi$ denoting the operator

$$\Xi = \nabla_t + \gamma_1 \nabla_1 + \gamma_2 \nabla_2 + \rho_3[\mathfrak{a}_3, \cdot\,] + \rho_1[\mathfrak{a}_1, \cdot\,] + \rho_2[\mathfrak{a}_2, \cdot\,].$$

(II.4.1)

The preceding decomposition of $\mathcal{D}$ and the fact that $[\nabla_3, \mathcal{D}] = 0$ leads to a 'partial' Bochner-Weitzenboch formula that has the form

$$\int_{(0,\infty)\times\mathbb{R}^2\times S^1} |\mathcal{D}\psi|^2 = \int_{(0,\infty)\times\mathbb{R}^2\times S^1} |\nabla_3\psi|^2 + \int_{(0,\infty)\times\mathbb{R}^2\times S^1} |\Xi\psi|^2 \ .$$

(II.4.2)

This formula has this immediate implication:

*If $\psi \in \mathrm{kernel}(\mathcal{D})$, then $\nabla_3\psi \equiv 0$ and $\Xi\psi \equiv 0$.*

(II.4.3)

This is to say that the respective kernels of $\mathcal{D}$ and $\Xi$ are identical.

There is a second implication with regards to the image of $\mathcal{D}$. To say more about this, note first that the complexifications of both $\mathbb{H}$ and $\mathbb{L}$ can be written as direct sums of orthogonal subspaces indexed by $\mathbb{Z}$ such that the elements in the index k summand (for $k \in \frac{2\pi}{\ell}\mathbb{Z}$) have the form $e^{ikx_3}\xi$ with $\xi$ being an 8-component vector of $\mathfrak{su}(2)$ valued functions on $(0,\infty) \times \mathbb{R}^2$ (which is to say that it depends only the variables $(t, z_1, z_2)$). Because $\mathcal{D}$ commutes with $\nabla_3$, the operators $\mathcal{D}$ and $\Xi$ map the index k summand in $\mathbb{H}$ to the index k summand in $\mathbb{L}$. It then follows from this and (II.4.2) that $\mathcal{D}$ has closed range mapping $\mathbb{H}$ to $\mathbb{L}$ if and only if the restriction of $\Xi$ to the $k = 0$ summand in $\mathbb{H}$ has closed range as a map to the corresponding $k = 0$ summand in $\mathbb{L}$. (The operator $\Xi$ is determined by its action on the $k = 0$ summand because elements in any $k \neq 0$ have the form $e^{ikx_3}\xi$ with $\xi$ being $S^1$-invariant, and because $\Xi$ acts on these elements to give $e^{ikx_3}(\Xi\xi)$.)

The following is also a consequence of (II.4.2): If $\{\eta_n\}_{n\in\mathbb{N}}$ is a convergent sequence in any given $k \in \frac{2\pi}{\ell}\mathbb{Z}$ indexed summand of $\mathbb{L}$ and such that each $\eta_k$ has the form $\mathcal{D}\psi_k$ for $\psi_k$ in the k-indexed summand of $\mathbb{H}$, then $\{k\psi_k\}_{k\in\mathbb{N}}$ converges in $\mathbb{L}$.

With regards to $\Xi$: Because the operator $\Xi$ acts on $e^{ikx_3}\xi$ to give $e^{ikx_3}(\Xi\xi)$, it can and henceforth will be viewed as an operator taking the $S^1$-invariant subspace in $\mathbb{H}$ to the $S^1$-invariant subspace in $\mathbb{L}$. Also: When restricting to this subspace (which will be done henceforth when discussing $\Xi$ without further notice or new notation), integrations such as those in (I.1.9) (I.1.10), (I.2.3), (I.2.8), (II.2.8) and (II.4.2) can and will be restricted to a fixed $x_3 \in S^1$ slice of $(0,\infty)\times\mathbb{R}^2\times S^1$. This is to say that all integrations unless stated to the contrary will be over the domain $(0,\infty) \times \mathbb{R}^2$.

**b) The coordinate rescaling symmetry**

As remarked above, it is sufficient for proving closed range and finite dimensional cokernel to do that for the operator $\Xi$ on $(0,\infty)\times\mathbb{R}^2$ in (II.4.1). This operator is covariant with respect to the 1-parameter group of coordinate rescaling diffeomorphisms that is defined as follows: The parameter space is the space of positive real numbers; and the diffeomorphism defined by any such number (call it $\lambda$) is the map taking $(t, z) \in (0, \infty) \times \mathbb{R}^2$ to $(\lambda t, \lambda z)$. With this map denoted by $\phi_\lambda$, the covariance of $\Xi$ is this: It obeys $\Xi(\phi_\lambda^*(\psi)) = \lambda^{-1}\phi_\lambda^*(\Xi\psi)$. This can be said in terms of the generator of the action, the operator $t\nabla_t + z_1\nabla_1 + z_2\nabla_2$. Denoting the latter by $x\nabla_x$, then $[x\nabla_x, \Xi] = -\Xi$. Something more is true, which is this: The $\gamma$ and $\rho$ matrices commute with $x\nabla_x$, whereas any of the covariant derivatives operators $\nabla_t, \nabla_1$ and $\nabla_2$ and any of the commutator operators $[\mathfrak{a}_1, \cdot\,], [\mathfrak{a}_2, \cdot\,]$ and $[\mathfrak{a}_3, \cdot\,]$ obey $[x\nabla_x, (\cdot)] = -(\cdot)$. (The preceding remarks about these operators follow directly from the depiction of $(A, \mathfrak{a})$ in (II.1.2).)

This coordinate rescaling property can be exploited by first introducing the fiberwise orthogonal endomorphism $U \equiv \frac{1}{x}(t + z_1\gamma_1 + z_2\gamma_2)$ and then noting that $\Xi\psi$ for any given $\psi$ can be written as

$$\Xi\psi = (\nabla_x + \tfrac{1}{x}\Omega)U\psi$$

(II.4.4)

with $\Omega$ denoting a certain first order, symmetric operator mapping $\mathbb{H}$ to $\mathbb{L}$ that commutes with both $\nabla_x$ and multiplication by $x$. More is said about $\Omega$ later in this subsection.

With regards to U: This endomorphism extends to define a unitary transformation of $\mathbb{L}$ and an invertible, bounded map from $\mathbb{H}$ to $\mathbb{H}$. It is not unitary for the $\|\cdot\|_\mathbb{H}$ norm because it doesn't commute with the covariant derivative. It is none-the-less bounded as a map from $\mathbb{H}$ to $\mathbb{H}$ because $|\nabla U|$ is bounded by a constant multiple of $\frac{1}{x}$ and multiplication by $\frac{1}{x}$ defines a bounded map from $\mathbb{H}$ to $\mathbb{L}$ (see (II.2.2)).

The significance of $\Omega$ comes from the next proposition. To set the stage, let $S^+$ again denote the hemisphere in $(0,\infty)\times\mathbb{R}^2$ where $x = 1$; and letting $\mathbb{V}$ denote $\oplus_8 \mathfrak{su}(2)$, reintroduce the Hilbert space $\mathbb{H}_S$ from Section 2c. Likewise define $\mathbb{L}_S$. The operator $\Xi$ defines a bounded map from $\mathbb{H}_S$ to $\mathbb{L}_S$ because it is a first order operator and because the parts that involve commutators with components of $\mathfrak{a}$ can be bounded using (II.2.5).

**Proposition II.4.1**: *The operator $\Xi$ has closed range as a map from $\mathbb{H}$ to $\mathbb{L}$ if there exists a positive number (call it $\varepsilon$) such that*

$$\int_{S^+} |\Omega\xi - \tfrac{1}{2}\xi|^2 \geq \varepsilon \int_{S^+} |\xi|^2$$

*for all $\xi \in \mathbb{H}_S$. In any event, $\Xi$ has trivial kernel.*

This proposition is proved momentarily.

By way of a parenthetical remark: If there is no ε for the proposition, then Ξ won't have closed range.

The next proposition makes some formal observations about Ω which imply the following: Either Proposition II.4.1 has its ε or there exists $\xi \in \mathbb{H}_S$ such that $\Omega\xi = \frac{1}{2}\xi$.

**Proposition II.4.2**: *The bilinear form on $\mathbb{H}_S$ defined by the pairing*

$$(\zeta,\xi) \to \int_{S^+} \langle \Omega\zeta, \Omega\xi \rangle$$

*has an $\mathbb{L}_S$-orthonormal basis of eigensections whose corresponding set of eigenvalues is a discrete set bounded from below with no accumulation points and finite multiplicities.*

This proposition is proved after the proof of Proposition II.4.1.

***Proof of Proposition II.4.1***: The proof has two parts. The first part shows that if the conditions of the proposition are met, then there exists ε > 0 such that

$$\int_{[0,\infty)\times\mathbb{R}^2} \frac{1}{x^2}|\psi|^2 \leq \frac{1}{\varepsilon} \int_{[0,\infty)\times\mathbb{R}^2} |\Xi\psi|^2$$

(II.4.5)

for all $\psi \in \mathbb{H}$. This part also proves that Ξ has trivial kernel. The second part of the proof uses the preceding inequality to prove that $\|\psi\|_{\mathbb{H}}$ is bounded by a ψ-independent multiple of $\|\Xi\psi\|_{\mathbb{L}}$ when ψ is from $\mathbb{H}$. That implies (directly) that Ξ has closed range.

*Part 1*: Introduce ξ to denote Uψ. Then (II.4.4) leads to the identity

$$\int_{[0,\infty)\times\mathbb{R}^2} \left(|\nabla_x\xi|^2 + \frac{2}{x}\langle \nabla_x\xi, \Omega\xi\rangle + \frac{1}{x^2}|\Omega\xi|^2\right) = \int_{[0,\infty)\times\mathbb{R}^2} |\Xi\psi|^2 \, .$$

(II.4.6)

To see where this leads, write $\frac{2}{x}\langle \nabla_x\xi,\Omega\xi\rangle$ as $\frac{1}{x}\langle \nabla_x\xi,\Omega\xi\rangle + \frac{1}{x}\langle \nabla_x\xi,\Omega\xi\rangle$ and then integrate by parts with respect to the derivatives in Ω on the first summand and integrate by parts with respect to $\nabla_x$ for the second to rewrite the integral in (II.4.6) as

$$\int_{[0,\infty)\times\mathbb{R}^2} \left(|\nabla_x\xi|^2 - \frac{1}{x^2}\langle \xi, \Omega\xi\rangle + \frac{1}{x^2}|\Omega\xi|^2\right) = \int_{[0,\infty)\times\mathbb{R}^2} |\Xi\psi|^2 \, .$$

(II.4.7)

Now write (II.4.7) as

$$\int_{[0,\infty)\times\mathbb{R}^2} (|\nabla_x \xi|^2 - \tfrac{1}{4x^2}|\xi|^2) + \int_{[0,\infty)\times\mathbb{R}^2} \tfrac{1}{x^2}|\Omega\xi - \tfrac{1}{2}\xi|^2 = \int_{[0,\infty)\times\mathbb{R}^2} |\Xi\psi|^2$$

(II.4.8)

This in turn can be written (using another integration by parts) as

$$\int_{[0,\infty)\times\mathbb{R}^2} |\nabla_x \xi + \tfrac{1}{2x}\xi|^2 + \int_{[0,\infty)\times\mathbb{R}^2} \tfrac{1}{x^2}|\Omega\xi - \tfrac{1}{2}\xi|^2 = \int_{[0,\infty)\times\mathbb{R}^2} |\Xi\psi|^2 \; .$$

(II.4.9)

If there exists ε from the proposition, then (II.4.9) leads directly to (II.4.5).

With regard to the kernel of $\Xi$ in $\mathbb{H}$: The identity in (II.4.9) proves that $\Xi$ has trivial kernel because if $\psi$ is from $\mathbb{H}$ and annihilated by $\Xi$, then $U\psi$ must be annihilated by $\nabla_x + \tfrac{1}{2x}$; and that can happen only if $|\psi|$ can be written as $\tfrac{1}{\sqrt{x}}$ times a non-negative, $x$-independent function. But nothing like that (except 0) can come from $\mathbb{H}$ because the square of the norm of the $x$-derivative doesn't have finite integral on $(0,\infty) \times \mathbb{R}^2$.

*Part 2*: Now suppose that (II.4.5) holds. Use $\varsigma$ to denote the function on the domain $(0,\infty) \times \mathbb{R}^2$ given by the rule $\varsigma(t,z) = \chi(\tfrac{|z|}{t} - 1)$ with $\chi$ as described in Section I.2d. This function $\varsigma$ is equal to 1 where $|z|$ is less than t and it is equal to zero where $|z|$ is greater than 2t.

Given $\psi$ from $\mathbb{H}$, write $\psi$ as $\psi_0 + \psi_1$ with $\psi_0 = \varsigma\psi$ and $\psi_1 = (1-\varsigma)\psi$. Having done this, then write

$$|\Xi\psi|^2 = |\Xi\psi_0|^2 + |\Xi\psi_1|^2 + 2\varsigma\langle\Xi\psi,\Xi\psi_1\rangle + 2\langle\mathfrak{S}_\Xi(D\varsigma)\psi,\Xi\psi_1\rangle$$

(II.4.10)

with $D\varsigma$ denoting the $(0,\infty)\times\mathbb{R}^2$ exterior derivative of $\varsigma$ and with $\mathfrak{S}_\Xi(D\varsigma)$ denoting the symbol of the operator $\Xi$. Since $|\mathfrak{S}_\Xi(D\varsigma)|$ is bounded by a constant multiple of $\tfrac{1}{x}$, integration of (II.4.10) with (II.4.5) and the Cauchy-Schwarz inequality leads to the following observation: There exists a $\psi$-independent number (to be denoted by $c$) with the property that

$$\int_{[0,\infty)\times\mathbb{R}^2} |\Xi\psi_0|^2 + \int_{[0,\infty)\times\mathbb{R}^2} |\Xi\psi_1|^2 \le c \int_{[0,\infty)\times\mathbb{R}^2} |\Xi\psi|^2 \; .$$

(II.4.11)

To see where this now leads: First invoke the Bochner-Weitzenboch identity in (I.2.8) for the integral of $|\Xi\psi_0|^2$ (take both $\xi$ and $\psi$ to be $\psi_0$). Because $\varsigma$ has support where $|z| < 2t$, the norm of the endomorphism that is denoted by $\mathbb{X}$ in (I.2.8) is bounded by an $x$-independent multiple of $\tfrac{1}{x^2}$ (because $\tfrac{1}{t^2} \le \tfrac{16}{x^2}$ where $|z| < 2t$). Because of this and because of (II.4.5), the identity in (I.2.8) leads to the inequality

$$\|\psi_0\|_{\mathbb{H}}^2 \leq \int_{[0,\infty)\times\mathbb{R}^2} |\Xi\psi_0|^2 + c_0 \int_{[0,\infty)\times\mathbb{R}^2} |\Xi\psi|^2$$

(II.4.12)

with $c_0$ being independent of $\psi$. Because of (II.4.11) and (II.4.12), $\|\psi_0\|_{\mathbb{H}}$ is bounded by a $\psi$-independent multiple of $\|\Xi\psi\|_{\mathbb{L}}$.

To see about the $|\Xi\psi_1|^2$ integral in (II.4.11), return for the moment to the formula in (II.1.2). The key observation is that there is an automorphism of $P/\{\pm 1\}$ on the $|z| > t$ part of $(0,\infty)\times\mathbb{R}^2$ to be denoted by g such that $gAg^{-1} + gdg^{-1}$ can be written as $A^{NP} + \mathfrak{r}_A$ and $g\mathfrak{a}g^{-1}$ can be written as $\mathfrak{a}^{NP} + \mathfrak{r}_\mathfrak{a}$ where $(A^{NP}, \mathfrak{a}^{NP})$ is the Nahm pole solution and where $\mathfrak{r}_A$ and $\mathfrak{r}_\mathfrak{a}$ have norms bounded by $\frac{t}{x|z|}$. This implies in turn that

$$g\Xi g^{-1} = \Xi_{NP} + \mathfrak{r}$$

(II.4.13)

where $\Xi_{NP}$ is the Nahm pole solution's version of $\Xi$ and where $\mathfrak{r}$ is an endomorphism with norm bounded by a constant multiple of $\frac{t}{x|z|}$ also. Use this decomposition with (II.4.5) and (II.4.11) to see that

$$\int_{[0,\infty)\times\mathbb{R}^2} |\Xi_{NP}(g\psi_1)|^2 \leq c_2 \int_{[0,\infty)\times\mathbb{R}^2} |\Xi\psi|^2 \ .$$

(II.4.14)

Now use (II.4.14) with the conclusions from Section II.3b (and the bounds for $\mathfrak{r}_A$ and $\mathfrak{r}_\mathfrak{a}$, and again (II.4.5)) to see that $\|\psi_1\|_{\mathbb{H}}$ is bounded by a $\psi$-independent multiple of $\|\Xi\psi\|_{\mathbb{L}}$.

The conclusions of the preceding two paragraphs imply the desired conclusion that $\|\psi\|_{\mathbb{H}}$ is bounded by a $\psi$-independent multiple of $\|\Xi\psi\|_{\mathbb{L}}$ because $\|\psi\|_{\mathbb{H}}$ is bounded by a constant multiple of the sum of $\|\psi_0\|_{\mathbb{H}}$ and $\|\psi_1\|_{\mathbb{H}}$ and $\|\frac{1}{x}\psi\|_{\mathbb{L}}$ which have all been bounded at this point by a $\psi$-independent multiple of $\|\Xi\psi\|_{\mathbb{L}}$.

*Proof of Proposition II.4.2*: If the bilinear form on $\mathbb{H}_S$ given by the rule

$$\mathcal{E}(\psi,\eta) = \int_{S^+} (\langle\Omega\psi - \tfrac{1}{2}\psi, \Omega\eta - \tfrac{1}{2}\eta\rangle - \tfrac{1}{4}\langle\psi,\eta\rangle)$$

(II.4.15)

satisfies the conditions for Proposition II.2.1 then so does the bilinear form in the proposition; and if that is the case, then Proposition II.4.2 follows as an instance of Proposition II.2.1. To see about these conditions, compare (II.4.8) with the Bochner-Weitzenboch formula in (I.2.8) to see that

$$\mathcal{E}(\psi,\psi) = \int_{S^+} (|\nabla^S \psi|^2 + \|[\mathfrak{a},\psi]\|^2 + \langle \psi, \mathbb{X}\psi \rangle)$$

(II.4.16)

with $\mathbb{X}$ depicted in (II.2.7). By virtue of (II.4.16), the endomorphism $[\mathfrak{a}_i, [\cdot, \mathfrak{a}_i]] + \mathbb{X}$ plays the role of what is denoted by W in Proposition II.2.1. Therefore, by virtue of (II.1.3), the corresponding version of $t^2|W|$ is bounded as required by Proposition II.2.1. It remains only to check that W is bounded from below on $S^+$. In this regard: The only issue concerns the behavior of W as $t \to 0$ on $S^+$. To study the small t behavior, use the last bullet in (II.1.3) to see that W differs by a bounded endomorphism from the Nahm pole's version where t is small. Since the latter is $\{\mathfrak{a}_{NP}, \{\mathfrak{a}_{NP}, \cdot\}\} + \frac{1}{t}\{\mathfrak{a}_{NP}, \cdot\}$ which is non-negative (see (II.3.3)), the endomorphism W is bounded from below as required.

### d) Algebraic symmetries

A pair of commuting constant endomorphisms of $\oplus_8 \mathfrak{su}(2)$ commute with $\Xi$ also (and with $\nabla_x$ and $\mathfrak{L}$). The first is denoted by Q and it is defined as follows:

$$Q = \rho_1 \rho_2 - [\sigma_3, \cdot] \,.$$

(II.4.17)

The eigenvalues of Q are $\pm 3i$ and $\pm i$. They are imaginary because Q is anti-symmetric. Also because of this: The eigenvectors are sections of the complexification of $\oplus_8 \mathfrak{su}(2)$ (which is $\oplus_8 \mathfrak{sl}(2;\mathbb{C})$); and if an eigenvector has eigenvalue $3i$ or $i$, then the hermitian conjugate of the eigenvector (hermitian conjugate each summand in $\oplus_8 \mathfrak{sl}(2;\mathbb{C})$) has eigenvalue $-3i$ or $-i$ as the case may be.

The second endomorphism is denoted by L. It is defined by first writing any element $\psi \in \oplus_8 \mathfrak{su}(2)$ as $\psi_0 \sigma_3 + \psi_\perp$ where $\psi_0 \in \mathbb{R}^8$ and with $\langle \sigma_3 \psi_\perp \rangle = 0$. Then

$$L\psi = -\rho_1 \rho_2 \gamma_3 (-\psi_0 \sigma_3 + \psi_\perp).$$

(II.4.18)

This can also be written as $-\rho_1 \rho_2 \gamma_3\, e^{i\frac{\pi}{2}\sigma_3} \psi\, e^{-i\frac{\pi}{2}\sigma_3}$. The endomorphism L is symmetric so its eigenvectors are real; and it has square 1 so its eigenvalues are $\pm 1$.

### e) Rotations of the $\mathbb{R}^2$ factor of $(0,\infty) \times \mathbb{R}^2$

Let $(A, \mathfrak{a})$ denote one of Witten's model solutions from (II.1.2). Rigid rotations about the origin in the $\mathbb{R}^2$ factor of $(0,\infty) \times \mathbb{R}^2$ take $(A, \mathfrak{a})$ to a new pair that is equivalent to the original via an automorphism of P (a gauge transformation). As a consequence of this, the operator

$$\mathfrak{L} = z_1\nabla_2 - z_2\nabla_1 + \tfrac{(m+1)}{2}\tfrac{\sinh(\Theta)}{\cosh(\Theta)}\tfrac{\cosh((m+1)\Theta)}{\sinh((m+1)\Theta)})[\sigma_3,\cdot] - \tfrac{1}{2}\gamma_1\gamma_2 - \tfrac{1}{2}\rho_1\rho_2$$
(II.4.19)

commutes with $\Xi$. It also commutes with the operator $\nabla_x$ and multiplication by $x$. The eigenvalues of $\mathfrak{L}$ take values $i\mathbb{Z}$. (They are imaginary because $\mathfrak{L}$ is anti-symmetric.)

This symmetry plays no explicit role in subsequent arguments.

## 5. The range of $\Xi$

This section proves that $\Xi$ has closed range which implies in turn that $\mathcal{D}$ has closed range. Here is how the proof will go: By virtue of Propositions II.4.1 and II.4.2, it is sufficient to prove that $\tfrac{1}{2}$ is not an eigenvalue of $\Omega$ on $\mathbb{H}_S$. Now in general, if $\lambda$ is an eigenvalue of $\Omega$ on $\mathbb{H}_S$ and if $\xi$ is the corresponding eigenvector, then $x^{-\lambda}U^{-1}\xi$ is annihilated by $\Xi$. (Remember that $U = \tfrac{1}{x}(t+z_1\gamma_1+z_2\gamma_2)$.) This follows from (II.4.4). Therefore, to show that any given $\lambda \in \mathbb{R}$ is not an eigenvalue of $\Omega$, it is sufficient to show that $\Xi$ does not annihilate any map from $(0,\infty)\times\mathbb{R}^2$ to $\oplus_8\mathfrak{su}(2)$ that has the form $x^{-\lambda}\eta$ with $\eta$ in $\mathbb{H}_S$. This is what will be done for the case $\lambda = \tfrac{1}{2}$ and for all $\lambda \in [0, \tfrac{3}{2}]$. The following lemma makes a formal statement to this effect.

**Lemma II.5.1**: *The operator $\Omega$ has no eigenvectors from $\mathbb{H}_S$ with eigenvalue in the closed interval $[0, \tfrac{3}{2}]$.*

The proof of this lemma is contained in the subsequent subsections.

### a) The role of $\mathfrak{b}_3$ and $\mathfrak{c}_t$

The zero'th order operator $\mathbb{X}$ depicted in (II.2.7) that appears in the Bochner-Weitzenboch formula (I.2.5)) has zero's in the third row and last row, and in the third column and last column. As a consequence, it doesn't see the components $\mathfrak{b}_3$ and $\mathfrak{c}_t$ in (I.1.5) which are the respective third and eighth entries of the column vector. Therefore, if $\psi$ is a map from $(0,\infty)\times\mathbb{R}^2$ to $\oplus_8\mathfrak{su}(2)$ with $\Xi\psi = 0$ and if $\mathfrak{t}$ denotes either of the components $\mathfrak{b}_3$ or $\mathfrak{c}_t$, then

$$-(\nabla_t^2+\nabla_1^2+\nabla_2^2)\mathfrak{t} + [\mathfrak{a}_i,[\mathfrak{t},\mathfrak{a}_i]] = 0 .$$
(II.5.1)

Suppose in addition that $\psi$ has the form $x^{-\lambda}\eta$ with $\eta$ from $\mathbb{H}_S$. Then $\mathfrak{t}$ will have the form $x^{-\lambda}\mathfrak{f}$ with $\mathfrak{f}$ being an $\mathfrak{su}(2)$-valued function on $S^+$ which is in the completion of the space of smooth, compactly supported $\mathfrak{su}(2)$-valued functions on $S^+$ using the norm whos square is the $S^+$ integral of $|\nabla^S(\cdot)|^2$. With this understood, take the inner product of both sides of (II.5.1) with $\mathfrak{f}$ and then integrate over $S^+$ to see that

$$-\lambda(\lambda-1)\int_{S^+}|\mathfrak{f}|^2 + \int_{S^+}(|\nabla^S\mathfrak{f}|^2 + \tfrac{\cosh^2\Theta}{\sinh^2\Theta}\|[\mathfrak{ta},\mathfrak{f}]\|^2) = 0 .$$

(II.5.2)

(Note that the term $\tfrac{\cosh^2\Theta}{\sinh^2\Theta}\|[\mathfrak{ta},\mathfrak{f}]\|^2$ in the right most integral is $x^2\|[\mathfrak{a},\mathfrak{f}]\|^2$.)

Now invoke (II.2.4) to see that this leads to an inequality that forces $\mathfrak{f}$ to be zero unless $\lambda(\lambda-1) > 2$, which is to say that either $\lambda > 2$ or $\lambda < -1$. This implies in particular that if $\lambda$ is in the interval $[-1, 2]$, then both $\mathfrak{b}_3$ and $\mathfrak{c}_t$ must vanish if $\psi$ has the form $x^{-\lambda}\eta$ with $\eta$ from $\mathbb{H}_S$. An instance of this is the case of interest which is when $\lambda = \tfrac{1}{2}$.

Supposing now that $\psi$ is a smooth map to $\otimes_8 \mathfrak{su}(2)$ with $\mathfrak{b}_3$ and $\mathfrak{c}_t$ being zero, then the action of $\Xi$ on $\psi$ will be written in terms of $\mathfrak{su}(2)\otimes_\mathbb{R}\mathbb{C}$ ($\equiv \mathfrak{sl}(2;\mathbb{C})$) valued functions

$$\mathfrak{b} = \tfrac{1}{2}(\mathfrak{b}_1+i\mathfrak{b}_2) \text{ and } \beta = \tfrac{1}{2}(\mathfrak{c}_3+i\mathfrak{b}_t) \text{ and } \mathfrak{o} = \tfrac{1}{2}(\mathfrak{c}_1-i\mathfrak{c}_2) .$$

(II.5.3)

The action of $\Xi$ on $\psi$ when written in terms of these $\mathfrak{sl}(2;\mathbb{C})$ valued functions is an $\mathfrak{sl}(2;\mathbb{C})$-valued vector with four components:

- $\nabla_t \mathfrak{b} + i(\nabla_1+i\nabla_2)\beta - i\alpha[\sigma_3,\mathfrak{b}]$
- $\nabla_t \beta + i(\nabla_1-i\nabla_2)\mathfrak{b} + i\alpha[\sigma_3,\beta] + i[\varphi^*,\mathfrak{o}]$
- $\nabla_t \mathfrak{o} - i\alpha[\sigma_3,\mathfrak{o}] + i[\varphi,\beta]$.
- $i(\nabla_1+i\nabla_2)\mathfrak{o} + [\varphi,\mathfrak{b}]$ .

(II.5.4)

With regards to notation: What is denoted here by $\alpha$ is the $\mathbb{R}$-valued function that is defined by writing $\mathfrak{a}_3$ from (II.1.2) as $\alpha\sigma_3$. Meanwhile, the $\mathfrak{sl}(2;\mathbb{C})$ valued function $\varphi$ is $\mathfrak{a}_1-i\mathfrak{a}_2$ with $\mathfrak{a}_1$ and $\mathfrak{a}_2$ from (II.1.2), and $\varphi^*$ is $\mathfrak{a}_1+i\mathfrak{a}_2$ (this is -1 times the Hermitian conjugate of $\varphi$).

By way of an explanation for (II.5.4): The top bullet in (II.5.4) is obtained from (I.1.2) by taking $\mathfrak{p}_1+i\mathfrak{p}_2$, and the second bullet in (II.5.4) is obtained from (I.1.2) by taking $\mathfrak{q}_3+i\mathfrak{p}_t$. The third and fourth bullets in (II.5.4) are obtained from (I.1.2) by respectively taking $\mathfrak{q}_1-i\mathfrak{q}_2$ and $\mathfrak{q}_t-i\mathfrak{p}_3$.

### b) The joint eigenspaces of Q and L

Assume now that $\psi$ is a smooth map from $(0,\infty)\times\mathbb{R}^2$ to $\otimes_8 \mathfrak{sl}(2;\mathbb{C})$ with $\mathfrak{b}_3$ and $\mathfrak{c}_t$ components being zero. Since the endomorphisms Q and L commute with $\Xi$, the map $\psi$ can be written as a sum of eight elements from $\mathbb{H}_S\otimes_\mathbb{R}\mathbb{C}$ with distinct terms being pointwise orthogonal (with respect to the Hermitian metric on $\oplus_8 \mathfrak{sl}(2;\mathbb{C})$) and in distinct (Q,L) eigenspaces. Moreover if $\Xi$ annihilates $\psi$, then it annihilates each of these terms in the (Q,L) eigenspace decomposition. And if $\psi$ can be written as $x^{-\lambda}\eta$ with $\eta$ coming from $\mathbb{H}_s$ (or just being independent $x$), then each term can be written this way also. With all of this clear, what follows momentarily is a list of four of the eight possible cases for the joint (Q, L) projection of $\psi$ and the effect of $\Xi$ on each. The remaining four cases are obtained from the ones listed below by taking the Hermitian conjugate of each entry of the listed ones. The effect of $\Xi$ on each of these Hermitian conjugated elements is likewise obtained from what is written below via Hermitian conjugation.)

The upcoming list of four cases refers to the decomposition of $\mathfrak{sl}(2;\mathbb{C})$ into eigenspaces of the endomorphism $[\frac{i}{2}\sigma_3, \cdot]$. This decomposition is written as

$$\mathfrak{sl}(2;\mathbb{C}) = L^+ \oplus \mathbb{C}\sigma_3 \otimes L^-$$

(II.5.5)

with $L^+$ and $L^-$ denoting the respective +1 and -1 eigenspaces of $[\frac{i}{2}\sigma_3, \cdot]$. In this regard: The commutator of two elements in $L^+$ is zero and likewise that of two elements in $L^-$. Meanwhile, the commutator of an element in $L^+$ with one in $L^-$ is in the span of $\sigma_3$. By the same token, the trace pairing $(\mathfrak{u}, \mathfrak{v}) \to \langle \mathfrak{u}\mathfrak{v}\rangle = -\frac{1}{2}\mathrm{trace}(\mathfrak{u}\mathfrak{v})$ is zero when either both $\mathfrak{u}$ and $\mathfrak{v}$ are in $L^+$ or in $L^-$; and it is a perfect pairing between $L^+$ and $L^-$. By the same token, the action of Hermitian conjugation is a $\mathbb{C}$-antilinear identification between $L^+$ and $L^-$.

Supposing that $\mathfrak{t}$ denotes an $\mathfrak{sl}(2;\mathbb{C})$ valued function on $(0,\infty)\times\mathbb{R}^2$, then it will be written using the decomposition of (5.5) as $\mathfrak{t}_+ + \mathfrak{t}_0\sigma_3 + \mathfrak{t}_-$. By way of an example: The $\mathfrak{sl}(2;\mathbb{C})$ valued function $\varphi$ that appears in (5.4) and in (1.12) is an $L^+$ valued function whereas $\varphi^*$ in (5.4) an $L^-$-valued function

CASE 1: The $\mathfrak{sl}(2;\mathbb{C})$-valued functions $\mathfrak{b}$ and $\beta$ are in $L^+$ valued (and denoted by $\mathfrak{b}_+$ and $\beta_+$) and $o$ is zero. The lower two bullets in (II.5.4) vanish and the first two are

- $\nabla_t \mathfrak{b}_+ - 2\alpha \mathfrak{b}_+ + i(\nabla_1 + i\nabla_2)\beta_+$
- $\nabla_t \beta_+ + 2\alpha\beta_+ + i(\nabla_1 - i\nabla_2)\mathfrak{b}_+$

(II.5.6)

CASE 2: The $\mathfrak{sl}(2;\mathbb{C})$-valued functions $\mathfrak{b}$ and $\beta$ are proportional to $\sigma_3$ (they are $\mathfrak{b}_0\sigma_3$ and $\beta_0\sigma_3$) whereas $o$ is $L^+$-valued (it is written as $o_+$). The equations in (II.5.4) are

- $\frac{\partial}{\partial t}\mathfrak{b}_0 + i(\frac{\partial}{\partial z_1} + i\frac{\partial}{\partial z_2})\beta_0$
- $\frac{\partial}{\partial t}\beta_0 + i(\frac{\partial}{\partial z_1} - i\frac{\partial}{\partial z_2})\mathfrak{b}_0 - 2\langle\varphi^* o_+\rangle$

- $\nabla_t o_+ - 2\alpha o_+ - 2\varphi \beta_0$.
- $i(\nabla_1 + i\nabla_2) o_+ + 2\varphi b_o$ .

(II.5.7)

CASE 3: The $\mathfrak{sl}(2;\mathbb{C})$-valued functions $b$ and $\beta$ are $L^-$-valued (they are $b_-$ and $\beta_-$) whereas $o$ is proportional to $\sigma_3$ (it is written as $o_0\sigma_3$). The equations in (II.5.4) are

- $\nabla_t b_- + 2\alpha b_- + i(\nabla_1 + i\nabla_2)\beta_-$
- $\nabla_t \beta_- - 2\alpha\beta_- + i(\nabla_1 - i\nabla_2) b_- + 2o_0 \varphi^*$
- $\frac{\partial}{\partial t} o_0 + 2\langle \varphi \beta_- \rangle$.
- $i(\frac{\partial}{\partial z_1} + i\frac{\partial}{\partial z_2}) o_0 - 2\langle \varphi b_- \rangle$ .

(II.5.8)

CASE 4: Both $b$ and $\beta$ are zero and $o$ is $L^-$-valued (it is written as $o_-$). The top two equations in (II.5.4) are zero and the lower two are

- $\nabla_t o_- + 2\alpha o_-$.
- $(\nabla_1 + i\nabla_2) o_-$ .

(II.5.9)

The remaining subsections in this Section II.5 study these four cases with the extra assumption that the original $\psi$ is annihilated by $\Xi$ and that it has the form $\psi = x^{-\lambda}\eta$ with $\lambda$ being constant and with $\eta$ being from $\mathbb{H}_S$. This is to say that expressions in (II.5.6)-(II.5.9) are assumed to vanish and that each of the $\mathfrak{sl}(2;\mathbb{C})$ valued functions $b_{(\cdot)}, \beta_{(\cdot)}$ and $o_{(\cdot)}$ that appear in these expressions is the product of $x^{-\lambda}$ with an $\mathfrak{sl}(2;\mathbb{C})$ valued function from the Hilbert space completion of the space of compactly supported $\mathfrak{sl}(2;\mathbb{C})$-valued functions on $S^+$ using the norm whose square is the $S^+$ integral of $|\nabla^S(\cdot)|^2$. (The latter Hilbert space is denoted by $H_S$; thus $\mathbb{H}_S = \oplus_8 H_S$). The four cases listed above are considered in the order 4, 1, 2, 3 which is the order of complexity.

c) CASE 4

CASE 4 is by far the easiest to deal with: If the expression in the top bullet of (II.5.9) is zero, then $\langle \varphi o_- \rangle$ is independent of t because $\varphi$ obeys $\nabla_t \varphi = 2\alpha \varphi$ (see the top bullet in (II.1.5). Likewise, if the expression in the second bullet of (II.5.9) is zero, then $\langle \varphi o_- \rangle$ is a holomorphic function on each constant t slice of $(0,\infty) \times \mathbb{R}^2$ because the top bullet in (II.1.5) also says that $(\nabla_1 + i\nabla_2)\varphi = 0$. Thus, $\langle \varphi o_- \rangle$ is a constant in t, holomorphic function (call it $h(z)$). This implies in turn that $|o_-|$ can be written as $|h(z)||\varphi|^{-1}$. Since $\varphi$ vanishes at $z = 0$ with order $m$, it follows that $h(z)$ must vanish to order at least $m$ there also. Meanwhile: The complex coordinate z can be written as $x \frac{1}{\cosh\Theta} e^{i\phi}$ with $\phi$ being the longitudinal coordinate on $S^+$, and t can be written as $x \frac{\sinh\Theta}{\cosh\Theta}$. As for $|\phi|^{-1}$, it can be written as $x \frac{\sinh((m+1)\Theta)}{(m+1)\cosh(\Theta)}$. Thus, $|o_-|$ has the form

$$|o_-| = x |h(x \tfrac{1}{\cosh\Theta} e^{i\phi})| \tfrac{\sinh((m+1)\Theta)}{(m+1)\cosh(\Theta)}.$$

(II.5.10)

Under the given assumptions, the function $|o_-|$ is a homogeneous function of $x$. This being the case, then $h$ must be homogeneous polynomial of some degree (call it $p$) no less than $m$; and then (II.5.10) implies that $|o_-| \sim x^{p+1}$. Thus, the eigenvalue $\lambda$ of $\Omega$ must have the form $\lambda = -(p+1)$ which is no greater than $-(m+1)$.

### d) CASE 1

The assumption here is that the expressions in (II.5.6) vanish on $(0, \infty) \times \mathbb{R}^2$. To see the ramifications: Take the (Hermitian) inner product of the top bullet in (II.5.6) with $b_+$ and the lower bullet with $\beta_+$. Then subtract the result of doing this to the top bullet from the result of doing this to the lower bullet. These actions leads to the identity

$$\tfrac{\partial}{\partial t}(|\beta_+|^2 - |b_-|^2) + i(\tfrac{\partial}{\partial x_1} - i\tfrac{\partial}{\partial x_2})\langle \beta_+^* b_+\rangle - i(\tfrac{\partial}{\partial x_1} + i\tfrac{\partial}{\partial x_2})\langle b_+^* \beta_+\rangle = -4\alpha(|\beta_+|^2 + |b_+|^2)$$

(II.5.11)

where $t^*$ for $t \in \mathfrak{sl}(2;\mathbb{C})$ is shorthand for $-t^\dagger$. (The $(\cdot)^*$ notation is introduced here to avoid possible cognitive dissonance from the fact the Hermitian inner product on $\mathfrak{sl}(2;\mathbb{C})$ that comes from the inner product on $\mathfrak{su}(2)$ is $-\langle(\cdot)^\dagger(\cdot)\rangle$ which has the dissonance provoking minus sign out front.)

Let $v$ denote the 1-form

$$v = (|\beta_+|^2 - |b_-|^2)dt + i(\langle \beta_+^* b_+\rangle - \langle b_+^* \beta_+\rangle)dz_1 + (\langle \beta_+^* b_+\rangle + \langle b_+^* \beta_+\rangle)dz_2$$

(II.5.12)

The identity in (II.5.11) asserts that $d*v = -4\alpha(|\beta_+|^2 + |b_+|^2)$. Because $|v|$ is $|\beta_+|^2 + |b_+|^2$, this is saying that

$$d*v = -4\alpha |v|.$$

(II.5.13)

Now fix $x > 0$ and integrate both sides of (II.5.13) over the $x = x$ hemisphere. Having done this, then integrate by parts on the left and use the formula in (I.1.2) for $\alpha$ (it is $-|a_3|$) in the right hand integral. The result is an identity that says

$$\tfrac{d}{dx} \int_{\{x=x\}} *\nu = \tfrac{2}{x} \int_{\{x=x\}} \frac{(m+1)\cosh((m+1)\Theta)}{\sinh((m+1)\Theta)} |\nu| \;.$$

(II.5.14)

(These $S^+$ integrations and the integration by parts are allowed if it is assumed that $\beta_+$ and $\flat_+$ on any given $x = x$ hemisphere are in the Hilbert space $H_S$. See Section II.2a. This is the case here.)

If it is assumed that $\beta_+$ and $\flat_+$ each have the form $x^{-\lambda} t$ with $t$ being from $H_S$, then (II.5.14) says this:

$$(1-\lambda) \int_{S^+} *\nu = \int_{S^+} \frac{(m+1)\cosh((m+1)\Theta)}{\sinh((m+1)\Theta)} |\nu| \;.$$

(II.5.15)

The latter identity can hold only if $|1 - \lambda| > 1$ which is to say that $\lambda$ is either negative or it is greater than 2. (This is because the function multiplying $|\nu|$ in the right hand integral is greater than 1 and $|\nu|$ is no less than the norm of its $dx$ component.)

**e) CASE 2**

The assumption in this case is that the expressions that appear in (II.5.7) are zero. To start the analysis for this case, write $o_+$ where $z \neq 0$ as $\hat{o} \varphi$ with $\hat{o}$ being a $\mathbb{C}$-valued function. Since $\nabla_t \varphi = 2\alpha\varphi$ and $(\nabla_1 + i\nabla_2)\varphi = 0$, the second and third bullets of (II.5.7) when set equal to zero and written using $\hat{o}$ say that

$$\tfrac{\partial}{\partial t} \hat{o} = 2\beta_0 \quad and \quad -i(\tfrac{\partial}{\partial z_1} + i\tfrac{\partial}{\partial z_2})\hat{o} = 2\flat_0$$

(II.5.16)

The top bullet in (II.5.7) is necessarily zero if $\hat{o}$ obeys (II.5.6) and the lower bullet says

$$-(\tfrac{\partial^2}{\partial t^2} + \tfrac{\partial^2}{\partial z_1^2} + \tfrac{\partial^2}{\partial z_2^2})\hat{o} + 4|\varphi|^2 \hat{o} = 0 \;.$$

(II.5.17)

With regards to $\hat{o}$ near $z = 0$: The $\mathfrak{sl}(2;\mathbb{C})$-valued function $\varphi$ near $z = 0$ vanishes as $z^m$ and so in principle, $\hat{o}$ near $z = 0$ can have a meromorphic pole which looks like $\tfrac{1}{z^p}$ to leading order with p being positive but not greater than $m$. Any such pole is compatible with the right hand identity in (II.5.16). Compatibility with the left hand identity in (II.5.16) requires that the $z \to 0$ limit of $z^p \hat{o}$ at fixed t be independent of t as $z \to 0$ (and likewise for any lower order poles of $\hat{o}$). As explained in the next paragraph, poles can be ruled out in certain instances if $o_+$ is a homogeneous of $x$.

If it is assumed that $o_+$ has the form $x^{-\lambda} u$ with $u$ being an $\mathfrak{sl}(2;\mathbb{C})$-valued function on $S^+$, then $\hat{o}$ has the form $x^{1-\lambda} \hat{u}$ with $\hat{u}$ being an $\mathfrak{sl}(2;\mathbb{C})$ valued function that is defined on the $|z| > 0$ part of $S^+$ (this is because φ can be written as $\frac{1}{x}$ times an $\mathfrak{sl}(2;\mathbb{C})$ valued function on $S^+$). As a consequence, if $z^p \hat{o}$ is non-zero at $z = 0$, then it is a constant multiple of $t^{p+1-\lambda}$. Note in particular that this is not compatible with the left hand equation in (II.5.16) unless λ is equal to $p+1$. Since $p$ is a positive integer not greater than $m$, there can be no pole in ô unless λ is from the set $\{2, \ldots, m+1\}$.

Supposing that ô has no pole at $z = 0$, then the $\mathfrak{sl}(2;\mathbb{C})$-valued function $\hat{u}$ (which comes from ô via $\hat{o} = x^{1-\lambda}\hat{u}$) comes from $H_S$. (There are no issues with regards to $\hat{u}$ where $t \to 0$ on $S^+$ because |φ| diverges in this limit.) Granted that $\hat{u}$ is in $H_S$, then taking the inner product of both sides of (II.5.17) with $\hat{u}$ and integrating over $S^+$ implies after an integration by parts that

$$-\lambda(\lambda-1) \int_{S^+} |\hat{u}|^2 + \int_{S^+} (|\nabla^S \hat{u}|^2 + 2\frac{(m+1)^2 \cosh^2 \Theta}{\sinh^2((m+1)\Theta)}|\hat{u}|^2) = 0.$$

(II.5.18)

What with (II.2.4), this forces $\hat{u}$ to vanish unless $\lambda^2 - \lambda > 2$, thus unless $\lambda > 2$ or $\lambda < -1$.

**f) CASE 3**

The assumption here is that the expressions in (II.5.8) are zero. To analyze the implications of this, introduce by way of notation $q$ to denote $\frac{\varphi^*}{|\varphi|^2} o_0$ which is a function on the complement of the $z = 0$ locus in $(0,\infty) \times \mathbb{R}^2$ with values in the $L^-$ summand of the decomposition of $\mathfrak{sl}(2;\mathbb{C})$ that is depicted by (II.5.5). The third and fourth bullets of (II.5.8) say in effect that

- $\beta_- = (\nabla_t + 2\alpha) q$.
- $\flat_- = -i(\nabla_1 + i\nabla_2) q$.

(II.5.19)

The assertion that the top bullet expression in (II.5.8) vanishes is redundant given (II.5.19) because the operators $(\nabla_t + 2\alpha)$ and $(\nabla_1 + i\nabla_2)$ commute when acting on an $L^-$-valued function. By way of a contrast, the second bullet in (II.5.9) says that:

$$-(\nabla_t^2 + \nabla_1^2 + \nabla_2^2) q + (4\alpha^2 + 2|\varphi|^2) q = 0.$$

(II.5.20)

(This identity is derived using the expression for $B_3$ in the third bullet of (II.1.5).)

Now suppose that $\beta_-$, $\flat_-$ and $o_0$ are homogeneous with respect to their $x$-dependence, that each can be written as $x^{-\lambda} u$ with $u$ being an $L^-$-valued function on $S^+$ (in the case of $\beta_-$ or $\flat_-$) or a $\mathbb{C}$-valued function on $S^+$ (in the case of $o_0$). It then follows that $q$ has the form $x^{-\lambda-1} \mathfrak{q}$ with $\mathfrak{q}$ being an $L^-$-valued function that is defined on the complement of the $z = 0$ point in $S^+$. Suppose for the moment that $q$ is bounded as $z \to 0$. Then $q$ will be in $H_S$ because $o_0$ is and because $\frac{1}{|\varphi|}$ is $\mathcal{O}(t)$ as $t \to 0$ in $S^+$. Granted this, take the Hermitian inner product of both sides of (II.5.20) with $q$, integrate the result over $S^+$ and then integrate by parts to obtain the following identity:

$$-\lambda(\lambda+1) \int_{S^+} |q|^2 + \int_{S^+} (|\nabla^S q|^2 + \frac{(m+1)^2(\cosh^2((m+1)\Theta) + \cosh^2\Theta)}{\sinh^2((m+1)\Theta)}|q|^2) = 0.$$

(II.5.21)

What with (II.2.4) and Item a) of the second bullet in (II.1.3), this last identity can hold only in the event that $\lambda^2 + \lambda > 2 + (m+1)^2$. This is to say that $\lambda$ must be outside the interval between $-\frac{1}{2}((9+4(m+1)^2)^{1/2}+1)+$ and $\frac{1}{2}((9+4(m+1)^2)^{1/2}-1)$. In particular, no matter the value of $m$ (assuming $m > 0$), the number $\lambda$ must be either less than $-\frac{7}{2}$ or greater than $\frac{3}{2}$.

The preceding analysis is valid provided that $q \equiv \frac{\varphi^*}{|\varphi|^2} o_0$ is bounded as $z \to 0$ at any fixed, positive value of $t$. As explained next, this is the case when $o_0$ is $x$-dependent, but homogeneous with respect this dependence. This is to say that it has the form $x^{-\lambda}\hat{o}$ with $\lambda \neq 0$ and with $\hat{o}$ being a $\mathbb{C}$-valued function on $S^+$. Assuming this form for $o_0$, write it near $z = 0$ as the product of $x^{-\lambda}$ times a power series in the variables $\frac{z}{t}$ and $\frac{\bar{z}}{t}$. Now use this power series in the third bullet of (II.5.8) near $z = 0$ so see that $\frac{1}{z^m} o_0$ must have a well defined $z \to 0$ limit since it is assumed that $\lambda$ isn't zero. This constant time $|z|^m$ bound for $|o_0|$ near $z = 0$ implies directly that $q$ is bounded as $z \to 0$ because $|\varphi|$ is bounded from below by a non-zero multiple of $|z|^m$ as $z \to 0$.

## 6. The cokernel of $\mathcal{D}$

Because $\mathcal{D}$ has closed range, a proof that the formal $L^2$ adjoint $\mathcal{D}^\dagger$ has trivial kernel in $\mathbb{L}$ completes the proof that $\mathcal{D}$ is Fredholm with trivial kernel and cokernel. With regards to $\mathcal{D}^\dagger$: The latter operator can be written as $\Xi^\dagger + \gamma_3 \nabla_3$ with $\Xi^\dagger$ denoting the formal $L^2$ adjoint of the operator $\Xi$. Moreover, any element in the $\mathbb{L}$-kernel of $\mathcal{D}^\dagger$ can be written as a Fourier series with respect to the $x_3$ coordinate with each Fourier component being in the $\mathbb{L}$ kernel of $\mathcal{D}^\dagger$ also. In particular, if $k \in \frac{2\pi}{\ell}\mathbb{Z}$ and if $e^{ikx_3}\xi$ is the k'th Fourier mode of a kernel element in $\mathbb{L}$ (so $\xi$ is independent of $x_3$ and also in $\mathbb{L}$), then $\xi$ must obey the equation

$$\Xi^\dagger \xi + i k \gamma_3 \xi = 0.$$

(II.6.1)

Section II.6a explains why the k = 0 kernel of $\mathcal{D}^\dagger$ in $\mathbb{L}$ is trivial. Section II.6b uses this fact to prove that there are no k ≠ 0 kernel elements in the $\mathbb{L}$ kernel of $\mathcal{D}^\dagger$.

### a) The $x_3$-invariant kernel

A k = 0 element in the kernel of $\mathcal{D}^\dagger$ from $\mathbb{L}$ is an element from $\mathbb{L}$ in the kernel of $\Xi^\dagger$. Supposing that η is of that sort, write $\Xi$ as in (II.4.4) to see that η obeys the equation

$$(-\nabla_x + \tfrac{1}{x}(\Omega - 2))\eta = 0 .$$

(II.6.2)

Meanwhile, by virute of η being from $\mathbb{L}$, its restriction to all but a measure zero set of constant $x$ hemispheres in $(0,\infty) \times \mathbb{R}^2$ has finite $L^2$ norm (the square of the $L^2$ norm is depicted in (II.2.10).) This implies that η can be written as a generalized Fourier sum with the terms being pairwise orthogonal on the constant $x$ slices and with each term being the product of a function of $x$ and an eigensection of $\Omega$. Each such term must also be annihilated $\Xi^\dagger$ which implies this: If ξ now denotes one of these terms in the sum, and if λ is the relevant $\Omega$-eigenvalue, then ξ must be a homogeneous function of $x$ of the form $x^{\lambda-2} u$ with $u$ being the relevant $\Omega$ eigensection on $S^+$. But this sort of section can't be in $\mathbb{L}$ because it's integral on $(0, \infty) \times \mathbb{R}^2$ is not finite (the integral will diverge as $x \to 0$ if λ is -1 or smaller; and it will diverge as $x \to \infty$ if λ is -1 or larger.)

### b) The $x_3$-dependent kernel

The k ≠ 0 case of (II.6.1) is complicated by the fact that $\gamma_3$ does not commute with $\Xi$. What happens instead is this:

$$\gamma_3 \Xi^\dagger = -\Xi \gamma_3 .$$

(II.6.3)

To see how to deal with this, use (II.4.4) to write (II.6.1) as a pair of equations:

- $\Xi^\dagger \xi + k\zeta = 0$ ,
- $\Xi \zeta - k\xi = 0$

(II.6.4)

where $\zeta \equiv i\gamma_3 \xi$. Now invoke (II.4.4) with the fact that $\nabla_x^\dagger = -\nabla_x - \tfrac{2}{x}$ to rewrite these as

- $(-\nabla_x - \tfrac{2}{x} + \tfrac{1}{x}\Omega)\xi + k\varsigma = 0$ .
- $(\nabla_x + \tfrac{1}{x}\Omega)\varsigma - k\xi = 0$ ,

(II.6.5)

where $\varsigma \equiv i\, U\gamma_3 \xi$.

Here is the reason for introducing (II.6.5): The pair (ξ, ς) can be written at any fixed, positive $x$ as a linear combination of the form

$$(\xi, \varsigma) = \sum_{\eta \in \Lambda} (a_\eta, b_\eta) \eta$$

(II.6.6)

where $\Lambda$ denotes an orthonormal basis of $H_S$ eigensections for the operator $\Omega$ and where the corresponding $(a_\eta, b_\eta)$ is a pair of functions on $(0, \infty)$. This decomposition is useful because $(\xi,\varsigma)$ obeys (II.6.5) if and only if each term in (II.6.6) does also; and the equation in (II.6.5) for any such term amounts to an ordinary differential equation for the pair $(a_\eta, b_\eta)$. Indeed, if $(a, b)$ are the $x$-dependent coefficients of a given eigensection of $\Omega$, and if $\lambda$ is the associated eigenvalue $\lambda$, then the corresponding version of (II.6.5) says that

- $-\frac{d}{dx} a + \frac{\lambda - 2}{x} a + kb = 0$ .
- $\frac{d}{dx} b + \frac{\lambda}{x} b - ka = 0$ .

(II.6.6)

The preceding system of ordinary differential equations has two linearly independent solutions for any value of $\lambda$. But even so, no solution will have square integrable norm on $(0, \infty)$ with respect to the measure $x^2 dx$ unless $\frac{1}{2} < \lambda < \frac{3}{2}$. (This claim is proved in the next subsection.)

This constraint on $\lambda$ precludes the existence of a non-trivial kernel to $\mathcal{D}^\dagger$ in $\mathbb{L}$ because of what is said by Lemma II.5.1.

**c) The constraint on $\lambda$**

What follows is a proof of the assertion that $\frac{1}{2} < \lambda < \frac{3}{2}$ is required if (II.6.6) has a solution with $|a|^2 + |b|^2$ being $x^2 dx$–integrable on $(0, \infty)$. Take the Hermitian inner product of the top equation with $x^3 a$ and that of the lower equation with $x^3 b$. Add the two resulting identities and their complex conjugates to see that

$$\tfrac{1}{2} x^3 \tfrac{d}{dx}(|b|^2 - |a|^2) + x^2((\lambda - 2)|a|^2 + \lambda|b|^2) = 0 .$$

(II.6.7)

This can be integrated on the domain $(0, \infty)$ and then integration by parts leads to this:

$$(\lambda - \tfrac{1}{2}) \int_0^\infty |a|^2 x^2 dx + (\lambda - \tfrac{3}{2}) \int_0^\infty |b|^2 x^2 dx = 0 .$$

(II.6.8)

(The integration by parts is justified because $\liminf_{x \to 0} x^3 \mu^2$ and $\liminf_{x \to \infty} x^3 \mu^2$ must both vanish when the function $\mu^2$ has finite integral on $(0, \infty)$.) The asserted bounds for $\lambda$ follow from (II.6.8) because the two terms in (II.6.8) will have the same sign otherwise.

By way of a parenthetical remark: If $\lambda = 1$, then

$$(a, b) = \tfrac{1}{x}e^{-kx}(1, 1) \quad and \quad (a, b) = \tfrac{1}{x}e^{kx}(1, -1)$$

(II.6.7)

are two linearly independent solutions to the system depicted in (II.6.6); and the square of the norm of the former is $x^2 dx$–integrable on $(0,\infty)$.

### III. THE OPERATOR $\mathcal{D}$ ON $(0,\infty) \times Y$

Suppose now that Y is a compact, oriented Riemannian 3-manifold. The sections in this 'lecture series' consider the operator $\mathcal{D}$ on $(0,\infty) \times Y$ as defined by a solution with singular asymptotic conditions as $t \to 0$ that are determined by a knot or link in Y. The solution is also required to converge in a suitable sense as $t \to \infty$ so as to define a flat $Sl(2;\mathbb{C})$ connection.

### 1. The Kapustin-Witten linearization on $(0,\infty) \times Y$

Let Y denote a compact, oriented Riemanian 3-manifold. The upcoming Theorem III.1 concerns the operator $\mathcal{D}$ as depicted in (1.2) for the case when $(A,\mathfrak{a})$ is a pair of connection on a principal SU(2) bundle over $(0,\infty) \times Y$ (denoted by P) and ad(P)-valued section of T*Y over the same domain $(0,\infty) \times Y$. Theorem III.1 refers to the $(0,\infty) \times Y$ versions of the Hilbert spaces $\mathbb{H}$ and $\mathbb{L}$ that are defined in Section I.1b.

Thereom III.1 considers only pairs $(A, \mathfrak{a})$ that obey three constraints; two of these constrain $(A,\mathfrak{a})$ where t is much less than 1, and the third constrains $(A, \mathfrak{a})$ where t is much greater than 1. The first constraint says in effect that $(A,\mathfrak{a})$ should look like the Nahm pole solution from (II.1.2) at distance greater than $\mathcal{O}(t)$ from the knot.

**CONSTRAINT 1**: *There is an isometric isomorphism to be denoted by $\tau$ from TY to ad(P) on Y–K and, given $\varepsilon > 0$, there exists $R_\varepsilon > 1$ and a positive time $t_\varepsilon$ such that when $\tau$ is viewed as an* ad(P)*-valued 1-form, then*

$$|\mathfrak{a} + \tfrac{1}{2t}\tau| + |\nabla_A \tau| < \tfrac{\varepsilon}{t}$$

*on the part of $(0,t_\varepsilon) \times Y$ where the distance to K is greater than $R_\varepsilon t$.*

The second constraint says that $(A,\mathfrak{a})$ should look like one of Witten's model solutions from (I.1.2) at distance less than $\mathcal{O}(t)$ from the knot or link when t is very small. To set notation: The knot or link is denoted by K; it will be called a 'knot' even if it has more than one component. (It is compact, 1-dimensional submanifold.) The option that $K = \emptyset$ is allowed in which case the second constraint is vacuous.

The second constraint refers to a fixed radius tubular neighorhood of K which it denotes by $\mathcal{N}_K$. The radius of $\mathcal{N}_K$ is denoted by $r_0$. The radius $r_0$ disk centered at the origin in $\mathbb{R}^2$ is denoted by $D_0$. The disk $D_0$ is given its Euclidean metric in what follows. Each component of K labels a copy of the circle $S^1$ which is the metric from its identification as $\mathbb{R}/(\ell\mathbb{Z})$ with $\ell$ denoting the length of the component in question.

More notation: Suppose for the moment that X is a space whose first homology with $\mathbb{Z}/2$ coefficients is non-zero. Relevant examples are $(0,\infty) \times D_0 \times S^1$ and its subspace where $|z| \geq Rt$. Let $P_1$ and $P_2$ denote principal bundles over X. Pairs $(A^1, \mathfrak{a}^1)$ and $(A^2, \mathfrak{a}^2)$ of connection on $P_{*=1,2}$ and section of $\text{ad}(P_{*=1,2}) \otimes T^*X$ can be equivalent over any given small radius ball via an isomorphism between $P_1$ and $P_2$ but not globally equivalent. In this case, there is none-the-less an isomorphism between these two bundles that is defined modulo $\{\pm 1\}$ that identifies the two pairs (the subgroup $\{\pm 1\}$ is the center of SU(2)). Said formally: An isomorphism modulo $\{\pm 1\}$ is an isomorphism between the principal SO(3) bundle $P_1/\{\pm 1\}$ and $P_2/\{\pm 1\}$.

A further bit of notation: The product principal SU(2) bundle on $\mathbb{R}^2 \times S^1$ is denoted by $P_0$; and, given a positive integer to be denoted by *m*, the integer *m* version of the model Kapustin-Witten solution from (1.12) is denoted by $(A^{(m)}, \mathfrak{a}^{(m)})$. This is a pair of connection on $P_0$ and $\text{ad}(P_0)$-valued section of $T^*(\mathbb{R}^2 \times S^1)$ over $(0,\infty) \times \mathbb{R}^2 \times S^1$.

The second constraint follows

**CONSTRAINT 2**: *There is a data set consisting of the following assignment to each component of* K*:*
a) *A positive integer to be denoted by* m.
b) *An orientation preserving diffeomorphism from the component's version of* $D_0 \times S^1$ *to* $\mathcal{N}_K$ *which maps* $\{0\} \times S^1$ *to* K *with isometric differential on* $\{0\} \times S^1$.
c) *An isomorphism modulo* $\{\pm 1\}$ *over* $(0,\infty) \times D_0 \times S^1$ *between between* $P_0$ *and the pull-back of* P *via the diffeomorphism in Item b).*
*This data set is such that given* $\varepsilon > 0$ *and* $R > 1$, *there exists a positive time* $t_{\varepsilon,R}$ *with the following significance: For any component of* K, *if* $(A,\mathfrak{a})$ *is viewed as a pair of connection on* $P_0$ *and* $\text{ad}(P_0)$-*valued 1-form over the domain* $(0,\infty) \times D_0 \times S^1$ *using pull-back via the diffeomorphism from Item b) and the isomorphism in Item c), then*

$$|\mathfrak{a} - \mathfrak{a}^{(m)}| \leq \tfrac{\varepsilon}{t} \quad \text{and} \quad |A - A^{(m)}| \leq \tfrac{\varepsilon}{t}$$

*on the part of* $(0,\infty) \times D_0 \times S^1$ *where* $t \leq t_{\varepsilon,R}$ *and* $|z| \leq Rt$.

There is one large t constraint on $(A, \mathfrak{a})$. To set the stage: Suppose for the moment that $P^\infty$ is a principal SU(2) bundle over Y, that $A^\infty$ is a connection on $P^\infty$, and that $\mathfrak{a}^\infty$ is an $\text{ad}(P^\infty)$-valued 1-form on Y. These define a corresponding $Sl(2;\mathbb{C})$ connection which is $A^\infty + i\mathfrak{a}^\infty$. This connection is denoted by $\mathbb{A}$ and its associated exterior covariant derivative by $d_\mathbb{A}$. The latter maps $ad(P^\infty) \otimes_\mathbb{R} \mathbb{C}$ valued functions to like valued 1-forms, and these 1-forms to $ad(P^\infty) \otimes_\mathbb{R} \mathbb{C}$ valued 2-forms, and so on up to 3-forms. (The bundle $ad(P^\infty) \otimes_\mathbb{R} \mathbb{C}$ is associated to $P^\infty$ via the adjoint representation of SU(2) on $\mathfrak{sl}(2;\mathbb{C})$.) The formal $L^2$ adjoint of $d_\mathbb{A}$ is denoted by $d_\mathbb{A}^\dagger$. An instance of this gives a linear operator from $ad(P^\infty) \otimes_\mathbb{R} \mathbb{C}$ valued 1-forms to $ad(P^\infty) \otimes_\mathbb{R} \mathbb{C}$ valued functions. If $\mathbb{A}$ is a flat SL(2;$\mathbb{C}$) connection, then $d_\mathbb{A}^2 \equiv 0$ in which case the associated cohomology groups can be defined by the usual rule whereby

$$\{H^p = \ker(d_\mathbb{A})/\mathrm{im}(d_\mathbb{A})\}_{p=0,1,2,3}.$$

(III.1.1)

The connection $\mathbb{A}$ is said to be *irreducible* when $H^0 = \{0\}$. It is said to be *regular* when $H^1 = \{0\}$. (As with DeRham cohomology, $H^2$ is isomorphic to $H^1$ and $H^3$ to $H^0$.)

To continue the stage setting: Keep in mind that any principal SU(2) bundle over $(0,\infty) \times Y$ is isomorphic to one that is pulled back from Y via the projection map that sends (t, y) to y. This projection map is denoted by $\pi$ in what follows.

**CONSTRAINT 3**: *There exists a data set consisting of*:
a) *A principal SU(2) bundle* $P^\infty \to Y$.
b) *A regular flat SL(2;$\mathbb{C}$) connection on* $P^\infty \times_{SU(2)} SL(2;\mathbb{C})$ *to be denoted by* $\mathbb{A}$.
c) *Given* $\varepsilon > 0$, *a positive time* $t_\varepsilon$.
d) *An isomorphism between* P *and* $\pi^*P^\infty$.
*This data has the following significance: The isomorphism in Item d) identifies* $A + i\mathfrak{a}$ *with an SL(2;$\mathbb{C}$) connection on* $[t_\varepsilon, \infty) \times Y$ *that differs pointwise from* $\pi^*\mathbb{A}$ *by at most* $\frac{\varepsilon}{t}$.

The upcoming Theorem IV.1 in the last 'lecture series' gives conditions that guarantee the $o(\frac{1}{t})$ convergence of $A + i\mathfrak{a}$ to a flat SL(2:$\mathbb{C}$) connection when $(A, \mathfrak{a})$ obeys the Kapustin-Witten equations.

What follows is the promised Theorem III.1.

**Theorem III.1**: *Suppose that* $(A, \mathfrak{a})$ *is a pair of connection on the bundle* P *and* ad(P)*-valued section of* $T^*Y$ *over the domain* $(0, \infty) \times Y$ *that obeys* CONSTRAINTS 1, 2 *and* 3. *The corresponding operator* $\mathcal{D}$ *defines a bounded, Fredholm map from* $\mathbb{H}$ *to* $\mathbb{L}$ *whose index is positive and at most twice the complex dimension of the space* $H^0$ *with the latter defined by the limit flat SL(2;$\mathbb{C}$)-connection* $\mathbb{A}$ *from* CONSTRAINT 3. *In particular*:
- *If this limit flat connection is irreducible, then the index is zero*.
- *If the limit flat connection is reducible, and if* $(A, \mathfrak{a})$ *obeys the equations in (I.2.1) in addition to* CONSTRAINTS 1-3, *then the index is the (complex) dimension of* $H^0$.

The proof that $\mathcal{D}$ has closed range and finite dimensional kernel is Section III.3. Section III.2 supplies some preliminary observations for that proof. The proof that $\mathcal{D}$ is Fredholm occupies Section III.4; and Section III.5 computes the index.

With regards to the second bullet: If $\dim_{\mathbb{C}} H^0 > 0$, then the positive index is due to the existence of automorphisms of P on $(0, \infty) \times Y$ that limit to nontrivial, $\mathbb{A}$-covariantly constant automorphisms as $t \to \infty$. This is sometimes called the *extended kernel*.

By way of a second comment: The small t constraints (CONSTRAINTS 1 and 2) for Proposition III.1 can be replaced by suitable curvature integral constraints, for example those for small time in [T?] in the case when the knot K is absent and $(A, \mathfrak{a})$ obeys the Kapustin-Witten equations.

## 2. Constructions where t is small

This section states and proves two lemmas about $\mathcal{D}$ where t is small. (They are used in the subsequent sections to prove Theorem III.1.)

### a) Isolating the small t part of $\mathcal{D}$

A cut-off function is used to separate the analysis of $\mathcal{D}$ at very small values of t from the analsysis at large values of t. To elaborate: Fix $\delta \in (0, \frac{1}{2}]$ and use it to define cut-off functions on $(0, \infty)$ to be denoted by $\chi_\delta$ by the rule

$$t \to \chi_\delta(t) = \chi(2 + \tfrac{\ln t}{|\ln \delta|}) \ .$$

(III.2.1)

To be sure: The function $\chi_\delta$ is 1 where $t \leq \delta^2$ and it is zero where $t \geq \delta$.

If $\psi$ is from $\mathbb{H}$, then it can be written as $\psi_0 + \psi_1$ with

$$\psi_0 = \chi_\delta \psi \quad and \quad \psi_1 = (1 - \chi_\delta) \psi \ .$$

(III.2.2)

These definitions lead to the identities

- $\mathcal{D}\psi_0 = \chi_\varepsilon \mathcal{D}\psi + (\tfrac{d}{dt} \chi_\delta)(\psi_0 + \psi_1)$ ,
- $\mathcal{D}\psi_1 = (1 - \chi_\varepsilon)\mathcal{D}\psi - (\tfrac{d}{dt} \chi_\delta)(\psi_0 + \psi_1)$ .

(III.2.3)

The latter then lead to the inequality:

$$|\mathcal{D}\psi|^2 \geq c^{-1}(|\mathcal{D}\psi_0|^2 + |\mathcal{D}\psi_1|^2) - c \tfrac{1}{|\ln \delta|^2 t^2}(|\psi_0|^2 + |\psi_1|^2)$$

(III.2.4)

where $c$ denotes here and in what follows a number that is greater than 1 and independent of both $\varepsilon$ and $\psi$. (Its value can be assumed to increase between successive appearances.). Integrate both sides of (III.2.4) over the domain $(0,\infty) \times Y$ and use Hardy's inequality (see (I.2.3)) to derive the next inequality:

$$\int_{(0,\infty)\times Y} |\mathcal{D}\psi|^2 \geq c^{-1} \int_{(0,\infty)\times Y} |\mathcal{D}\psi_0|^2 + c^{-1} \int_{(0,\infty)\times Y} |\mathcal{D}\psi_1|^2 - c \frac{1}{|\ln\delta|^2} \|\psi_0\|_{\mathbb{H}}^2 - c \frac{1}{|\ln\delta|^2} \int_{(\frac{1}{T^2},\infty)\times Y} \frac{1}{t^2}|\psi_1|^2 \ .$$

(III.2.5)

The $|\mathcal{D}\psi_1|^2$ integral in (III.2.5) will be studied independently in subsequent subsections. The following lemma states what is needed about the $|\mathcal{D}\psi_0|^2$ integral.

**Lemma III.2.1**: *Given a pair* $(A,\mathfrak{a})$ *obeying* CONSTRAINTS 1 *and* 2, *there exists* $\kappa > 1$ *with the following significance: If* $\delta < \frac{1}{\kappa}$ *and supposing that* $\psi \in \mathbb{H}$, *then*

$$\|\mathcal{D}\psi_0\|_{\mathbb{L}} \geq \tfrac{1}{\kappa} \|\psi_0\|_{\mathbb{H}} \ .$$

*Proof of Lemma III.2.1*: When $\varepsilon$ is very small, the domain of the left integral is where $(A,\mathfrak{a})$ is everywhere close to either one of the $(A^{(m)},\mathfrak{a}^{(m)})$ pairs or a pair that looks near any given point like the Nahm pole solution. This behavior is exploited with the help of a second cut-off function to separate the analysis near the K from the analysis away from K. This second function requires the choice of a number $r \in (0, \frac{1}{100}r_0)$. Having chosen r, define the function $\varpi_r$ on Y by the rule

$$\varpi_r(\cdot) = \chi(\tfrac{\mathrm{dist}(\cdot,K)}{r} - 1) \ .$$

(III.2.6)

This function is equal to 1 where the distance to K is less than r and it is equal to zero where the distance to K is greater than 2r.

Now write $\psi_0$ as $\psi_{00} + \psi_{0K}$ where

$$\psi_{00} = (1 - \varpi_r)\psi_0 \quad and \quad \psi_{0K} = \varpi_r \psi_0 \ .$$

(III.2.7)

The same sort of analysis that lead from (III.2.2) to (III.2.5) leads in this case to this:

$$\int_{(0,\infty)\times Y} |\mathcal{D}\psi_0|^2 \geq c^{-1} \int_{(0,\infty)\times Y} |\mathcal{D}\psi_{00}|^2 + c^{-1} \int_{(0,\infty)\times Y} |\mathcal{D}\psi_{0K}|^2 - c \frac{\delta^2}{r^2} \|\psi_0\|_{\mathbb{H}}^2 \ .$$

(III.2.8)

With regards to the right most term on the right hand side of this inequality: This term with the factors of $c^{-1}$ in front of the two integrals on the right hand side of (III.2.8) account for terms with derivatives of $\varpi_r$ multiplying $\psi_0$ that arise when $\mathcal{D}\psi_0$ is compared with the sum of $\mathcal{D}\psi_{00}$ and $\mathcal{D}\psi_{0K}$. The fact that all such terms are accounted for in this way follows because the $(0,\infty)\times Y$ integral of $|d\varpi_t|^2|\psi_0|^2$ is bounded by $c\frac{1}{r^2}\|\psi_0\|_{\mathbb{L}}^2$ and thus by $c\frac{\delta^2}{r^2}\|\frac{1}{t^2}\psi_0\|_{\mathbb{L}}^2$ since $t \le \varepsilon$ on the support of $\psi_0$. Meanwhile, $\frac{\delta^2}{r^2}\|\frac{1}{t^2}\psi_0\|_{\mathbb{L}}^2$ is bounded by $\frac{\delta^2}{r^2}\|\psi_0\|_{\mathbb{H}}^2$ via Hardy's inequality.

The next lemma summarizes what will be said about the two integrals on the right hand side of (III.2.8).

**Lemma III.2.2**: *Given a pair* $(A,\mathfrak{a})$ *obeying the conditions in* CONSTRAINTS *1 and 2, there exists* $\kappa > 1$ *with the following significance: If* $r < \frac{1}{\kappa}$ *and* $\delta < \frac{1}{\kappa^2}$ *and supposing that* $\psi \in \mathbb{H}$, *then*
- $\|\mathcal{D}\psi_{00}\|_{\mathbb{L}} \ge \frac{1}{\kappa}\|\psi_{00}\|_{\mathbb{H}}$.
- $\|\mathcal{D}\psi_{0K}\|_{\mathbb{L}} \ge \frac{1}{\kappa}\|\psi_{0K}\|_{\mathbb{H}}$.

The proof of lemma III.2.2 occupies the next subsection. Assume it to be true and, as argued directly, Lemma III.2.1 is a consequence.

To obtain Lemma III.2.1, note first that

$$\|\psi_{00}\|_{\mathbb{H}}^2 + \|\psi_{0K}\|_{\mathbb{H}}^2 \ge c^{-1}\|\psi_0\|_{\mathbb{H}}^2 - c\frac{\delta^2}{r^2}\|\psi_0\|_{\mathbb{H}}^2.$$

(III.2.9)

In this regard: All terms with derivatives of $\varpi_t$ multiplying $\psi_0$ that arise when comparing the $\mathbb{H}$ norm of $\psi_0$ with those of $\psi_{00}$ and $\psi_{0K}$ are accounted for by the right most term on the right hand side of (III.2.9). (The justification for this claim is the same as for the analogous term in (III.2.9).)

If r is sufficiently small and then $\delta$ is chosen sufficiently small (given r), then the inequalities in Lemma III.2.2 when used with (III.2.8) and then (III.2.9) directly give the assertion of Lemma III.2.1.

**b) Proof of Lemma III.2.2**

The proof of this lemma has eight parts. Part 1 proves the top bullet and Parts 2-8 prove the second bullet.

*Part 1*: The analysis for the top bullet of Lemma III.2.2 starts by rewriting the operator $\mathcal{D}$ using the data from CONSTRAINT 1. To set the stage for this, fix a positive $\varepsilon$ and $R > 1$ so as to obtain the number $t_{\varepsilon,R}$ and the section $\tau$ from that constraint. Let $X_{\varepsilon,R}$ denote the subset in $(0, \infty) \times Y$ where $t < t_{\varepsilon,R}$ and where the distance to K is greater than $Rt$. Viewing $\tau$ first as an ad(P)-valued 1-form on $X_{\varepsilon,R}$, write $\mathfrak{a}$ as $-\frac{1}{2t}\tau + \mathfrak{v}$. Viewing $\tau$ next as an isometric homomorphism from TY to ad(P), and letting $\Gamma$ denote the push-forward via $\tau$ of the Levi-Civita connection on TY, write A as $\Gamma + \mathrm{A}$. By assumption the norms of $\mathfrak{v}$ and $\mathrm{A}$ are at most $c\frac{\varepsilon}{t}$ on $X_{\varepsilon,R}$. As a consequence, the opertor $\mathcal{D}$ on $X_{\varepsilon,R}$ can be written as

$$\mathcal{D} = \tfrac{\partial}{\partial t} + \gamma_i \nabla_{\Gamma i} - \tfrac{1}{2t}\rho_i[\tau_i, \cdot\,] + \mathfrak{r}$$

(III.2.10)

where the notation has $\nabla_\Gamma$ denoting the covariant derivative on ad(P) valued tensors that is defined by $\Gamma$; it has $\mathfrak{r}$ denoting an endomorphism whose norm obeys $|\mathfrak{r}| < c\frac{\varepsilon}{t}$.

Supposing that the support of $\psi_{00}$ is in $X_{\varepsilon,R}$, then the next inequality is a direct consequence of what was just said about $\mathcal{D}$.

$$\int_{(0,\infty)\times Y} |\mathcal{D}\psi_{00}|^2 \geq c^{-1} \int_{(0,\infty)\times Y} \left( |\tfrac{\partial}{\partial t}\psi_{00}|^2 + \sum_{i=1}^{3} |\nabla_{\Gamma i}\psi_{00}|^2 \right)$$
$$+ c^{-1} \int_{(0,\infty)\times Y} \left( \tfrac{1}{t}\langle \psi_{00}, \{-\tfrac{1}{2t}\tau, \psi_{00}\}\rangle + |\{\tfrac{1}{2t}\tau, \psi_{00}\}|^2 \right) - c\varepsilon^2 \int_{(0,\infty)\times Y} \tfrac{1}{t^2}|\psi_{00}|^2$$

(III.2.11)

Granted this inequality then (II.3.3) and Hardy's inequality in (I.2.3) lead from (III.2.11) to the following if $\varepsilon < c^{-1}$ (which will henceforth be assumed):

$$\int_{(0,\infty)\times Y} |\mathcal{D}\psi_{00}|^2 \geq c^{-1} \int_{(0,\infty)\times Y} \left( |\tfrac{\partial}{\partial t}\psi_{00}|^2 + \sum_{i=1}^{3} |\nabla_{\Gamma i}\psi_{00}|^2 \right).$$

(III.2.12)

Moreover, because $|\mathrm{A}|$ is bounded by $c\frac{\varepsilon}{t}$ and because of Hardy's inequality in (I.2.3), the derivatives $\frac{\partial}{\partial t}$ and $\nabla_\Gamma$ that appear in (III.2.12) can be replaced by the A-covariant derivative if $\varepsilon < c^{-1}$. (This is a larger version of $c$ than before; it is henceforth assumed that $\varepsilon$ obeys this new bound.) Last but not least: Because of the same Hardy's inequality and because $|\mathfrak{a}| \leq c\frac{1}{t}$, the latter version of (III.2.12) leads to a final inequality:

$$\|\mathcal{D}\psi_{00}\|_{\mathbb{L}} \geq c^{-1}\|\psi\|_{\mathbb{H}}.$$

(III.2.13)

Keep in mind though that this inequality is predicated on $\varepsilon$ being less than $c^{-1}$ and on the support of $\psi_{00}$ being in $X_{\varepsilon,R}$ (which is the part of $(0, t_{\varepsilon,R}) \times Y$ where the distance to K is greater than R.)

*Part 2*: This part of the proof starts the story for the second bullet of Lemma III.2.2. The analysis uses what is said by CONSTRAINT 2 in a number of ways. To start, the diffeomorphism in Item b) of CONSTRAINT 2 and the isomorphism in Item c) of CONSTRAINT 2 will be used henceforth (without comment for the most part) to identify a given component of $(0,\infty) \times \mathcal{N}_K$ with $(0,\infty) \times D_0 \times S^1$ and the bundle P over this component of $(0,\infty) \times \mathcal{N}_K$ with the product principal SU(2) bundle over $(0,\infty) \times D_0 \times S^1$. In this way, the pair $(A, \mathfrak{a})$ becomes a pair consisting of a connection on the product bundle over $(0,\infty) \times D_0 \times S^1$ and $\mathfrak{su}(2)$-valued section over this same domain of $T^*(D_0 \times S^1)$.

The CONSTRAINT 2 diffeomorphism and isomorphism can also be used to view the operator $\mathcal{D}$ as a differential operator over the domain $(0,\infty) \times D_0 \times S^1$. Viewed in this light, it is almost one of the versions of $\mathcal{D}$ that are discussed in Sections II.4-II.6. To elaborate: First, the metric that $\mathcal{D}$ uses is not the Euclidean metric, it is the pull-back of the metric from Y via the diffeomorphism in Item b) of CONSTRAINT 2. Even so, this metric differs pointwise from the Euclidean metric by at most $c|z|$ because it is assumed to be isometric along $\{0\} \times S^1$. What follows is a consequence: Let $\mathcal{D}^0$ denote the Euclidean metric version of $\mathcal{D}$ as defined using the pull-back of $(A,\mathfrak{a})$ via the diffeomorphism from Item b) of CONSTRAINT 2 and the isomorphism of Item c) of CONSTRAINT 2. Then

$$\mathcal{D} = \mathcal{D}^0 + \mathfrak{r}_1(\nabla) + \mathfrak{r}_0$$

(III.2.14)

where $\mathfrak{r}_0$ and $\mathfrak{r}_1$ are endomorphisms whose norms obey $|\mathfrak{r}_1| \leq c|z|$ and $|\mathfrak{r}_0| \leq c$. (Here, $z = z_1 + iz_2$ is the complex Euclidean coordinate for $D_0$.)

*Part 3*: At the same time, the pair $(A,\mathfrak{a})$ differs little from the model $(A^{(m)}, \mathfrak{a}^{(m)})$ when $(A,\mathfrak{a})$ is viewed on the $|z| < Rt$ part of the domain in $(0, \infty) \times D_0 \times S^1$ via the diffeomorphism from Item b) of CONSTRAINT 2 and the isomorphism in Item c) of CONSTRAINT 2. To be precise here: The pair $(A,\mathfrak{a})$ differs from the model solution by at most $c \frac{\varepsilon}{t}$ where $t \leq t_{\varepsilon,R}$ and $|z| \leq Rt$.

*Part 4*: Meanwhile, the pair $(A,\mathfrak{a})$ is described by CONSTRAINT SET 1 where $|z|$ is between $Rt$ and $r$ and $t < t_{\varepsilon,R}$ if it is assumed that $\delta < t_{\varepsilon,R}$). Supposing that such is the case (which will henceforth be done), then $(A,\mathfrak{a})$ can be written on this part of $(0, \infty) \times D_0 \times S^1$ so as to look very much like $(\Gamma_0, -\frac{1}{2t}\tau_0)$ with $\tau_0$ denoting a *Euclidean metric* isometry from the tangent space of $(D_0 - 0) \times S^1$ to $\mathrm{ad}(P_0)$ and with $\Gamma_0$ denoting the (flat) connection that is induced on $\mathrm{ad}(P_0)$ from the Euclidean metric's Levi-Civita connection on $T(D_0 \times S^1)$. To say this precisely: The isometry $\tau_0$ can be chosen so that $(A,\mathfrak{a})$ differs from $(\Gamma_0, -\frac{1}{2t}\tau_0)$ by at most $c \frac{\varepsilon + r}{t}$ where $t < t_{\varepsilon,r}$ and where $Rt \leq |z| \leq r$.

*Part 5*: There is also a Euclidean metric isometry from $T(D_0 \times S^1)$ to $ad(P_0)$ to be denoted by $\tau_{(m)}$ that writes $(A^{(m)}, \mathfrak{a}^{(m)})$ as $(\Gamma_0, -\frac{1}{2t}\tau_{(m)}) + \mathfrak{h}$ with the norm of $\mathfrak{h}$ obeying $|\mathfrak{h}| \leq c \frac{1}{R^2 t}$ where $|z| \geq R t$ on $(0, \infty) \times \mathbb{R}^2 \times S^1$.

*Part 6*: Assuming that $\delta < t_{\varepsilon, R}$, then the observations in Parts 3-5 imply this: There is an isomorphism modulo $\{\pm 1\}$ on the $|z| < r$ part of $(0, \delta] \times D_0 \times S^1$ between $P_0$ and the pull-back of P via the diffeomorphism from Item b) of CONSTRAINT 2 that writes $(A, \mathfrak{a})$ as a pair (which will still denoted by $(A, \mathfrak{a})$) that obeys

$$(A, \mathfrak{a}) = (A^{(m)}, \mathfrak{a}^{(m)}) + \mathfrak{w} \tag{III.2.15}$$

with $\mathfrak{w}$ obeying $|\mathfrak{w}| \leq c(\varepsilon + r + \frac{1}{R^2})\frac{1}{t}$.

*Part 7*: With regards to $\psi_{0K}$: The modulo $\{\pm\}$ isomorphism from the Part 6 maps $\psi_{0K}$ to some $\psi'$ and it maps $\mathcal{D}\psi'$ to an element that can be written as

$$\mathcal{D}^{0(m)}\psi' + \mathfrak{k}_1(\nabla)\psi' + \mathfrak{k}_0\psi' \tag{III.2.16}$$

where $\mathcal{D}^{0(m)}$ is the Euclidean space version of $\mathcal{D}$ that is defined by $(A^{(m)}, \mathfrak{a}^{(m)})$ as described in Sections II.4–II.6; and where $\mathfrak{k}_1$ and $\mathfrak{k}_0$ denote endomorphisms with norms obeying $|\mathfrak{k}_1| \leq c|z|$ (which is less than $cr$) and $|\mathfrak{k}_0| \leq c(\varepsilon + r + \frac{1}{R^2})\frac{1}{t}$.

*Part 8*: According to what Sections II.4-II.6 say, and by virtue of Hardy's inequality,

$$\|\mathcal{D}^{0(m)}\psi'\|_{\mathbb{L}0} \geq c^{-1}\|\psi'\|_{\mathbb{H}0(m)}, \tag{III.2.17}$$

where $\|\cdot\|_{\mathbb{L}0}$ denotes the $\mathbb{L}$-norm (the $L^2$ norm) that is defined by the Euclidean metric and where $\|\cdot\|_{\mathbb{H}0(m)}$ denotes the $\mathbb{H}$-norm as defined by the Euclidean metric and $(A^{(m)}, \mathfrak{a}^{(m)})$. Therefore, if $\varepsilon < c^{-1}$ and $r < c^{-1}$ and $R > c$, then by virtue of (III.2.14)– (III.2.17), and again by Hardy's inequality:

$$\|\mathcal{D}\psi_{0K}\|_{\mathbb{L}} \geq c^{-1}\|\psi_{0K}\|_{\mathbb{H}}. \tag{III.2.18}$$

As before, this inequality is predictated on $\varepsilon$ being small, less than $c^{-1}$ and also $r < c^{-1}$ and $R > c^{-1}$; and then $\delta < c^{-1} t_{\varepsilon, E}$.

## 3. Constructions where t is large

This section states and proves two lemmas about $\mathcal{D}$ where t is large.

### a) Isolating the large t part of $\mathcal{D}$

This subsection isolates the large t part of the $|\mathcal{D}\psi_1|^2$ integral that appears in (III.2.5). This is done with the help of yet another another cut-off function. The definition of this new function requires the choice of $T > 1$. Given that, set $\varsigma_T$ to denote the function

$$t \to \varsigma_T(t) = \chi(2 - \tfrac{\ln t}{\ln T}) \,.$$

(III.3.1)

To be sure: This function is equal to 1 where $t \geq T^2$ and it is equal to zero where $t \leq T$. Supposing that $\psi_1$ is from $\mathbb{H}$, write it as $\psi_{11} + \psi_{1\infty}$ where

$$\psi_{11} = (1 - \varsigma_T)\psi_1 \quad and \quad \psi_{1\infty} = \varsigma_T \psi_1 \,.$$

(III.3.2)

There are identities for $\mathcal{D}\psi_{11}$ and $\mathcal{D}\psi_{1\infty}$ that mirror (III.2.3) that lead to the following large t mirror of (III.2.5):

$$\int_{(0,\infty)\times Y} |\mathcal{D}\psi_1|^2 \geq c^{-1} \int_{(0,\infty)\times Y} |\mathcal{D}\psi_{1\infty}|^2 + c^{-1} \int_{(0,\infty)\times Y} |\mathcal{D}\psi_{11}|^2 - c\tfrac{1}{|\ln T|^2}\|\psi_{1\infty}\|_{\mathbb{H}}^2$$
$$- c\tfrac{1}{|\ln T|^2} \int_{[\tfrac{1}{T^2},T^2]\times Y} \tfrac{1}{t^2}|\psi_{11}|^2 \,.$$

(III.3.3)

The next lemma says what is needed with regards to the $|\mathcal{D}\psi_{1\infty}|^2$ integral in (III.3.3).

**Lemma III.3.1**: *Given a pair* $(A, \mathfrak{a})$ *obeying the conditions in* CONSTRAINT 3, *there exists* $\kappa > 1$ *with the following significance: If* $T > \kappa$ *and supposing that* $\psi \in \mathbb{H}$, *then*

$$\|\mathcal{D}\psi_{1\infty}\|_{\mathbb{L}} \geq \tfrac{1}{\kappa}\|\psi_{1\infty}\|_{\mathbb{H}} \,.$$

*Proof of Lemma III.3.1*: Use the isomorphism from CONSTRAINT 3 between P and $\pi^*P^\infty$ to identify these two bundles. Write the $\pi^*$-pull-back of the flat connection from CONSTRAINT 3 as $A^\infty + i\mathfrak{a}^\infty$ with $A^\infty$ denoting a connection on P over $(0, \infty) \times Y$ and with $\mathfrak{a}^\infty$ denoting an ad(P)-1-form on $(0,\infty)\times Y$ with no dt component.

Fix $\varepsilon > 0$ but small and having done that, suppose henceforth that T is greater than the version of $t_\varepsilon$ that appears in CONSTRAINT 3. Now use the automorphism from CONSTRAINT 3 to write $(A, \mathfrak{a})$ as $(A^\infty, \mathfrak{a}^\infty) + \mathfrak{z}$ with $\mathfrak{z}$ obeying $|\mathfrak{z}| \leq c \frac{\varepsilon}{\mathfrak{t}}$. This same isomorphism writes $\mathcal{D}$ where $t \geq t_\varepsilon$ as

$$\mathcal{D} = \mathcal{D}^\infty + \mathfrak{x}$$

(III.3.4)

with $\mathcal{D}^\infty$ denoting the $(A^\infty, \mathfrak{a}^\infty)$ version of $\mathcal{D}$ and with $\mathfrak{x}$ denoting an endomorphism with norm bounded also by $c \frac{\varepsilon}{\mathfrak{t}}$.

Let $\|\cdot\|_{\mathbb{H}\infty}$ denote the version of the $\mathbb{H}$-norm from (I.1.9) that is defined using the pair $(A^\infty, \mathfrak{a}^\infty)$. By virtue Hardy's inequality (and because $|\mathfrak{z}| \leq c \frac{\varepsilon}{\mathfrak{t}}$), the $\|\cdot\|_{\mathbb{H}}$ and $\|\cdot\|_{\mathbb{H}\infty}$ norms are uniformly equivalent on the space of compactly supported sections of $\mathbb{W}$ that are supported where $t > t_\varepsilon$. (The bundle $\mathbb{W}$ is $\oplus_2 (\mathrm{ad}(P) \oplus (\mathrm{ad}(P) \otimes T^*Y))$.) More to the point, the following inequalities hold on this subspace of sections:

$$(1 - c\varepsilon)\|\cdot\|_{\mathbb{H}} \leq \|\cdot\|_{\mathbb{H}\infty} \leq (1 + c\varepsilon)\|\cdot\|_{\mathbb{H}}.$$

(III.3.5)

As explained momentarily, Lemma III.3.1 follows from the inequality in (III.3.5), the $c \frac{\varepsilon}{\mathfrak{t}}$ bound on the norm of $\mathfrak{x}$ in (III.3.4) and the following lemma about $\mathcal{D}^\infty$ if $\varepsilon < c^{-1}$ and if T is chosen so as to be greater than $t_\varepsilon$.

**Lemma III.3.2**: *Let $\mathbb{A}$ denote a regular, flat $SL(2;\mathbb{C})$ connection on Y that can be written as $A^\infty + i\mathfrak{a}^\infty$ with $A^\infty$ being a connection on a principal $SU(2)$ bundle over Y and with $\mathfrak{a}^\infty$ being a 1-form on Y with values in that associated Lie algebra bundle. If, in addition, $d_{A^\infty} * \mathfrak{a}^\infty = 0$, then there exists $\kappa > 1$ with the following significance: Pull-back the principal $SU(2)$ bundle and $(A^\infty, \mathfrak{a}^\infty)$ to define the operator $\mathcal{D}^\infty$ and the $(A^\infty, \mathfrak{a}^\infty)$ version of the Hilbert space $\mathbb{H}$. If $\psi'$ is in this space, then $\|\mathcal{D}^\infty \psi'\|_{\mathbb{L}} \geq \frac{1}{\kappa} \|\psi'\|_{\mathbb{H}\infty}$.*

This lemma is proved in the next subsection. Assume it to be true for now.

To finish the proof of Lemma III.3.1: If the conditions of Lemma III.3.2 are met by the limit flat connection $\mathbb{A}$ from CONSTRAINT 3, then (III.3.4) leads directly to this:

$$\|\psi_{1\infty}\|_{\mathbb{H}\infty} - c\varepsilon \|\tfrac{1}{\mathfrak{t}} \psi_{1\infty}\|_{\mathbb{L}} \leq c \|\mathcal{D}\psi_{1\infty}\|_{\mathbb{L}}$$

(III.3.6)

if it is assumed that $T > t_\varepsilon$. If $\varepsilon < c^{-1}$, then (III.3.6) with Hardy's inequality leads in turn to the inequality $\|\psi_{1\infty}\|_{\mathbb{H}\infty} \leq c \|\mathcal{D}\psi_{1\infty}\|_{\mathbb{L}}$; and that with (III.3.5) gives the inequality that is asserted by Lemma III.3.1 if $\varepsilon < c^{-1}$.

With regards to the conditions in Lemma III.3.6: They are all guaranteed for the limit flat connection coming from CONSTRAINT 3. Indeed, the only condition that is not explicit in CONSTRAINT 3 is the $d_{A^\infty} * \mathfrak{a}^\infty = 0$ condition. But this follows directly from the fact that $d_A * \mathfrak{a} = 0$ and that $(A, \mathfrak{a})$ after an automorphism limits pointwise as $t \to \infty$ to $(A^\infty, \mathfrak{a}^\infty)$. (There is no need for an assumption about first derivatives converging as $t \to \infty$ because the pointwise convergence assumption in CONSTRAINT 3 already implies that

$$\lim_{t \to 0} \int_{\{t\} \times Y} \langle \mathfrak{a} \wedge * d_A \mathfrak{v} \rangle = \int_Y \langle \mathfrak{a}^\infty \wedge * d_{A^\infty} \mathfrak{v} \rangle$$

(III.3.7)

for any fixed section $\mathfrak{v}$ of $\mathrm{ad}(P^\infty)$ over Y.)

### b) Proof of Lemma III.3.2

The proof of the lemma has eleven parts (they are short).

*Part 1*: Let $P^\infty$ denote a principal SU(2) bundle over Y and let $\mathrm{ad}(P^\infty)_{\mathbb{C}}$ denote the complexification of $\mathrm{ad}(P^\infty)$, thus $\mathrm{ad}(P^\infty) \otimes_{\mathbb{R}} \mathbb{C}$. (It is the associated $\mathfrak{sl}(2; \mathbb{C})$ bundle via the adjoint representation of SU(2) on $\mathfrak{sl}(2; \mathbb{C})$.) Supposing that $\mathbb{A}$ denotes a given flat $SL(2; \mathbb{C})$ connection on $P^\infty \times_{SU(2)} SL(2; \mathbb{C})$, let $d_{\mathbb{A}}$ denote the corresponding exterior covariant derivative. (Remember that this operator has square zero because $\mathbb{A}$ is flat.) Let $d_{\mathbb{A}}^\dagger$ denote the formal adjoint of $d_{\mathbb{A}}$ as defined using the $L^2$-inner product.

*Part 2*: Just to be sure, the $L^2$ inner product between an $\mathrm{ad}(P^\infty)_{\mathbb{C}}$ valued functions $\mathfrak{v}$ and $\mathfrak{w}$ is the pairing

$$\int_Y \langle \mathfrak{v}^* \wedge * \mathfrak{w} \rangle$$

(III.3.8)

where $\mathfrak{v}^*$ is used here to denote $-\mathfrak{v}^\dagger$. The - sign appears with the Hermitian adjoint operation because the linear form $\langle \cdot \rangle$ on the space of $2 \times 2$ complex matrices is $-\frac{1}{2}$ times the trace. The pairing in ((III.3.8) is a positive definite inner product on $C^\infty(Y; \mathrm{ad}(P^\infty)_{\mathbb{C}})$ with this sign convention. The $L^2$ pairing between $\mathrm{ad}(P^\infty)_{\mathbb{C}}$-valued differential forms and tensors is defined analogously with the metric used for the tensor part.

*Part 3*: Of particular interest is the instance where $d_{\mathbb{A}}^\dagger$ maps $C^\infty(Y; \mathrm{ad}(P^\infty)_{\mathbb{C}} \otimes T^*Y)$ back to $C^\infty(Y; \mathrm{ad}(P^\infty)_{\mathbb{C}})$. Hodge theory for the deRham complex with local coefficients says that there is an orthogonal (with respect to the $L^2$-inner product) decomposition

$$C^\infty(Y; ad(P^\infty)_\mathbb{C} \otimes T^*Y) = \ker(d_\mathbb{A}^\dagger) \oplus im(d_\mathbb{A}) \ .$$

(III.3.9)

Hodge theory also asserts that the operator $*d_\mathbb{A}$ maps the $\ker(d_\mathbb{A}^\dagger)$ summand in (III.3.10) to itself; and while doing this, it define a self-adjoint, unbounded operator on the $L^2$ completion of $\ker(d_\mathbb{A}^\dagger)$, a self-adjoint operator with pure point spectrum having no accumulation points and finite multiplicities. As a consequence, there exists a positive number (to be denoted by $\lambda$) with the following significance: If $\eta$ is from $C^\infty(Y; ad(P^\infty)_\mathbb{C} \otimes T^*Y)$ and if it is $L^2$-orthogonal to the kernel of $d_\mathbb{A}$ on the $\ker(d_\mathbb{A}^\dagger)$ summand in (III.3.10), then

$$\int_Y (|d_\mathbb{A}\eta|^2 + |d_\mathbb{A}^\dagger \eta|^2) \geq \lambda \int_Y |\eta|^2 \ .$$

(III.3.11)

If the connection $\mathbb{A}$ is regular (which means that $H^1 = \{0\}$), then (III.3.11) holds for all $\eta \in C^\infty(Y; ad(P^\infty)_\mathbb{C} \otimes T^*Y)$ because $H^1$ is isomorphic to the kernel of $d_\mathbb{A}$ on $\ker(d_\mathbb{A}^\dagger)$.

*Part 4*: If $\mathbb{A}$ is irreducible, which is to say that the kernel of $d_\mathbb{A}$ on $C^\infty(Y; ad(P^\infty)_\mathbb{C})$ is $\{0\}$, then Hodge theory says that there is a positive number (also denoted by $\lambda$) such that

$$\int_Y |d_\mathbb{A} v|^2 \geq \lambda \int_Y |v|^2 \ .$$

(III.3.12)

for all $v \in C^\infty(Y; ad(P^\infty)_\mathbb{C})$. If $\mathbb{A}$ is not irreducible, then (III.3.12) holds provided that $\eta$ is $L^2$-orthogonal to the kernel of $d_\mathbb{A}$ in $C^\infty(Y; ad(P^\infty)_\mathbb{C})$.

*Part 5*: The group $SL(2;\mathbb{C})$ has an involution that is defined as follows: Let g denote a given element in $SL(2;\mathbb{C})$, thus a $\mathbb{C}$-valued, $2 \times 2$ matrix with determinant 1. The involution sends g to the Hermitian conjugate of $g^{-1}$. The matrix $(g^{-1})^\dagger$ is denoted subsquently by $g^*$.) This involution induces an involution on the Lie algebra $\mathfrak{sl}(2;\mathbb{C})$ which sends a matrix $\mathfrak{v}$ to $-\mathfrak{v}^\dagger$ (which is the same $(\cdot)^*$ used in (III.3.9)).

The involution $(\cdot)^*$ on $SL(2;\mathbb{C})$ fixes the $SU(2)$ subgroup; moreover: If h is in $SU(2)$ and g is in $SL(2;\mathbb{C})$, then $(hgh^{-1})^* = hg^*h^{-1}$. As a consequence, the bundle $P^\infty \times_{SU(2)} SL(2;\mathbb{C})$ inherits this involution. If $\mathbb{A}$ denotes a connection on this bundle, then the involution pulls it back to another connection: Write $\mathbb{A}$ as $A + i\mathfrak{a}$ with A being a connection on P and with $\mathfrak{a}$ being an ad(P)-valued 1-form; its pull-back is then $A - i\mathfrak{a}$.

(The action on $\mathbb{A}$ is denoted by $\mathbb{A}^*$). If $\mathbb{A}$ is a flat connection, then so is $\mathbb{A}^*$. In addition, if $\mathbb{A}$ is regular and/or irreducible, then likewise so is $\mathbb{A}^*$. Also: The respective versions of $\lambda$ that appear in (III.3.6) and (III.3.7) for the connection $\mathbb{A}$ and $\mathbb{A}^*$ are identical.

*Part* 6: The matrix $\mathbb{A}^*$ is introduced by virtue of the fact that the formal $L^2$-adjoint of the operator $*d_\mathbb{A}$ acting on $C^\infty(Y; \text{ad}(P^\infty)_\mathbb{C} \otimes T^*Y)$ is $*d_{\mathbb{A}^*}$ which is not $*d_\mathbb{A}$ unless $\mathbb{A}$ is an SU(2) connection on P (its Hermitian part vanishes). A manifestation to keep in mind is this: The image in $C^\infty(Y; \text{ad}(P^\infty)_\mathbb{C} \otimes T^*Y)$ of the operator $*d_\mathbb{A}$ is $L^2$ orthogonal to the image in $C^\infty(Y; \text{ad}(P^\infty)_\mathbb{C} \otimes T^*Y)$ of $d_{\mathbb{A}^*}$ which is not the image of $d_\mathbb{A}$ unless $\mathbb{A}$ is an SU(2) connection on P.

By the same token, the operator $d_\mathbb{A}^\dagger$ when written using an orthonormal frame for $T^*Y$ sends $\mathfrak{v} \in C^\infty(Y; \text{ad}(P^\infty)_\mathbb{C} \otimes T^*Y)$ to the section $-(\nabla_{\mathbb{A}^*_i} \mathfrak{v})_i$ of $\text{ad}(P^\infty)_\mathbb{C}$. This is to say that the covariant derivatives that appear in $d_\mathbb{A}^\dagger$ are $\mathbb{A}^*$-covariant derivatives.

*Part* 7: To put (III.3.11) and (III.3.12) into the context of Lemma III.3.2, let $\mathbb{A}$ denote a given $SL(2;\mathbb{C})$ connection on $P^\infty \times_{SU(2)} SL(2;\mathbb{C})$. Write this connection as $A + i\mathfrak{a}$ with A denoting a connection on $P^\infty$ and with $\mathfrak{a}$ denoting an $\text{ad}(P^\infty)$ valued section of $T^*Y$ over this same domain. Meanwhile, write an element in $C^\infty(Y; \text{ad}(P^\infty)_\mathbb{C} \otimes T^*Y)$ as $\mathfrak{b} + i\mathfrak{c}$ with $\mathfrak{b}$ and $\mathfrak{c}$ denoting $\text{ad}(P^\infty)$-valued 1-forms on Y. And, write an element in $C^\infty(Y; \text{ad}(P^\infty)_\mathbb{C})$ as $\mathfrak{c}_t + i\mathfrak{b}_t$ with $\mathfrak{c}_t$ and $\mathfrak{b}_t$ denoting $\text{ad}(P^\infty)$-valued functions on Y. Then,

- $d_\mathbb{A}(\mathfrak{b} + i\mathfrak{c}) = (d_A \mathfrak{b} - \mathfrak{c} \wedge \mathfrak{a} - \mathfrak{a} \wedge \mathfrak{c}, d_A \mathfrak{c} + \mathfrak{b} \wedge \mathfrak{a} + \mathfrak{a} \wedge \mathfrak{b})$.
- $d_\mathbb{A}^\dagger(\mathfrak{b} + i\mathfrak{c}) = (-\nabla_{Ai}\mathfrak{b}_i - [\mathfrak{a}_i, \mathfrak{c}_i], -\nabla_{Ai}\mathfrak{c}_i + [\mathfrak{a}_i, \mathfrak{b}_i])$.
- $d_{\mathbb{A}^*}(\mathfrak{b}_t + i\mathfrak{c}_t) = (d_A \mathfrak{c}_t + [\mathfrak{a}, \mathfrak{b}_t], d_A \mathfrak{b}_t - [\mathfrak{a}, \mathfrak{c}_t])$.

(III.3.13)

Comparing these identities with (I.1.2) identifies the $\gamma_i \nabla_i + \rho_i[\mathfrak{a}_i, \cdot]$ part of $\mathcal{D}$ on any constant t slice of $(0,\infty) \times Y$ (viewed as acting on sections of $\mathbb{W}$ over $\{t\} \times Y$) with the operator on

$$C^\infty(Y; \text{ad}(P^\infty)_\mathbb{C} \otimes T^*Y) \oplus C^\infty(Y; \text{ad}(P^\infty)_\mathbb{C})$$

(III.3.14)

that sends a pair $(\eta, \nu)$ to the $\mathbb{A} = A + i\mathfrak{a}$ version of

$$-(*d_\mathbb{A} \eta + d_{\mathbb{A}^*} \nu, d_\mathbb{A}^\dagger \eta).$$

(III.3.15)

Note in particular the appearance of $\mathbb{A}^*$ in this formula. (Remember that $\mathbb{W}$ is the vector bundle $\oplus_2(\mathrm{ad}(P)) \oplus (T^*Y \otimes \mathrm{ad}(P))$.)

*Part 8*: If $\mathbb{A}$ is a flat connection, then the square of the $L^2$ norm of what is depicted in (III.3.15) can be written using an integration by parts as

$$\int_Y (|d_\mathbb{A}\eta|^2 + |d_\mathbb{A}^\dagger \eta|^2 + |d_{\mathbb{A}^*}\nu|^2) \;.$$

(III.3.16)

With regards to the derivation: The integral of the Hermitian inner product between $*d_\mathbb{A}\eta$ and $d_{\mathbb{A}^*}\nu$ is zero when $d_\mathbb{A}^2 = 0$. Because these are mutually orthogonal, the square of the $L^2$ norm of their sum is the sum of the squares of their $L^2$ norms.

One other important observation: If $d_A *\mathfrak{a} = 0$, then an integration by parts writes the $|d_{\mathbb{A}^*}\nu|^2$ integral in ((III.3.16) as

$$\int_Y |d_{\mathbb{A}^*}\nu|^2 = \int_Y (|\nabla_A \nu|^2 + |[\mathfrak{a},\nu]|^2) \;.$$

((III.3.17)

This last identity holds in particular when $(A,\mathfrak{a})$ is the pair $(A^\infty, \mathfrak{a}^\infty)$ from Lemma III.3.2 because of the assumption in that lemma that $d_{A^\infty} * \mathfrak{a}^\infty$.

*Part 9*: With regards to the rest of ((III.3.16): There is a Bochner-Wietzenboch formula for $d_\mathbb{A} + d_\mathbb{A}^\dagger$ acting on $C^\infty(Y; \mathrm{ad}(P^\infty)_\mathbb{C} \otimes T^*Y)$ which is this:

$$\int_Y (|d_\mathbb{A}\eta|^2 + |d_\mathbb{A}^\dagger \eta|^2) = \int_Y (|\nabla_A \eta|^2 + |[\mathfrak{a},\eta]|^2) + \int_Y \langle \eta, \mathbb{F}\eta \rangle$$

((III.3.18)

with $\mathbb{F}$ denoting a sum of an endomorphism determined by $\mathbb{A}$ and one that is determined by the Ricci tensor of the metric on $Y$.

For the purposes at hand, it is enough to know that $|\mathbb{F}| \leq c$ when defined by a given flat, $SL(2;\mathbb{C})$ connection. Then, (III.3.11) and (III.3.18) lead to the inequality

$$\int_Y (|d_\mathbb{A}\eta|^2 + |d_\mathbb{A}^\dagger \eta|^2) \geq c^{-1} \int_Y (|\nabla_A \eta|^2 + |[\mathfrak{a},\eta]|^2) \;,$$

(III.3.19)

which is essential for what comes next.

*Part 10*: By way of a summary of Parts 7-9: If $\mathbb{A} = A + i\mathfrak{a}$ is a regular flat connection with $d_A * \mathfrak{a} = 0$, and if $\psi$ is a section of $\mathbb{W}$, then

$$\int_Y |\gamma_i \nabla_{Ai}\psi + \rho_i[\mathfrak{a}_i,\psi]|^2 \geq c^{-1} \int_Y (|\nabla_A \psi|^2 + |[\mathfrak{a},\psi]|^2)$$

(III.3.20)

which is suggestive of what is asserted by Lemma III.3.2.

*Part 11*: Let $\mathbb{A}$ again denote a connection on $P^\infty \times_{SU(2)} SL(2;\mathbb{C})$. Use pull-back by the projection map $\pi$ to view $\mathbb{A}$ as an $SL(2;\mathbb{C})$ connection over $(0,\infty) \times Y$ and write the latter now as $A^\infty + i\mathfrak{a}^\infty$ with $A^\infty$ being a connection on the $\pi$-pull-back of $P^\infty$ bundle over $(0,\infty) \times Y$ and with $\mathfrak{a}^\infty$ being section of the tensor product with $T^*Y$ of the corresponding Lie-algebra bundle. Use $\nabla^\infty$ to denote the $A^\infty$-covariant derivative.

Because $A^\infty$ and $\mathfrak{a}^\infty$ are pulled back by the projection map, the $A^\infty$-covariant derivative along the $(0,\infty)$ factor (thus, $\nabla^\infty_t$) commutes with the rest of $\mathcal{D}^\infty$ which is the $\gamma_i \nabla^\infty_i + \rho_i[\mathfrak{a}^\infty_i,\cdot]$ part. This fact and the fact that $\gamma^i \nabla^\infty_i + \rho^i[\mathfrak{a}^\infty_i,\cdot]$ is a symmetric operator have the following consequence: If $\psi$ is an element in the $(A^\infty, \mathfrak{a}^\infty)$ version of $\mathbb{H}$, then

$$\|\mathcal{D}^\infty \psi\|_{\mathbb{L}}^2 = \|\nabla^\infty_t \psi\|_{\mathbb{L}}^2 + \|(\gamma^i \nabla^\infty_i + \rho_i[\mathfrak{a}^\infty_i,\cdot])\psi\|_{\mathbb{L}}^2 .$$

(III.3.21)

(The identity in (III.3.21) holds whenever the pair $(A^\infty, \mathfrak{a}^\infty)$ can be identified via a principal bundle isomorphism with a pair pulled back by $\pi$ from Y. There is no requirement here that $A^\infty + i\mathfrak{a}^\infty$ be flat.)

This last identity with (III.3.20) invoked on each constant $\{t\}$ slice $\{t\} \times Y$ using $(A^\infty, \mathfrak{a}^\infty)$ leads directly to the assertion in Lemma III.3.2.

## 4. Closed range and finite dimensional kernel

This section proves a part of Theorem III.1 by showing that $\mathcal{D}$ has closed range and finite dimensional kernel. The argument has three parts.

*Part 1*: If $\psi \in \mathbb{H}$, then (III.2.5) holds and so does Lemma III.2.1 if $\delta < c^{-1}$. If this last version of $c$ is large, then (III.2.5) and Lemma III.2.1 lead to this:

$$\int_{(0,\infty) \times Y} |\mathcal{D}\psi|^2 \geq c^{-1} \|\psi_0\|_{\mathbb{H}}^2 + c^{-1} \int_{(0,\infty) \times Y} |\mathcal{D}\psi_1|^2 - c \frac{1}{|\ln T|^2} \int_{(\frac{1}{T^2},\infty) \times Y} \frac{1}{t^2} |\psi_1|^2 .$$

(III.4.1)

Granted the preceding, invoke (III.3.3) with Lemma III.3.1 to see that (III.4.1) leads to the next inequality if T is larger than some large version of $c$.

$$\int_{(0,\infty)\times Y} |\mathcal{D}\psi|^2 \geq c^{-1}(\|\psi_0\|_{\mathbb{H}}^2 + \|\psi_{1\infty}\|_{\mathbb{H}}^2) + c^{-1}\int_{(0,\infty)\times Y} |\mathcal{D}\psi_{11}|^2 - c\frac{1}{|\ln T|^2}\int_{[\frac{1}{T^2},T^2]\times Y} \frac{1}{t^2}|\psi_{11}|^2 \ .$$

(III.4.2)

The important point with regards to this inequality is that $\psi_{11}$ has compact support in the fixed region where $\delta^2 \leq t \leq T^2$. Keep this in mind.

  The operator $\mathcal{D}$ has a Bochner-Weitzenboch formula which has the form of that depicted in (I.2.5) but with $\mathbb{X}$ now denoting an endomorphism that is linear in the components of the curvature of A, the covariant derivative of $\mathfrak{a}$, $\mathfrak{a} \wedge \mathfrak{a}$ and the Riemann curvature tensor on Y. For the present purposes, it is enough to know that $|\mathbb{X}| \leq c_{\delta,T}$ on the support of $\psi_{11}$ with $c_{\delta,T}$ denoting here and in what follows a number greater than 1 that is independent of the choice of $\psi$ from $\mathbb{H}$. (Its value can be assumed to increase between successive appearances.) Using this Bochner-Weitzenboch formual with (III.4.5) leads to the following inequality:

$$\int_{(0,\infty)\times Y} |\mathcal{D}\psi|^2 \geq c^{-1}(\|\psi_0\|_{\mathbb{H}}^2 + \|\psi_{1\infty}\|_{\mathbb{H}}^2 + \|\psi_{11}\|_{\mathbb{H}}^2) - c_{\delta,T}\|\psi_{11}\|_{\mathbb{L}}^2 \ .$$

(III.4.3)

  There is one last point to be made: The sum of the three $\|\cdot\|_{\mathbb{H}}$ norms can be replaced by $\|\psi\|_{\mathbb{H}}$ if the $\mathbb{L}$ integral is replaced by an integral of $\psi$ over a larger but still compact domain:

$$\int_{(0,\infty)\times Y} |\mathcal{D}\psi|^2 \geq c^{-1}\|\psi\|_{\mathbb{H}}^2 - c_{\delta,T}\int_{[\frac{1}{2}\delta^2,2T^2]\times Y} |\psi|^2 \ .$$

(III.4.4)

This is because the covariant derivative of the product $f\psi$ with $f$ being a function (for example, the cut-off functions that are used to define $\psi_0$, $\psi_{11}$ and $\psi_{1\infty}$) is the sum of $f\nabla\psi$ and $df \otimes \psi$.

  *Part 2*: This part proves that the kernel of $\mathcal{D}$ is finite dimensional. To this end, note first that if $\mathcal{D}\psi = 0$ and $\psi$ is not identically zero then it must be non-zero at some points where t is between $\frac{1}{2}\delta^2$ and $2T^2$. This follows from (III.4.4).

  With preceding in mind, suppose to the contrary that the kernel of $\mathcal{D}$ is infinite dimensional so as to generate nonsense In this instance, there is a set $\{\psi_n\}_{n\in\mathbb{N}}$ in the kernel of $\mathcal{D}$ whose restriction to $[\frac{1}{2}\delta^2, 2T^2]\times Y$ is orthonormal in the following sense:

- $$\int_{[\frac{1}{2}\delta^2, 2T^2]\times Y} \langle \psi_n, \psi_m \rangle = 0 \text{ if } n \neq m.$$
- $$\int_{[\frac{1}{2}\delta^2, 2T^2]\times Y} |\psi_n|^2 = 1.$$

(III.4.5)

Now it follows from (III.4.4) that this $\{\psi_n\}_{n\in\mathbb{N}}$ sequence has uniformly bounded $\|\cdot\|_{\mathbb{H}}$-norm. Therefore, it has a weakly convergent subsequence in $\mathbb{H}$. Let $\Lambda \subset \mathbb{N}$ denote the labels of this subsequence and let $\psi_\infty$ denote the weak limit of the subsequence (an element in $\mathbb{H}$). By virtue of the Rellich theorem (see for example [F]) applied to the compact domain $[\frac{1}{2}\delta^2, 2T^2]\times Y$, the sequence $\{\psi_n\}_{n\in\Lambda}$ converges strongly in the $\mathbb{L}$-norm topology (the $L^2$ topology) on the space of sections of $\mathbb{W}$ over this domain. This strong convergence to the limit on $[\frac{1}{2}\delta^2, 2T^2]\times Y$ is nonsensical because it runs afoul of the conditions in (III.4.5). (The conditions in (III.4.5) imply that the integral over $[\frac{1}{2}\delta^2, 2T^2]\times Y$ of $|\psi_n - \psi_m|^2$ is equal to 2 when ever $n \neq m$, whereas: Convergence in the $\mathbb{L}$ topology says (by definition) that these integrals must limit to zero as the indices n and m from $\Lambda$ get ever larger.)

*Part 3*: This part proves that $\mathcal{D}$ has closed range. This is again a reductio ad absurdum argument. To start, let $\{\mathfrak{g}_n\}_{n\in\mathbb{N}}$ denote a convergent sequence in $\mathbb{L}$, each in the image of $\mathcal{D}$, with limit $\mathfrak{g}$. Suppose for the sake of argument that $\mathfrak{g}$ is not in the image of $\mathcal{D}$ to generate nonsense.

For each positive integer n, there is by assumption an element $\psi_n$ from $\mathbb{H}$ obeying $\mathcal{D}\psi_n = \mathfrak{g}_n$. An element from the kernel of $\mathcal{D}$ in $\mathbb{H}$ can be added to $\psi_n$ if needs be so that

$$\int_{[\frac{1}{2T^2}, 2T^2]\times Y} \langle \psi_n, \psi \rangle = 0$$

(III.4.6)

when $\psi$ is from the kernel of $\mathcal{D}$. This last condition on $\psi_n$ is assumed.

Suppose first that there is a subsequence of $\{\psi_n\}_\Lambda$ with an a priori bound on

$$\int_{[\frac{1}{2}\delta^2, 2T^2]\times Y} |\psi_n|^2.$$

(III.4.7)

By virtue of (III.4.4), that subsequence will have an apriori bound on its $\|\cdot\|_{\mathbb{H}}$ norm also. As a consequence, there will be a subsequence that converges weakly in $\mathbb{H}$. The limit of the latter subsequence will then be mapped by $\mathcal{D}$ to $\mathfrak{g}$, this because of the weak convergence: Let $\Lambda \subset \mathbb{N}$ denote the labels for the weakly convergent sequence and let $\psi$ denote the weak limit. Then, for any $\mathfrak{w} \in \mathbb{L}$, the sequence indexed by $\Lambda$ with n'th term

$$\int_{(0,\infty)\times Y} \langle (\mathcal{D}\psi_n - \mathcal{D}\psi), \mathfrak{w} \rangle$$

(III.4.8)

converges to zero (by virtue of the weak convergence). Therefore $\mathcal{D}\psi$ and $\mathfrak{g}$ have the same inner product with all elements in $\mathbb{L}$ which implies that $\mathcal{D}\psi = \mathfrak{g}$.

Granted the assumptions about $\{\psi_n\}_{n\in\mathbb{N}}$, the preceding conclusion is avoided only in the event that the sequence of integrals in (III.4.7) is increasing and unbounded. Let $\{N_n\}_{n\in\mathbb{N}} \subset (0, \infty)$ denote this sequence and define a new sequence in $\mathbb{H}$ by the rule whereby the n'th term is $\frac{1}{\sqrt{N_n}} \psi_n$. (This n'th element is denoted by $\psi'_n$ in what follows.) It then follows from (III.4.4) that this sequence $\{\psi'_n\}_{n\in\mathbb{N}}$ has bounded $\|\cdot\|_\mathbb{H}$ norm; and so it has a weakly convergent subsequence in $\mathbb{H}$ whose limit is denoted by $\psi'$. This element is non-zero because the Rellich lemma guarantees that the integral of $|\psi'|^2$ over the domain $[\frac{1}{2}\delta^2, 2T^2]\times Y$ is equal to 1. It is also in the kernel of $\mathcal{D}$ because $\mathcal{D}\psi'_n$ is $\frac{1}{\sqrt{N_n}} \mathfrak{g}_n$ which converges to zero in $\mathbb{L}$. This last conclusion is nonsensical because each $\psi'_n$ obeys (III.4.6) and therefore so does $\psi'$ (by virtue of the weak convergence again).

## 5. The cokernel of $\mathcal{D}$ and the $\mathbb{L}$-kernel of $\mathcal{D}$

The cokernel of $\mathcal{D}$ is canonically isomorphic to the $\mathbb{L}$-kernel of $\mathcal{D}^\dagger$. The goal is to prove that the kernel of the latter is finite dimensional. This is done by proving that the $\mathbb{L}$-kernel of $\mathcal{D}^\dagger$ is canonically isomorphic to the $\mathbb{L}$-kernel of $\mathcal{D}$ and that every element of the latter is in $\mathbb{H}$.

The isomorphism between the two $\mathbb{L}$-kernels is defined as follows: Let ¥ denote the automorphism of $\mathbb{W}$ that when viewed as $\oplus_2(\text{ad}(P)\otimes T^*((0,\infty)\times Y))$ sends any given element $(\mathfrak{b}_t dt + \mathfrak{b}, \mathfrak{c}_t dt + \mathfrak{c})$ to

$$(\mathfrak{c}_t dt - \mathfrak{c}, -\mathfrak{b}_t dt + \mathfrak{b})$$

(III.5.1)

This automorphism obeys $¥^2 = -1$. More to the point, it follows from (I.1.4) that

$$\mathcal{D}¥ = -¥\mathcal{D}^\dagger .$$

(III.5.2)

There is also a depiction of ¥ using the Clifford matrix notation in (1.8): $¥ = \gamma_1\gamma_2\gamma_3\rho_1\rho_2\rho_3$. Because this algebraic automorphism ¥ maps $\mathbb{H}$ to $\mathbb{H}$ and $\mathbb{L}$ to $\mathbb{L}$, the $\mathbb{H}$-kernel of $\mathcal{D}^\dagger$ is mapped isometrically by ¥ onto the $\mathbb{H}$-kernel of $\mathcal{D}$ and likewise ¥ maps the $\mathbb{L}$-kernel of $\mathcal{D}^\dagger$ isometrically onto the $\mathbb{L}$ kernel of $\mathcal{D}$.

As a consequence of what was just observed, the conclusion that $\mathcal{D}$ has finite dimensional cokernel follows with a proof that the $\mathbb{L}$-kernel of $\mathcal{D}$ is in $\mathbb{H}$.

### a) The $\mathbb{H}$-kernel of $\mathcal{D}$ and the $\mathbb{L}$-kernel of $\mathcal{D}$

The next lemma plays a central role in that proof that the $\mathbb{L}$-kernel of $\mathcal{D}$ is a subspace of the $\mathbb{H}$-kernel of $\mathcal{D}$.

**Lemma III.5.1**: *A smooth element from $C^\infty((0,\infty)\times Y;\mathbb{W})$ is in $\mathbb{H}$ if both of the following conditions are met:*

- $\mathcal{D}\eta$ *is in $\mathbb{L}$ ,*

- *The function on $(0,\infty)$ defined by the rule* $t \to \frac{1}{t^2} \int_t^{2t} \int_{\{\cdot\}\times Y} |\eta|^2$ *is bounded.*

Looking ahead to the application of this lemma: If $\eta$ is from the $\mathbb{L}$-kernel of $\mathcal{D}$, then the condition in the top bullet of this lemma is obeyed (!). Sections III.5b-III.5d establish the condition in the second bullet. In this regard, if $\eta$ is from $\mathbb{L}$, then the only issue with regards to an upper bound for the second bullet's function is where $t \ll 1$. To say more about this issue (and for use later), introduce by way of notation $\mathcal{N}_r$ to denote the radius r tubular neighborhood of K with it understood without saying that r is chosen from $(0, r_0]$. The upcoming Lemma III.5.2 in Section III.5c asserts that the integral of $\frac{1}{t^2}|\eta|^2$ is finite on any domain of the form $(0, s] \times (Y-\mathcal{N}_{16Rs})$ for a suitable, s-independent choice of R if $\eta$ is in the $\mathbb{L}$-kernel of $\mathcal{D}$. Section III.5d proves that the same function $\frac{1}{t^2}|\eta|^2$ has finite integral on $(0, \delta] \times \mathcal{N}_r$ for suitable choices of $\delta$ and r.

*Proof of Lemma III.5.1*: The first observation is this: If $t_1 > t_0 > 0$ are any given pair of times, then the $[t_0, t_1] \times Y$ integral of $|\nabla\eta|^2 + |[\mathfrak{a},\eta]|^2$ is finite. This is because $\eta$ is assumed to be smooth. Thus, a first issue is whether the function $|\nabla\eta|^2 + |[\mathfrak{a},\eta]|^2$ has finite integral where t is unbounded and also where t has no positive lower bound. Part 1 of the following considers the integral where t has no upper bound and Part 2 considers the integral where t has no positive lower bound. Part 3 proves that if the $(0,\infty)\times Y$ integral of $|\nabla\eta|^2 + |[\mathfrak{a},\eta]|^2$ is bounded, then $\eta$ is necessarily in $\mathbb{H}$.

*Part 1*: Fix times $T > 2T_0 > 2$ and let $X_T$ denote the function on $(0, \infty)$ that is given by the rule $t \to \chi(\frac{t}{T} - 1)\chi(2 - \frac{t}{T_0})$. This function is equal to 1 where $t$ is between $2T_0$ and $T$; and it is equal to zero where $t$ is less than $\frac{1}{2}T_0$ or greater than $2T$. Now write

$$\mathcal{D}(X_T \eta) = (\tfrac{d}{dt} X_T)\eta + X_T \mathcal{D}\eta \ .$$

(III.5.1)

The $(0,\infty) \times Y$ integral of the square of the norm of each term on the right hand side of (III.5.1) has a T-independent upper bound. This is because of the assumption from the lemma's first bullet for the $X_T \mathcal{D}\eta$ term and because of the assumption from the lemma's second bullet for the $(\tfrac{d}{dt} X_T)\eta$ term. As a consequence, the $(0,\infty) \times Y$ integral of the square of the norm of $\mathcal{D}(X_T \eta)$ has a T-independent upper bound.

With respect to CONSTRAINT 3: Write the $\pi$-pull back of CONSTRAINT 3's flat SL(2;$\mathbb{C}$) connection as $A^\infty + i\mathfrak{a}^\infty$ after using the isomorphism in Item d) of that constraint to identify P with $\pi^* P^\infty$ (the notation has $A^\infty$ denoting a connection on P and $\mathfrak{a}^\infty$ denoting a section of $ad(P) \otimes T^*Y$). Fix $\varepsilon > 0$ for the moment. If $T_0$ is sufficiently large (which will be assumed henceforth), then $(A, \mathfrak{a})$ differs from $(A^\infty, \mathfrak{a}^\infty)$ by a term with norm bounded by $\frac{\varepsilon}{t}$ where $t > \frac{1}{4}T_0$. With this fact in mind, write $\mathcal{D}$ as $\mathcal{D}^\infty + \mathfrak{e}$ where $t > \frac{1}{4}T_0$ with $\mathcal{D}^\infty$ denotes the $(A^\infty, \mathfrak{a}^\infty)$ version of $\mathcal{D}$ and where $\mathfrak{e}$ is an endomorphism with norm bounded by $c\frac{\varepsilon}{t}$. Then write $\mathcal{D}(X_T \eta)$ as

$$\mathcal{D}^\infty(X_T \eta) = \mathcal{D}(X_T \eta) - \mathfrak{e}(X_T \eta) \ .$$

(III.5.2)

The $(0,\infty) \times Y$ integral of the square of the norm of $\mathcal{D}(X_T \eta)$ is bounded by virtue of what is said in the preceding paragraph. As for the $\mathfrak{e}(X_T \eta)$ term, the $(0,\infty) \times Y$ integral of the square of its norm is at most the $[\frac{1}{2} T_0, 2T] \times Y$ integral of $c\frac{\varepsilon^2}{t^2}|\eta|^2$.

Meanwhile, by virtue of (III.3.20) and (III.3.21), and by virtue of Hardy's inequality (which is (I.2.3) in this instance),

$$\|\mathcal{D}^\infty(X_T \eta)\|_{\mathbb{L}}^2 \geq c^{-1}(\|\nabla^\infty(X_T \eta)\|_{\mathbb{L}}^2 + \|[\mathfrak{a}^\infty, X_T \eta]\|_{\mathbb{L}}^2) + c^{-1} \|\tfrac{1}{t} X_T \eta\|_{\mathbb{L}}^2 \ .$$

(III.5.3)

This bound with the identity in (III.5.2) and the afore-mentioned $c\varepsilon^2 \|\tfrac{1}{t}\eta\|_{\mathbb{L}}^2$ bound for $\|\mathfrak{e}(X_T \eta)\|_{\mathbb{L}}^2$ imply that

$$\|\mathcal{D}(X_T \eta)\|_{\mathbb{L}}^2 \geq c^{-1} \|\nabla^\infty(X_T \eta)\|_{\mathbb{L}}^2 + \|[\mathfrak{a}^\infty, X_T \eta]\|_{\mathbb{L}}^2 + (c^{-1} - c\varepsilon^2)\|\tfrac{1}{t} X_T \eta\|_{\mathbb{L}}^2 \ ,$$

(III.5.4)

which implies in turn that

$$\|\mathcal{D}(\chi_T\eta)\|_{\mathbb{L}}{}^2 \geq c^{-1} \|\nabla(\chi_T\eta)\|_{\mathbb{L}}{}^2 + \|[\mathfrak{a},\chi_T\eta]\|_{\mathbb{L}}{}^2 + (c^{-1} - c\varepsilon^2)\|\tfrac{1}{t}\chi_T\eta\|_{\mathbb{L}}{}^2$$

(III.5.5)

because $\nabla - \nabla^\infty$ and $[\mathfrak{a} - \mathfrak{a}^\infty, \cdot]$ are endomorphisms with norms bounded by $c\tfrac{\varepsilon}{t}$.

If $\varepsilon \leq c^{-1}$ (as is henceforth assumed) then (III.5.5) and what was said about (III.5.1) lead to a T-independent upper bound for the $[T_0, T] \times Y$ integral of $|\nabla\eta|^2 + |[\mathfrak{a},\eta]|^2$. As a consequence, the function $|\nabla\eta|^2 + |[\mathfrak{a},\eta]|^2$ has finite integral on $[t, \infty) \times Y$ as long as $t > 0$.

*Part 2*: Fix a small, positive number $\varepsilon$ to be much less than 1. Having done that, define $R_\varepsilon$ and $t_\varepsilon$ as in CONSTRAINT 1. Then fix $T > \tfrac{1}{t_\varepsilon}$, a positive number $\delta \in (0, \tfrac{1}{100T^2})$ and let $\varsigma_\delta$ denote the function on $(0, \infty)$ that is given by the rule $t \to \chi(2 - \tfrac{t}{\delta})\chi(2 + \tfrac{\ln t}{\ln T})$. This function is zero where $t < \delta$ and where $t > \tfrac{1}{T}$; and it equals one where $2\delta < t < \tfrac{1}{T^2}$. Set $\psi_0$ to denote $\varsigma_\delta \eta$. This is an element in $\mathbb{H}$ with compact support on $(\delta, t_\varepsilon) \times Y$. The element $\mathcal{D}\psi_0$ can be written as

$$\mathcal{D}\psi_0 = (\tfrac{d}{dt}\varsigma_\delta)\eta + \varsigma_\delta \mathcal{D}\eta .$$

(III.5.6)

The important point now is that the $(0, t_\varepsilon) \times Y$ integral of the square of the norm of each term on the right hand side of (III.5.6) has a $\delta$-independent upper bound. With regards to $\mathcal{D}\eta$, this is the assumption in the first bullet of Lemma III.5.1  Regarding the term with $(\tfrac{d}{dt}\varsigma_\delta)\eta$: This term has support only for t values obeying $\tfrac{1}{2}\delta < t < 2\delta$ and $\tfrac{1}{T^2} < t < \tfrac{1}{T}$. When t is in the latter region, the norm of $(\tfrac{d}{dt}\varsigma_\delta)\eta$ is independent of $\delta$; and when t is in the region between $\tfrac{1}{2}\delta$ and $2\delta$, the norm of $(\tfrac{d}{dt}\varsigma_\delta)\eta$ is bounded by $c\tfrac{1}{t}|\eta|$ and so the integral of its square is bounded by $c$ due to Lemma III.5.1's second bullet assumption.

Granted the preceding, invoke Lemma III.3.1 to see that if $T \geq c$, then there is a $\delta$-independent bound on the $\mathbb{H}$-norm of $\psi_0$. This implies in turn a $\delta$-independent bound on

$$\int_{(2\delta, \tfrac{1}{T^2}) \times Y} (|\nabla\eta|^2 + |[\mathfrak{a},\eta]|^2) ,$$

(III.5.7)

which implies in turn that the $(0, \tfrac{1}{T^2}) \times Y$ integral of $|\nabla\eta|^2 + |[\mathfrak{a},\eta]|^2$ is finite.

*Part 3*: This part of the proof explains why the following claim is true:

*If $\eta \in C^\infty((0,\infty) \times Y; \mathbb{W})$ and if the $(0,\infty) \times Y$ integral*
*of $|\nabla\eta|^2 + |[\mathfrak{a},\eta]|^2$ is finite, then $\eta \in \mathbb{H}$.*

(III.5.8)

The issue here is whether there is a sequence of compactly supported elements in $C^\infty((0,\infty)\times Y; \mathbb{W})$ that converge to $\eta$ with respect to the $\|\cdot\|_{\mathbb{H}}$-norm. As explained directly, such a sequence exists if the function of t given by the rule in the second bullet of Lemma III.5..1 has lim-inf equal to zero as $t \to \infty$ and also as $t \to 0$. Indeed, if that is so, then there is an unbounded, increasing sequence $\{T_n\}_{n\in\mathbb{N}}$ in $(1,\infty)$ and a decreasing sequence $\{\delta_n\}_{n\in\mathbb{N}} \subset (0,1)$ limiting to zero with the following property: The sequence in $\mathbb{H}$ whose n'th term, $\psi_n$, is $\chi(\frac{\delta_n}{t} - 1)\chi(\frac{t}{T_n} - 1)\eta$ is such that $\lim_{n\to\infty} \|\eta - \eta_n\|_{\mathbb{H}} = 0$. And, if this limit is zero, then $\eta$ is in $\mathbb{H}$ because $\mathbb{H}$ contains all of its $\mathbb{H}$-norm limit points.

The function of t given by the rule in the second bullet of Lemma III.5. has its lim-inf equal to zero as $t \to \infty$ and also as $t \to 0$ if the function $\frac{1}{t^2}|\eta|^2$ has finite integral on the domain $(0,\infty)\times Y$. To see if this integral is indeed finite, consider first the integral on $[1,\infty)\times Y$. Fix $T > 1$ and mimic the integration by parts proof in Section 2a of Hardy's inequality in (I.2.3) by writing

$$\int_{[1,T]} \tfrac{1}{t^2} f^2 \, dt + \tfrac{1}{T} f^2(T) = -2 \int_{[1,T]} \tfrac{1}{t} f(\tfrac{d}{dt} f) dt + f^2(1)$$

(III.5.9)

for any given function $f$ (in this case $|\eta|$) to get (via the Cauchy-Schwarz inequality for the right-hand integral) a T-independent upper bound on the $[1,T]\times Y$ integral of $\frac{1}{t^2}|\eta|^2$. (Keep in mind that $|d|\eta||$ is no greater than $|\nabla\eta|$.)

To see about a finite bound for the $(0,1]\times Y$ of $\frac{1}{t^2}|\eta|^2$, fix $\varepsilon \in (0, \frac{1}{1000})$ and then invoke CONSTRAINT 1 to find $R_\varepsilon > 1$ and $t_\varepsilon \in (0,1]$ such that $\mathfrak{a}$ and $-\frac{1}{2}\tau$ differ by at most $\frac{\varepsilon}{t}$ where $t < t_\varepsilon$ and $\mathrm{dist}(\cdot,K) > R_\varepsilon t$. As a consequence, if $\varepsilon < c^{-1}$, then $|[\mathfrak{a},\eta]|$ is greater than $\frac{1}{t^2}|\eta|^2$ where $t < t_\varepsilon$ and $\mathrm{dist}(\cdot,K) > Rt$ for $R \geq R_\varepsilon$.

As explained directly, a bound for the $t < t_\varepsilon$ and $\mathrm{dist}(\cdot,K) \leq Rt$ integral of $\frac{1}{t^2}|\eta|^2$ is a consequence of a standard *Poincare inequality*:

*There exists $\kappa > 0$ such that for any $\rho > 0$ and smooth function f where $|z| < 2\rho$ in $\mathbb{C}$,*

$$\tfrac{1}{\rho^2}\int_{|z|<\rho} f^2 \leq \kappa \left( \int_{|z|<2\rho} |df|^2 + \tfrac{1}{\rho^2} \int_{\rho<|z|<2\rho} f^2 \right)$$

(III.5.10)

This inequality is a consequence of the fact that the Dirichlet Laplacian on the $|z| < 2$ disk in $\mathbb{C}$ has only positive eigenvalues. Indeed, by rescaling z, one can replace $\rho$ by 1. Then invoke this positive lowest eigenvalue property for the function $\chi(|z|-1)f$. (The inequality can also be proved using yet another version of Hardy's inequality.)

To exploit (III.5.10), first use the diffeomorphism in CONSTRAINT 3 to identify $\mathcal{N}_K$ with $D_0\times S^1$. Having done that, apply (III.5.10) on each constant $(t,x_3)$ copy of $D_0$ (viewed as a subset of $\mathbb{C}$) for $t < t_\varepsilon$ and $x_3 \in S^1$ with $f = |\eta|$ and $\rho = cRt$. Then, integrate the inequality in (III.5.10) over $S^1$ and use the following facts:

- *The integral of $|d|\eta||^2$ on the part $\mathcal{N}_K$ where $\mathrm{dist}(\cdot, K) < 2cRt$ is no greater than $c$ times the integral of $|\nabla \eta|^2$ whole of $\{t\} \times Y$.*
- *The integral of $\frac{1}{t^2}|\eta|^2$ on the part $\mathcal{N}_K$ where $cRt < \mathrm{dist}(\cdot, K) < 2cRt$ is no greater than $c$ times the integral of $|[\mathfrak{a}, \eta]|^2$ over the whole of $\{t\} \times Y$.*

(III.5.11)

(The first follows because $|d|\eta|| \leq |\nabla \eta|$ and the second is a consequence of what was said in the preceding paragraph about $\mathfrak{a}$ and $-\frac{1}{2t}\tau$ differing by at most $\frac{\varepsilon}{t}$ where $t < t_\varepsilon$ and $\mathrm{dist}(\cdot, K) > R_\varepsilon t$.)

**b) The $\mathbb{L}$-kernel of $\mathcal{D}$ when $K = \emptyset$.**

The simplest case to consider is the case when there is no knot present. This is assumed here. To start the argument for this case, fix some positive but very small $\varepsilon$ and then use CONSTRAINT 1 to define a positive time $t_\varepsilon$ so $\mathfrak{a}$ differs from $-\frac{1}{2t}\tau$ by an $\mathrm{ad}(P)$ valued 1-form with norm at most $c\frac{\varepsilon}{t}$ and $|\nabla_A \tau| \leq c\frac{\varepsilon}{t}$ where $t < t_\varepsilon$. As noted in (3.3), the endomorphism $-\frac{1}{2t}\rho_i[\tau_i, \cdot]$ has eigenvalues $\pm\frac{\lambda}{t}$ with $\lambda$ either 1 or 2. This understood, write a given element from the $\mathbb{L}$-kernel of $\mathcal{D}$ (call it $\eta$) as $\eta = \sum_{\lambda=1,2}(\eta_\lambda^+ + \eta_\lambda^-)$ with $\eta_\lambda^+$ denoting the component of $\eta$ in the $+\frac{\lambda}{t}$ eigenspace of $-\frac{1}{2t}\rho_i[\tau_i, \cdot]$ and with $\eta_\lambda^-$ denoting the component in the $-\frac{\lambda}{t}$. The identity $\mathcal{D}\eta = 0$ when written where $t < t_\varepsilon$ using these components can be written as four separate equations (these are the $\mathcal{D}$ analogs of (II.3.7) with a non-flat metric):

- $\nabla_t \eta_\lambda^+ + \frac{\lambda}{t}\eta_\lambda^+ + \gamma_i \nabla_i \eta_\lambda^- + \mathfrak{r}_{\lambda+}\eta = 0$,
- $\nabla_t \eta_\lambda^- - \frac{\lambda}{t}\eta_\lambda^- + \gamma_i \nabla_i \eta_\lambda^+ + \mathfrak{r}_{\lambda-}\eta = 0$,

(III.5.12)

where $\{\mathfrak{r}_{\lambda\pm}\}_{\lambda=1,2}$ are endomorphisms with norms bounded by $c\frac{\varepsilon}{t}$.

Take the inner product of the top equation in (III.5.12) with $\eta_\lambda^+$ and the inner product of the lower equation with $\eta_\lambda^-$, then integrate the result over each $t \in (0, t_\varepsilon)$ slice $\{t\} \times Y$, then subtract the result of doing that to the lower bullet in (III.5.12) from the result of doing that to the top bullet in (III.5.12). The end product is an analog of (II.3.8):

$$\frac{d}{dt} \int_{\{t\}\times Y} (|\eta_\lambda^+|^2 - |\eta_\lambda^-|^2) = \frac{2\lambda}{t} \int_{\{t\}\times Y} (|\eta_\lambda^+|^2 + |\eta_\lambda^-|^2) + \mathcal{K}_\lambda(t)$$

(III.5.13)

where $\mathcal{K}_\lambda(\cdot)$ is a function on $(0, t_\varepsilon)$ that obeys the bound

$$|\mathcal{K}_\lambda(t)| \leq c \frac{\varepsilon}{t} \int_{\{t\} \times Y} |\eta|^2 .$$

(III.5.14)

Add the $\lambda = 1$ and $\lambda = 2$ versions of (III.5.14) and the result (with $\eta^\pm = \eta_1^\pm + \eta_2^\pm$) is this:

$$\frac{d}{dt} \int_{\{t\} \times Y} (|\eta^+|^2 - |\eta^-|^2) = \frac{2}{t} \int_{\{t\} \times Y} (|\eta_1^+|^2 + 2|\eta_2^+|^2 + |\eta_1^-|^2 + 2|\eta_2^-|^2) + \mathcal{K}(t)$$

(III.5.15)

with $\mathcal{K} = \mathcal{K}_1 + \mathcal{K}_2$.

Because of (III.5.14), the right hand side of (III.5.15) is positive when $\varepsilon < c^{-1}$ (unless $\eta \equiv 0$ at small t in which case $\eta$ is in $\mathbb{H}$). This implies that the integral

$$\int_{\{t\} \times Y} (|\eta^+|^2 - |\eta^-|^2)$$

(III.5.16)

is either positive on $(0, t_\varepsilon)$ or there exists $t_0 > 0$ such that (8.16) is negative for $t < t_0$. Suppose first that it is positive on $(0, t_0)$. Multiply both sides of (III.5.14) by $\frac{1}{t}$ to see that

$$\frac{d}{dt}\left(\frac{1}{t} \int_{\{t\} \times Y} (|\eta^+|^2 - |\eta^-|^2)\right) = \frac{1}{t^2} \int_{\{t\} \times Y} (|\eta_1^+|^2 + 3|\eta_2^+|^2 + 3|\eta_1^-|^2 + 5|\eta_2^-|^2) + \frac{1}{t}\mathcal{K}(t)$$

(III.5.17)

which implies that

$$\frac{d}{dt}\left(\frac{1}{t} \int_{\{t\} \times Y} (|\eta^+|^2 - |\eta^-|^2)\right) \geq (1 - c\varepsilon)\frac{1}{t^2} \int_{\{t\} \times Y} |\eta|^2 .$$

(III.5.18)

Now fix some small $\delta \in (0, t_\varepsilon)$ and integrate between $t = \delta$ and $t = t_\varepsilon$. Integrate by parts on the left and throw away the boundary term at $t = \delta$ (which is negative) to see that

$$\int_\delta^{t_\varepsilon} \frac{1}{t^2} \int_{\{t\} \times Y} |\eta|^2$$

(III.5.19)

has a $\delta$-independent upper bound. Therefore,

$$\int_{(0,t_\varepsilon) \times Y} \frac{1}{t^2}|\eta|^2$$

(III.5.20)

is finite. This implies that $\eta \in \mathbb{H}$ (see Lemma III.5.1).

Now suppose that there exists $t_0 \in (0, t_\varepsilon)$ such that (III.5.16) is negative for $t \in (0, t_0)$. In this case, (III.5.15) leads to this:

$$\frac{d}{dt} \int_{\{t\} \times Y} (|\eta^-|^2 - |\eta^+|^2) \geq -\frac{2-c\varepsilon}{t} \int_{\{t\} \times Y} (|\eta^-|^2 - |\eta^+|^2)$$

(III.5.21)

which implies that

$$\int_{\{t\} \times Y} (|\eta^-|^2 - |\eta^+|^2) \geq \frac{1}{t^{3/2}}$$

(III.5.22)

if $\varepsilon < c^{-1}$ and t is small. But now note that this bound runs afoul of the assumption that $\eta$ is from $\mathbb{L}$ since $|\eta|^2$ can't have finite integral on $(0, 1) \times Y$ in the event of (III.5.22).

### c) The $\mathbb{L}$-kernel of $\mathcal{D}$ on the complement of K

The lemma that follows describes the behavior of an element in the $\mathbb{L}$-kernel of $\mathcal{D}$ on the complement of small radius tubular neighborhoods of K.

To fix the notation: Having specified a small number (denoted by $\varepsilon$), this lemma refers to the numbers $t_\varepsilon$ and $R_\varepsilon$ that are supplied by CONSTRAINT 1. Also, given a positive number r, the lemma uses $\mathcal{N}_r$ to denote the subset of Y with distance at most r from K.

**Lemma III.5.2**: *Suppose that $\eta$ is from the $\mathbb{L}$-kernel of $\mathcal{D}$. There exists $\kappa > 1$ such that if $\varepsilon < \kappa^{-1}$ and $R \geq \kappa R_\varepsilon$ and $s < \frac{1}{4} t_\varepsilon$, then*

$$\int_{(0,s] \times (Y - \mathcal{N}_{16Rs})} \frac{1}{t^2} |\eta|^2 \leq \frac{\kappa}{R^2 s^2} \int_{(0,2s] \times (\mathcal{N}_{32Rs} - \mathcal{N}_{Rs})} |\eta|^2 + \frac{\kappa}{s^2} \int_{[s, 2s] \times (Y - \mathcal{N}_{Rs})} |\eta|^2$$

*Proof of Lemma III.5.2*: The proof has five parts.

*Part 1*: Fix a small, positive $\varepsilon$ much less than 1 and then use CONSTRAINT 1 to determine a corresponding $R_\varepsilon$ and $t_\varepsilon$. Fix $R > 100 R_\varepsilon$; and having fixed $t_0 \in (0, \frac{1}{100} t_\varepsilon)$ set r to be $R t_0$. No generality is lost by assuming that this number r is less than $\frac{1}{100}$ times the number $r_0$ which is the radius of the tubular neighborhood $\mathcal{N}_K$. With r in hand, reintroduce the function $\varpi_r$ from (III.2.6). Write an element $\eta$ from the $\mathbb{L}$-kernel of $\mathcal{D}$ where $t \leq t_\varepsilon$ as $\eta_{00} + \eta_{0K}$ with $\eta_{00} \equiv (1 - \varpi_r)\eta$ and $\eta_{0K} = \varpi_r \eta$. These obey

$$\mathcal{D}\eta_{00} = -\mathfrak{S}(d\varpi_r)(\eta_{00} + \eta_{0K}) \quad and \quad \mathcal{D}\eta_{0K} = \mathfrak{S}(d\varpi_r)(\eta_{00} + \eta_{0K})$$

(III.5.23)

where $\mathfrak{S}$ denotes the symbol of $\mathcal{D}$.

What is denoted by $\eta_{00}$ is supported on the complement of the radius r tubular neighborhood of K which is where $(A, \mathfrak{a})$ differs by at most $c\varepsilon \frac{1}{t}$ from $(\Gamma, -\frac{1}{2}\tau)$ if $t \leq t_0$. See CONSTRAINT 1 in this regard.. By virtue of this support restriction, $\eta_{00}$ can be written as $\eta^+ + \eta^-$ as done in the Subsection 8a, and the analysis that led to (III.5.15) now leads to this:

$$\frac{d}{dt} \int_{\{t\}\times Y} (|\eta^+|^2 - |\eta^-|^2) = \frac{2}{t} \int_{\{t\}\times Y} (|\eta_1^+|^2 + 2|\eta_2^+|^2 + |\eta_1^-|^2 + 2|\eta_2^-|^2) + \mathcal{K}_0(t) + \mathcal{K}_K(t)$$

(III.5.24)

with the norms of $\mathcal{K}_0$ and $\mathcal{K}_K$ obeying

$$|\mathcal{K}_0| \leq \frac{c(\varepsilon + \frac{1}{R})}{t} \int_{\{t\}\times Y} |\eta_{00}|^2 \quad and \quad |\mathcal{K}_K| \leq c\frac{1}{r} \int_{\{t\}\times Y} (1-\varpi_{r/2})^2 |\eta_{0K}|^2 .$$

(III.5.25)

With regards to (III.5.25): It is assumed henceforth that $\varepsilon$ is sufficiently small (meaning less than $c^{-1}$) and R is sufficiently large (meaning greater than $c$) so that the factor $c(\varepsilon + \frac{1}{R})$ that multiplies the left hand integral in (III.5.25) giving the bound for $|\mathcal{K}_0|$ is less than $\frac{1}{1000}$. Keep in mind that $t \leq t_0$ also.

*Part 2*: Suppose here that there is a decreasing sequence of times $\{t_n\}_{n\in\mathbb{N}} \subset (0, t_0]$ with limit 0 such that

$$\int_{\{t\}\times Y} (|\eta^+|^2 - |\eta^-|^2) \geq 0$$

(III.5.26)

when t is $t_n$ for each positive integer n. Granted this, integrate (III.5.24) from $t_n$ to $t_0$ for each $n \in \mathbb{N}$ to see that

$$\int_{(0, t_0]\times Y} \frac{1}{t}|\eta_{00}|^2 \leq c\frac{1}{r} \int_{(0, t_0]\times Y} (1-\varpi_{r/2})^2 |\eta_{0K}|^2 + \int_{\{t_0\}\times Y} |\eta_{00}|^2$$

(III.5.27)

*Part 3*: The other possibility is that there exists $t_1 \in (0, t_0]$ such that (III.5.26) fails at all $t < t$. To deal with this case, multiply both sides of (III.5.24) by t and, given $s \in (0, t_1]$ and $t_* \in (0, s)$, integrate both sides of the result from $t_*$ to $s$ to see (after an integration by parts) this:

$$\int_{[t_*,s]\times Y} |\eta_{00}|^2 \leq c \frac{s}{r} \int_{(0,t_0]\times Y} (1-\varpi_{r/2})^2 |\eta_{0K}|^2 + c t_* \int_{\{t_*\}\times Y} |\eta_{00}|^2$$

(III.5.28)

The $|\eta_{00}|^2$ integral on the right hand side needs to be dealt with. To this end:

*Supposing that $f$ is a non-negative function on $(0, 1]$ with finite integral, there is a sequence $\{t_n\}_{n\in\mathbb{N}} \subset (0, 1]$ with limit zero such that $\lim_{n\in\mathbb{N}} t_n f(t_n) = 0$.*

(III.5.29)

Indeed, if this were not the case, there would exist some positive $\delta$ and a positive $t_*$ such that $f(t) \geq \frac{\delta}{t}$ on $(0, t_*)$ which is nonsensical if $f$ has finite integral on $(0, 1]$.

Granted (III.5.29), then (III.5.28) with a judiciously chosen sequence of $t_*$'s limiting to zero leads to the next bound:

$$\int_{(0,t]\times Y} |\eta_{00}|^2 \leq c \frac{t}{r} \int_{(0,t_0]\times Y} (1-\varpi_{r/2})^2 |\eta_{0K}|^2$$

(III.5.30)

which holds for $t \leq t_1$. The latter bound implies in turn the following: For each $t \in (0, t_1]$, there exists $t' \in [\frac{1}{2} t, t]$ such that

$$\int_{\{t'\}\times Y} |\eta_{00}|^2 \leq c \frac{1}{r} \int_{(0,t_0]\times Y} (1-\varpi_{r/2})^2 |\eta_{0K}|^2 \ .$$

(III.5.31)

With this understood, fix a judiciously chosen, decreasing sequence $\{t_n\}_{n\in\mathbb{N}} \subset (0, t_0]$ with limit zero so that (III.5.31) holds with $t' = t_n$. Go back to (III.5.24) and integrate both sides from $t_n$ to $t_0$ and then integrate by parts. Having done that, take the $n \to \infty$ limit of the result to see that (III.5.27) holds in this case also.

*Part 4*: The plan now is to redo the analysis in Parts 2 and 3 with $r$ replaced by $4r$. This is done because the new version of $\eta_{0K}$ will be equal to the old version of $\eta_{00}$ where $(1-\varpi_{4r})$ is nonzero. Because of that, and due to (III.5.28), the new version of $\eta_{0K}$ obeys

$$\int_{(0,t_0]\times Y} \frac{1}{t}(1-\varpi_{4r})^2 |\eta_{0K}|^2 \leq c \frac{1}{r} \int_{(0,t_0]\times Y} (1-\varpi_{r/2})^2 \varpi_{4r}^2 |\eta_{0K}|^2 + c \int_{\{t_0\}\times Y} |\eta_{00}|^2 \ .$$

(III.5.32)

With this inequality understood, suppose once more that there exists a decreasing sequence of times $\{t_n\}_{n\in\mathbb{N}} \subset (0, t_\epsilon]$ with limit zero such that (III.5.26) holds at each such time. Multiply (III.5.24) by $\frac{1}{t}$ and then, having fixed a positive integer n, integrate both sides of the resulting identity on the domain $[t_n, t_\varepsilon]$. Integrate by parts on the left hand side, throw away the $t_n$ boundary term because it is negative, and the result with (III.5.25) and (III.5.32) leads to an inequality that says this:

$$\int_{[t_n, t_0] \times Y} \tfrac{1}{t^2} |\eta_{00}|^2 \le c \tfrac{1}{r^2} \int_{(0, t_0] \times Y} (1 - \varpi_{r/2})^2 \varpi_{4r}^2 |\eta_{0K}|^2 + c (\tfrac{1}{r} + \tfrac{1}{t_0}) \int_{\{t_0\} \times Y} |\eta_{00}|^2 .$$

(III.5.33)

Since the right hand side of this is independent of n and since $\lim_{n \to 0} t_n = 0$, the inequality in (III.5.33) leads to this a priori bound:

$$\int_{[t_n, t_0] \times Y} \tfrac{1}{t^2} |\eta_{00}|^2 \le c \tfrac{1}{r^2} \int_{(0, t_0] \times Y} (1 - \varpi_{r/2})^2 \varpi_{4r}^2 |\eta_{0K}|^2 + c (\tfrac{1}{r} + \tfrac{1}{t_0}) \int_{\{t_0\} \times Y} |\eta_{00}|^2 .$$

(III.5.34)

Now suppose to the contrary that there exists $t_1 \subset (0, t_0]$ such that the integral on the left hand side of (III.5.26) is negative at all times $t < t_1$. No generality is lost by assuming that either the left hand side of (III.5.26) is zero at $t = t_1$ or $t_1 = t_0$. Note that in either case, the $n = 1$ version of the inequality in (III.5.33) holds. For $t < t_1$, the identity in (III.5.24) with (III.5.25) leads to a differential inequality for $f \equiv \int_{\{t\} \times Y} (|\eta^-|^2 - |\eta^+|^2)$:

$$\tfrac{d}{dt} f < -\tfrac{1.9}{t} f + |\mathcal{K}_{0K}| .$$

(III.5.35)

Supposing that $s < t \le t_1$ and given (III.5.32), this leads (via integration) to the following:

$$f(s) \ge (\tfrac{t}{s})^{1.9} \big( f(t) - c \tfrac{t}{r^2} \int_{(0, t_0] \times Y} (1 - \varpi_{r/2})^2 \varpi_{4r}^2 |\eta_{0K}|^2 - c \tfrac{t}{r} \int_{\{t_0\} \times Y} |\eta_{00}|^2 \big) .$$

(III.5.36)

Because $\eta$ is from $\mathbb{L}$, the preceding inequality can hold only in the event that the right hand side is negative which requires that

$$f(t) \le c \tfrac{t}{r^2} \int_{(0, t_0] \times Y} (1 - \varpi_{r/2})^2 \varpi_{4r}^2 |\eta_{0K}|^2 + c \tfrac{t}{r} \int_{\{t_0\} \times Y} |\eta_{00}|^2$$

(III.5.37)

for any $t \le t_1$ (which is either $t_0$ or the smallest time in $(0, t_\varepsilon]$ where $f$ is zero).

With (III.5.37) in hand, multiply both sides of (III.5.24) by $\tfrac{1}{t}$ and integrate between any given $n \in \mathbb{N}$ version of $t_n = \tfrac{1}{n} t_1$ and $t_0$. Then, invoke (III.5.36) with $t = t_n$ and (III.5.25) and (III.5.32) to see that (III.5.33) holds. Since n can be any positive integer, so does (III.5.34).

*Part 5*: To complete the proof of Lemma III.5.2: The element $\eta_{00}$ that appears in (III.5.24) is equal to $(1 - \varpi_{4r})\eta$ and thus it is equal to $\eta$ where $\text{dist}(\cdot, K) \geq 8Rt_0$. With this understood, then the inequality asserted by the lemma follows by taking the average of the various $t_0 \in [s, 2s]$ versions of (III.5.24) with the average defined by integration over this interval in $\mathbb{R}$.

### d) The $\mathbb{L}$-kernel of $\mathcal{D}$ near K

The focus in this subsection is the behavior of an element (call it $\eta$ again) in the $\mathbb{L}$-kernel of $\mathcal{D}$ near K. There are three parts to the analysis.

*Part 1*: The following observation is needed for the later parts of the analysis.

*If $\eta$ is from the $\mathbb{L}$-kernel of $\mathcal{D}$, then $\zeta \equiv \text{dist}(\cdot, K)\eta$ is in $\mathbb{H}$.*

(III.5.38)

This is proved momentarily; what follows directly is a corollary:

*If $\eta$ is from the $\mathbb{L}$-kernel of $\mathcal{D}$, then the $(0,1] \times Y$ integral of $\text{dist}^2(\cdot, K)|\nabla\eta|^2$ is finite.*

(III.5.39)

What with (III.5.38), this follows from Lemma III.3.1 by virtue of the inequality $\text{dist}(\cdot, K)|\nabla\eta| \leq c(|\nabla(\text{dist}(\cdot, K)\eta)| + |\eta|)$.

Lemma III.5.1 can be invoked using $\zeta$ for Lemma III.5.2's version of $\eta$ to prove (III.5.38) if its two requirements are met. With regards to the requirement from the top bullet of Lemma III.5.1: The element $\mathcal{D}\zeta$ is in $\mathbb{L}$ because $|\mathcal{D}\zeta| \leq c|\eta|$ and $\eta$ is in $\mathbb{L}$. To see about the requirement from the second bullet, fix a small positive number $\varepsilon$ to use for Lemma III.5.2 and a likewise a number R. Take $t_\varepsilon$ from CONSTRAINT 1. Supposing that $t \in (0, \frac{1}{4}t_\varepsilon]$, then the $[t, 2t] \times Y$ integral of $\frac{1}{t^2}|\zeta|^2$ is no greater than the sum of its integrals over two regions: The first region is the part of $[t, 2t] \times Y$ where the distance to K is at most $100Rt$; and the second is the part of $[t, 2t] \times Y$ where the distance to K is greater than $32Rt$. The integral of $\frac{1}{t^2}|\zeta|^2$ over the first region is at most $cR^2$ times that of $|\eta|^2$ since $\frac{1}{t}|\zeta|$ is $\frac{\text{dist}(\cdot, K)}{t}|\eta|$ which is at most $100R|\eta|$. The next two paragraphs consider the integral of $\frac{1}{t^2}|\zeta|^2$ over the $\text{dist}(\cdot, K) \geq 32Rt$ part of $[t, 2t] \times Y$.

To bound the $\text{dist}(\cdot, K) \geq 32Rt$ part contribution, first write

$$\int_{[t,2t]\times(Y-\mathcal{N}_{32Rt})} |\zeta|^2 = \sum_n \int_{[t,2t]\times(\mathcal{N}_{32(n+1)Rt}-\mathcal{N}_{32nRt})} |\zeta|^2 + \int_{[t,2t]\times(Y-\mathcal{N}_{\frac{1}{100R}t_\varepsilon})} |\zeta|^2$$

(III.5.39)

where the sum is over the integers starting from 1 to the least integer (call it $n_t$) such that $32n_t R t$ is greater than $\frac{1}{100R} t_\varepsilon$. Supposing that n is a positive integer from $\{1, \ldots, n_t\}$, then the corresponding term in the sum on the right hand side of (III.5.39) is at most

$$c n^2 R^2 t^2 \int_{[t,2t]\times(\mathcal{N}_{32(n+1)Rt}-\mathcal{N}_{32nRt})} |\eta|^2 .$$

(III.5.40)

With this understood, invoke Lemma III.5.2 with $s = 8nRt$ to bound the preceding by

$$c t^2 \int_{(0,t_\varepsilon]\times(\mathcal{N}_{256nRt}-\mathcal{N}_{8nRt})} |\eta|^2 + c R^2 t^2 \int_{[8nRt, 16nRt]\times(Y-\mathcal{N}_{8Rt})} |\eta|^2 .$$

(III.5.41)

Sum the various integer n versions of this last bound to see that

$$c(1+R^2)t^2 \int_{(0,t_\varepsilon]\times Y} |\eta|^2$$

(III.5.42)

is an upper bound for the sum that appears on the right hand side of (III.5.39). Lemma III.5.2 with a suitable choice of s supplies a similar bound for the right most integral on the right most side of (III.5.39).

The bound in (III.5.42) for the sum and right-most integral on the right hand side of (III.5.39) lead directly to a t-independent upper bound for the $[t, 2t]\times Y$ integral of the function $\frac{1}{t^2}|\zeta|^2$ which is what is required by the second bullet of Lemma III.51.

*Part 2*: Return now to the milieu of Lemma III.2.2 and Section III.2b for the analysis near K. In what follows, the parameters r and $\delta$ are chosen to obey the constraints of Lemmas III.2.1 and III.2.2. Part 2 of Section III.2b uses CONSTRAINT 2's diffeomorphism to identify $(0, \infty)\times \mathcal{N}_K$ with $(0, \infty)\times D_0 \times S^1$ with $D_0$ denoting the radius $r_0$ disk about the origin in $\mathbb{R}^2$ as $D_0 \times S^1$. The isomorphism modulo $\{\pm 1\}$ from CONSTRAINT 2 is used to write the operator $\mathcal{D}$ near K as in (III.2.14) with $|\mathfrak{r}_1| \le c|z|$ and $|\mathfrak{r}_0| \le c$. This isomorphism writes $(A, \mathfrak{a})$ as in (III.2.15) on a small $\delta$ version of $(0, \delta]\times D_0 \times S^1$; it also identifies $\eta$ with an element $\eta´$ which obeys the equation

$$\mathcal{D}^{0(m)}\eta´ + \mathfrak{k}_1(\nabla)\eta´ + \mathfrak{k}_0\eta´ = 0$$

(III.5.43)

with $\mathfrak{k}_1$ and $\mathfrak{k}_0$ denoting endomorphisms that obey $|\mathfrak{k}_1| \le c|z|$ and $|\mathfrak{k}_0| \le c(\varepsilon + r + \frac{1}{R^2})\frac{1}{t}$.

The function $\chi_\delta$ is defined in (III.2.1) and define the $\mathbb{R}^2 \times S^1$ version of $\varpi_{\Diamond r}$ by analogy with (7.9) as $\varpi_{\Diamond r} = \chi(\frac{|z|}{r} - 1)$. Let $\eta_\Diamond = \chi_{4\delta} \varpi_{\Diamond r/4} \eta'$. This has compact support where (8.43) is valid, and as a consequence of (8.43) and what is said by (8.39), it obeys

$$\mathcal{D}^{0(m)} \eta_\Diamond = \mathfrak{w}$$

(III.5.44)

with $\mathfrak{w}$ being in the Euclidean metric's version of $\mathbb{L}$ for $(0,\infty) \times \mathbb{R}^2 \times S^1$.

*Part 3*: The operator $\mathcal{D}^{0(m)}$ maps the Euclidean metric's version of $\mathbb{H}$ for $(0,\infty) \times \mathbb{R}^2 \times S^1$ isomorphically to the corresponding version of $\mathbb{L}$. As a consequence, there is a unique element in this version of $\mathbb{H}$ (to be denoted by $u$) obeying $\mathcal{D}^{0(m)} u = \mathfrak{w}$. As a consequence, $\mathcal{D}^{0(m)}(\eta_\Diamond - u) = 0$, so $\eta_\Diamond - u$ is in the $\mathbb{L}$-kernel of $\mathcal{D}^{0(m)}$, which means that $\eta_\Diamond = u$ because the $\mathbb{L}$-kernel of $\mathcal{D}^{0(m)}$ is isomorphic to the $\mathbb{L}$-kernel of its formal adjoint which is trivial. (See (III.5.2).) Thus, $\eta_\Diamond$ is in the Euclidean metric's version of $\mathbb{H}$.

Given that $\eta_\Diamond$ has compact support where $t \leq \frac{1}{T}$ and where $|z| < r$, it follows that the product of $\chi_T \varpi_r \eta$ is in the $(0,\infty) \times Y$ version of $\mathbb{H}$. This implies in turn (via Hardy's inequality) that the $(0, 1] \times Y$ integral of the square of the norm of $\frac{1}{t} \chi_T \varpi_r \eta$ is finite.

## 6. The index of $\mathcal{D}$

This section completes the proof of Theorem III.1 by explaining why the index of $\mathcal{D}$ (as a map from $\mathbb{H}$ to $\mathbb{L}$) is the complex dimension of the version of the vector space $H^0$ that is defined by the $t \to \infty$ limit flat $SL(2;\mathbb{C})$ connection on Y.

To say more about what is to come: As noted in Section III.5, the $\mathbb{L}$ kernel of $\mathcal{D}^\dagger$ is isomorphic to the $\mathbb{L}$-kernel of $\mathcal{D}$ which is a subspace of the $\mathbb{H}$-kernel of $\mathcal{D}$. As a consequence, the index of $\mathcal{D}$ is the dimension of the complementary subspace. As explained Section III.6a, this complementary subspace is $\{0\}$ when the limit flat connection is irreducible. Sections III.6b-III.6d treat the case when the limit flat connection is not irreducible: These sections construct an isomorphism (over $\mathbb{R}$) between the complementary subspace and a real subspace of the complex vector space $H^0$ that generates the latter as a vector space over $\mathbb{C}$.

### a) The issue with reducibility

Assume for now that the pair $(A, \mathfrak{a})$ obeys the assumptions of Theorem III.1. By virtue of what is said in CONSTRAINT 3, the $(A, \mathfrak{a})$ version of the operator $\mathcal{D}$ can be written as in (III.3.4) where t is large. By way of a reminder: The bundle P is identified via a bundle isomorphism with the pull-back of a principal SU(2) bundle on Y, pull-back defined by the projection map from $(0,\infty)\times Y$ to Y. (Both the pull-back bundle and the original on Y will be denoted henceforth by P.) This isomorphism is such that $(A, \mathfrak{a})$ can be written as $(A^\infty, \mathfrak{a}^\infty) + \mathfrak{v}$ with $\mathbb{A} = A^\infty + i\mathfrak{a}^\infty$ defining a flat $SL(2;\mathbb{C})$ connection on Y and with the norm of $\mathfrak{v}$ being bounded by $\frac{\varepsilon}{t}$ when t is greater than some time $t_\varepsilon$. What is written in (III.3.4) as $\mathcal{D}^\infty$ is the version of $\mathcal{D}$ that is defined by $(A^\infty, \mathfrak{a}^\infty)$. Meanwhile, the norm of the endomorphism $\mathfrak{r}$ that appears in (III.3.4) is pointwise at most $c\frac{\varepsilon}{t}$ when t is greater than $t_\varepsilon$.

Supposing that $\psi$ is in the $\mathbb{H}$ kernel of $\mathcal{D}$, write $\psi$ at fixed $t \geq t_\varepsilon$ as a section of the bundle depicted in (III.3.14) to be denoted by $(\eta, \nu)$. Then write the $\gamma_i \nabla_{A^\infty_i} + \rho_i[\mathfrak{a}^\infty_i, \cdot\,]$ part of $\mathcal{D}^\infty \psi$ as in (III.3.15). Since the limit flat $SL(2;\mathbb{C})$ connection is regular, the inequality in (III.3.11) holds and so

$$\int_{\{t\}\times Y} |\gamma_i \nabla_{A^\infty_i} \psi + \rho_i[\mathfrak{a}^\infty_i, \psi]|^2 \geq c^{-1} \int_{\{t\}\times Y} |\eta|^2 + c^{-1} \int_{\{t\}\times Y} |d_{\mathbb{A}^*}\nu|^2 ,$$

(III.6.1)

where $t > t_\varepsilon$. Note in particular that if $\mathbb{A}$ is not just regular but also irreducible, then (III.3.12) holds and then (III.6.1) implies that

$$\int_{\{t\}\times Y} |\gamma_i \nabla_{A^\infty_i} \psi + \rho_i[\mathfrak{a}^\infty_i, \psi]|^2 \geq c^{-1} \int_{\{t\}\times Y} |\psi|^2$$

(III.6.2)

where $t > t_\varepsilon$.

This last bound with (III.3.4) and Lemma III.2.1 plus Hardy's inequality (I.2.3) lead directly to this: If $\mathbb{A}$ is regular and irreducible, and if $\psi$ is in the $\mathbb{H}$-kernel of $\mathcal{D}$, then

$$\|\psi\|_{\mathbb{H}} \geq c^{-1} \|\psi\|_{\mathbb{L}} .$$

(III.6.3)

Thus the $\mathbb{H}$ and $\mathbb{L}$ kernels of $\mathcal{D}$ are identical when $\mathbb{A}$ is regular and irreducible, so $\mathcal{D}$ has index 0 in this case

### b) The $t \to \infty$ asymptotics of $\nu$

To analyze the case when $\mathbb{A}$ is reducible: Let $\mathbb{L}_Y$ denote the Hilbert space completion of the space of smooth sections over Y of $\oplus_2 (ad(P) \oplus (ad(P) \otimes T^*Y))$ using the norm whose square is defined by the rule

$$\psi \to \int_Y |\psi|^2 \ .$$

(III.6.4)

(The bundle $\oplus_2 (\mathrm{ad}(P) \oplus (\mathrm{ad}(P) \otimes T^*Y))$ over Y is denoted henceforth as $\mathbb{W}_Y$.)

The operator $\gamma_i \nabla_{A^\infty_i} + \rho_i[\mathfrak{a}^\infty_i, \cdot\,]$ can be viewed as an unbounded, essentially self-adjoint operator acting on $\mathbb{L}_Y$. As such, it has pure point spectrum with no accumulation points and finite multiplicities. By virtue of (7.35)-(7.37), its kernel consists of the sections of the ad(P) summands in $\mathbb{W}_Y$ that are $A^\infty$-covariantly constant and commute with $\mathfrak{a}^\infty$. Let $\Pi^0$ denote the $\mathbb{L}_Y$-orthogonal projection onto this kernel. (The dimension of this kernel is $2\dim_\mathbb{C} H^0$ which is at most six.) Also, let $\Pi^+$ denote the $\mathbb{L}_Y$-orthogonal projections onto the span of the eigenvectors of $\gamma_i \nabla_{A^\infty_i} + \rho_i[\mathfrak{a}^\infty_i, \cdot\,]$ with positive eigenvalue; and let $\Pi^-$ denote the corresponding projection to the span the eigenvectors of this same operator with negative eigenvalue.

Now let $f_\pm$ and $f_0$ denote the respective functions on $[t_\varepsilon, \infty)$ whose values are the $\mathbb{L}_Y$-norms on $\{t\} \times Y$ of $\Pi^+\psi$, $\Pi^-\psi$ and $\Pi^0\psi$. These functions and $\Pi^0\psi$ itself obey (by virtue of (III.3.4)) the differential inequalities written below when $\varepsilon < c^{-1}$ (the inequalities use $\lambda$ to denote half of the smallest of the absolute values of the non-zero eigenvalues of the operator $\gamma_i \nabla_{A^\infty_i} + \rho_i[\mathfrak{a}^\infty_i, \cdot\,]$).

- $\frac{d}{dt} f_+ + \lambda f_+ \le c \frac{\varepsilon}{t} f_- + c \frac{\varepsilon}{t} f_0$ .
- $\frac{d}{dt} f_- - \lambda f_- \ge -c \frac{\varepsilon}{t} f_+ - c \frac{\varepsilon}{t} f_0$ .
- $|\nabla_{A^\infty_t} \Pi^0 \psi| \le c \frac{\varepsilon}{t} f_+ + c \frac{\varepsilon}{t} f_-$ .

(III.6.5)

With regards to the third bullet: There is no $c \frac{\varepsilon}{t} f_0$ term on the right hand side because $\Pi^0 \psi$ has only the contributions from the parts of $\psi$ proportional to dt, thus the $\mathfrak{b}_t$dt and $\mathfrak{c}_t$dt parts (and $\mathfrak{c}_t \equiv 0$); and one can see from the form of $\mathcal{D}$ in (I.1.2) that these parts are sent by the $\gamma_i \nabla_i + \rho_i[\mathfrak{a}_i, \cdot\,]$ part of $\mathcal{D}$ to ad(P)-valued sections of T*Y. With regards to the first two bullets in (III.6.5): There is no $c \frac{\varepsilon}{t} f_+$ in the top inequality when $\varepsilon < c^{-1}$ and no $c \frac{\varepsilon}{t} f_-$ in the second one, they are accounted for by the taking $\lambda$ to be *half* of the smallest absolute value of a non-zero eigenvalue instead being equal to the smallest absolute value.

The rest of this section first states and then proves two lemmas concerning the large t behavior of $f_\pm$ and $\Pi^0 \psi$.

**Lemma III.6.1**: *The functions $f_+$ and $f_-$ are bounded on $[t_\varepsilon, \infty)$; and both $f_+^2$ and $f_-^2$ have finite integral. Meanwhile, the $\Pi^0 \psi$ is bounded and it has a unique $t \to \infty$ limit.*

Note that this lemma asserts in part that $|\Pi^+\psi|^2$ and $|\Pi^-\psi|^2$ have finite integral on $[t_\varepsilon, \infty) \times Y$. Note also that $\Pi^0\psi$ takes its values in a $\dim_\mathbb{C} H^0$-dimensional subspace of the image of $\Pi^0$ (half the dimension of this image) because the $\mathfrak{c}_t$ component of $\psi$ is identically zero.

**Lemma III.6.2**: *If* $\lim_{t\to\infty} \Pi^0\psi = 0$, *then* $|\Pi^0\psi|^2$ *has finite integral on* $[t_\varepsilon, \infty) \times Y$.

The introduction to this Section III.6 identified the index of $\mathcal{D}$ with the dimension of the quotient of the $\mathbb{H}$-kernel of $\mathcal{D}$ by its subspace with finite $\mathbb{L}$-norm. This lemma plus what Lemma III.6.1 says about $|\Pi^+\psi|^2$ and $|\Pi^-\psi|^2$ having finite integral on $[t_\varepsilon, \infty) \times Y$ implies that the dimension of this quotient (and hence the index of $\mathcal{D}$) is at most $\dim_\mathbb{C} H^0$.

***Proof of Lemma III.6.1***: To see about $f_+$: Multiply the top inequality in (III.6.5) by $f_+$. Fix $t > t_\varepsilon$ and integrate from $t_\varepsilon$ to $t$ to see that

$$f_+^2(t) + \lambda \int_{t_\varepsilon}^t f_+^2 \leq f_+^2(t_\varepsilon) + c\varepsilon \|\psi\|_\mathbb{H}^2 .$$

(III.6.6)

The bound on the right hand side arises because the $[t_\varepsilon, \infty)$ integrals of $\frac{1}{t^2} f_\pm^2$ and $\frac{1}{t^2} f_0^2$ are finite and bounded by $c\|\psi\|_\mathbb{H}^2$ courtesy of Hardy's inequality. The inequality in (III.6.6) establishes what the lemma claims with regards to $f_+$.

To see about $f_-$: Fix a small positive number to be denoted by $\delta$ (make it less than $\frac{1}{4}\lambda$) and then multiply both sides of the second inequality in (III.6.5) by $e^{-\delta t} f_-$. Having fixed a time $t > t_\varepsilon$, integrate the result over the interval $[t, \infty)$. The factor $e^{-\delta t}$ makes the integrals on both sides a priori finite; and it allows for an integration by parts on the left hand side integral with the only non-zero boundary contribution from the time $t$ endpoint of $[t, \infty)$. (Remember in this regard that $\frac{1}{t^2} f_-^2$ has finite integral on $[t_\varepsilon, \infty)$ and that $e^{-\delta t}$ is much less than $\frac{1}{t^2}$ when $t$ is very large). Doing this integration by parts leads to the inequality

$$e^{-\delta t} f_-^2(t) + \tfrac{1}{2}\lambda \int_{t_\varepsilon}^t e^{-\delta t} f_-^2 \leq f_+^2(t_\varepsilon) + c\varepsilon \|\psi\|_\mathbb{H}^2 .$$

(III.6.7)

Since the right hand side is independent of $\delta$, the $\delta \to 0$ limit can be taken on the left hand side to establish what the lemma claims for $f_-$.

The lemma's claims about $\Pi^0\psi$ follow from the lower inequality in (III.6.5) because $f_+$, $f_-$ and $\frac{1}{t^2}$ have finite integral on $(t_\varepsilon, \infty)$. Keep in mind in this regard that the 1-parameter family of ad($P_\infty$)-sections $\{\Pi^0\psi\}_{t>t_\varepsilon}$ all lie in a $\dim_\mathbb{C} H^0$ dimensional vector subspace of $C^\infty(Y; \mathrm{ad}(P_\infty))$. (The dimension is $\dim_\mathbb{C} H^0$ instead of $2\dim_\mathbb{C} H^0$ because the $\mathfrak{c}_t$ part of $\psi$ is identically zero.)

*Proof of Lemma III.6..2*: Suppose for the moment that $f$ is any given function on the interval $[t_\varepsilon, \infty)$. Fix $T > t_\varepsilon$. An integration by parts says that

$$\int_{t_\varepsilon}^T f^2 \, dt = -2 \int_{t_\varepsilon}^T f t (\tfrac{d}{dt} f) \, dt + f^2(T) - f^2(t_\varepsilon) ;$$

(III.6.8)

and this leads to another version of Hardy's inequality:

$$\int_{t_\varepsilon}^T f^2 \, dt \leq c \int_{t_\varepsilon}^T t^2 \, |\tfrac{d}{dt} f|^2 \, dt + c \, f^2(T) .$$

(III.6.9)

In particular, if $\lim_{T\to\infty} f(T) = 0$, then

$$\int_{t_\varepsilon}^T f^2 \, dt \leq c \int_{t_\varepsilon}^T t^2 \, |\tfrac{d}{dt} f|^2 \, dt .$$

(III.6.10)

Meanwhile, it follows from Lemma III.6.1 that the function $t^2 |\nabla_{A^\infty_t} \Pi^0 \psi|^2$ has finite integral on $[t_\varepsilon, \infty) \times Y$. Therefore, by virtue of (III.6.10) with $f(\cdot)$ taken to be the square root of the $\{\cdot\} \times Y$ of $|\Pi^0\psi|^2$, the function $|\Pi^0\psi|^2$ has finite integral on $[t_\varepsilon, \infty) \times Y$ in the event that $\lim_{t\to\infty} \Pi^0 \psi = 0$.

**c) Constructing elements in with $\lim_{t\to\infty} \Pi^0\psi \neq 0$**

The assumption here is that $(A, \mathfrak{a})$ obeys the equations in (I.2.1) in addition to CONSTRAINTS 1-3. This section constructs a $\dim_\mathbb{C} H^0$ dimensional subspace of elements in the $\mathbb{H}$-kernel of $\mathcal{D}$ with each member (except 0) having infinite $\mathbb{L}$-norm. To set the stage for the construction, suppose for the moment that $\phi$ is any given section of ad(P) over $(0, \infty) \times Y$. With $\phi$ in hand, define $\mathrm{B} = -\nabla_{At}\phi \, dt - d_A \phi$ and $\mathrm{C} = [\phi, \mathfrak{a}]$. The pair (B, C) is such that all but the $\mathfrak{p}_t$ component of (I.1.2) are zero (the component given by the second bullet). This is so because: For any $\varepsilon > 0$, the bundle automorphism $\exp(\varepsilon\phi)$ pulls $(A, \mathfrak{a})$ back to give another solution to (I.2.1). Meanwhile, the pair (B,C) is the first order in $\varepsilon$ change of $(A, \mathfrak{a})$ via the action of this bundle automorphism; and when (I.2.1) is viewed as an instance of (I.1.1), then the first, third and fourth components of (I.1.2) are the corresponding first order change of the components of (I.1.1)–which, as noted, is zero. The $\mathfrak{p}_t$ component of (I.1.2) is zero if and only if $\phi$ obeys the differential equation

$$\nabla^\dagger \nabla \phi + [\mathfrak{a}_i, [\phi, \mathfrak{a}_i]] = 0 \ .$$

(III.6.11)

(The notation uses $\nabla^\dagger$ to denote the formal $\mathbb{L}$-adjoint of $\nabla$.) The task for this section is to construct a $\dim_\mathbb{C} H^0$-dimensional space of $\phi$ that obey (III.6.11) with the corresponding $\psi$ inside $\mathbb{H}$ but without finite $\mathbb{L}$-norm when $\phi$ is not identically zero. The upcoming construction has four parts.

   *Part 1*: Fix a section over Y of ad(P) that is annihilated by both $d_{A^\infty}$ and $[\mathfrak{a}^\infty, \cdot]$; it will denoted by $x$. Then set $\phi_x$ to be the section of ad(P) over $(0, \infty) \times Y$ given by

$$\phi_x = t\chi(2 - \tfrac{t}{t_\varepsilon})x.$$

(III.6.12)

This is a section of ad(P) on $(0, \infty) \times Y$ which is zero where $t < t_\varepsilon$ and which is equal to $tx$ where $t > 2t_\varepsilon$. Although $\phi_x$ doesn't obey (III.6.11), it does have the following features: The corresponding $B = -\nabla_t \phi_x + d_A \phi_x$ and $C = [\phi_x, \mathfrak{a}]$ are pointwise bounded and, where $t > 2t_\varepsilon$, they obey

- $|d_A \phi_x| + |[\phi_x, \mathfrak{a}]| \leq c\varepsilon$.
- $|\nabla_t \phi_x - x| \leq c\varepsilon$ .

(III.6.13)

Although these properties are sufficient for the purpose at hand, the argumentsa are far simpler if an assertion of the upcoming Theorem IV.1 is assumed which is that $(A, \mathfrak{a})$ can be written where $t$ is large (after the action of a suitable automorphism of P) as

$$(A, \mathfrak{a}) = (A^\infty + \mathcal{A}, \mathfrak{a}^\infty + \mathfrak{v}) \quad with \quad |\mathcal{A}| + |\mathfrak{v}| \leq e^{-\delta t} \ .$$

(III.6.14)

As a consequence of this last bound, there exists some fixed $t_*$ such that when $t \geq t_*$, then

- $|d_A \phi_x| + |[\phi_x, \mathfrak{a}]| \leq e^{-\delta t/2}$.
- $|\nabla_t \phi_x - x| \leq e^{-\delta t/2}$.

(III.6.15)

The arguments that follow will cheat and assume that (III.6.14) holds and hence (III.6.15). (The upcoming proof of Theorem IV.1 does not use what is said in this section.)

*Part 2*:  Let $\mathbb{H}_0$ denote the Hilbert space completion of the space of smooth, compactly supported sections of ad(P) over $(0,\infty) \times Y$ using the norm whose square is

$$\int_{(0,\infty)\times Y} (|\nabla(\cdot)|^2 + |[\mathfrak{a},\cdot]|^2)$$

(III.6.16)

This norm is denoted by $\|\cdot\|_{\mathbb{H}}$ in what follows. Of particular import:  The section $\phi_x$ is <u>not</u> in the space $\mathbb{H}_0$. (This can be deduced from the fact that $\frac{1}{t^2}|\phi_x|^2$ does not have finite integral on $(0,\infty) \times Y$.)

The plan for Part 3 is to find an element in $\mathbb{H}_0$ to be denoted by $u$ so that $\phi = \phi_x + u$ obeys (III.6.11). By virtue of $u$ being in $\mathbb{H}_0$ and $\phi_x$ not, the element $\phi$ is not identically zero.  As explained in Part 4, the corresponding element ($\mathrm{B} = -\nabla_t \phi + d_A \phi$, $\mathrm{C} = [\phi, \mathfrak{a}]$) is in the $\mathbb{H}$-kernel of $\mathcal{D}$ but not the $\mathbb{L}$-kernel.

*Part 3*:  The desired section $u$ is the minimizer of the function on $\mathbb{H}_0$ that is denoted by $\mathfrak{Z}$ and whose value on any given element $w$ is given by the rule

$$\mathfrak{Z}(w) = \tfrac{1}{2} \int_{(0,\infty)\times Y} (|\nabla w|^2 + |[\mathfrak{a},w]|^2) + \mathfrak{X}(w)$$

(III.6.17)

where $\mathfrak{X}$ is shorthand for the following function of $w$:

$$\int_{(0,\infty)\times Y} (-\langle w, (\nabla^2_{A^\infty_t}\phi_x + [\mathfrak{v}_t, \nabla_{A^\infty_t}\phi_x])\rangle + \langle \nabla_{A_t} w, [\mathfrak{v}_t, \phi_x]\rangle + \langle d_A w, d_A \phi_x \rangle + \langle [\mathfrak{a}, w], [\mathfrak{a}, \phi_x]\rangle).$$

(III.6.18)

This function $\mathfrak{X}$ is observedly linear in $w$ and, by virtue of (III.6.14) and (III.6.15), it is bounded, so it is a bounded, linear functional on $\mathbb{H}_0$, hence continuous.  As a consequence, the function $\mathfrak{Z}$ is continuous on $\mathbb{H}_0$ and bounded from below on $\mathbb{H}_0$:

$$\mathfrak{Z}(w) \geq \tfrac{1}{4} \|w\|_{\mathbb{H}}^2 - c.$$

(III.6.19)

Because of these last observations, standard arguments in the calculus of variations imply that $\mathfrak{Z}$ has a unique minimizer in $\mathbb{H}_0$ and that this minimizer (the desired $u$) is such that $\phi \equiv \phi_x + u$ obeys (III.6.11).  (For existence:  Take a minimizing sequence $\{w_n\}_{n \in \mathbb{N}}$ for $\mathfrak{Z}$. As a consequence of (III.6.19), this sequence has bounded $\mathbb{H}$-norm.  Thus, it has a weakly convergent subsequence to a limit in $\mathbb{H}$.  The limit is non-zero because the infimum of $\mathfrak{Z}$ is negative due to it being the sum of a positive multiple of the *square* of the norm and a linear term. The minimizer a priori obeys the equations in in (III.6.11) by virtue of it being a critical point of $\mathfrak{Z}$. The minimizer is unique because the functional $\mathfrak{Z}$ is convex.)

*Part 4*: With $\phi$ as just defined, this last part of the subsection explains why the kernel element $\psi$ with components ($B = -\nabla_t\phi\, dt + d_A\phi$, $C = [\phi, \mathfrak{a}]$) is in $\mathbb{H}$ but not in $\mathbb{L}$. In this regard: The element $\psi$ is not in $\mathbb{L}$ because $\nabla_t\phi = \nabla_t\phi_x + \nabla_t u$ and $\nabla_t u$ is in $\mathbb{L}$ but $\nabla_t\phi_x$ is not (see (III.6.13) or (III.6.15)).

As explained directly, Lemma III.5.1 with $\eta$ in the latter being $\psi$ can be brought to bear to see that $\psi$ is in $\mathbb{H}$. Indeed, the criteria in the first bullet is met because $\mathcal{D}\psi = 0$. To see about the criteria in the second bullet, note first that $|\psi|^2$ has finite integral on $(0,1] \times Y$. Because of this, the arguments from Sections III.5c and III.5d can be repeated with only minor modifications to see that the function on $(0, 1]$ given by the rule

$$t \to \tfrac{1}{t^2} \int_t^{2t} \int_{\{\cdot\}\times Y} |\psi|^2$$

(III.6.20)

is bounded. Meanwhile, the function on $[1, \infty)$ given by the preceding rule is also bounded because $u$ is in $\mathbb{H}$ and because of (III.6.13).

## IV. THE BEHAVIOR OF KAPUSTIN-WITTEN SOLUTIONS ON $[1,\infty) \times Y$ AS $t \to \infty$

This last series of 'lectures' states and then proves a theorem about the $t \to \infty$ limit of a solution to (I.2.1) on $[1, \infty) \times Y$. This is Theorem IV.1.

To set the notation for Theorem IV.1, let $P^\infty$ denote a principal SU(2) bundle over $Y$ and let $(A^\infty, \mathfrak{a}^\infty)$ denote a pair of connection on $P^\infty$ and ad($P^\infty$)-valued 1-form on $Y$. The projection map from $[1,\infty) \times Y$ to $Y$ will be used to view $P^\infty$ as a bundle over $[1,\infty) \times Y$ and to likewise view $(A^\infty, \mathfrak{a}^\infty)$ as a pair defined over $[1,\infty) \times Y$. Theorem IV.1 also refers to the notion of a *regular* flat Sl(2; $\mathbb{C}$) connection on $Y$. This is a flat connection whose $H^1$ cohomology is zero. (The vector space $H^1$ is defined in (III.1.1).)

**Theorem IV.1**: *Let* $P \to [1, \infty) \times Y$ *denote a principal S*U(2) *bundle and let* $(A, \mathfrak{a})$ *denote a pair of connection on* P *and* ad(P)-*valued section of* T*Y *that over* $[1,\infty) \times Y$ *that obeys the equations in (I.2.1). If the lim-inf of the function on* $[1, \infty)$ *given by the rule*

$$t \to \int_{\{t\}\times Y} |\mathfrak{a}|^2$$

*is finite, then there exists the following data:*

- *A principal SU(2) bundle $P^\infty \to Y$ and a pair $(A^\infty, \mathfrak{a}^\infty)$ of connection on $P^\infty$ and $\mathrm{ad}(P^\infty)$-valued 1-form on Y with $A^\infty + i\mathfrak{a}^\infty$ being a flat $Sl(2;\mathbb{C})$ connection and with $d_{A^\infty} * \mathfrak{a}^\infty = 0$.*
- *An isomorphism between $P^\infty$ and P over $[1, \infty) \times Y$ that identifies $(A, \mathfrak{a})$ as a pair of connection on this pull-back and $\mathrm{ad}(P^\infty)$ section of $T^*Y$ over $[1, \infty) \times Y$.*

*Granted the second bullet's identification, these are such that*

$$\lim_{t \to \infty} \sup_{\{t\} \times Y} (|A - A^\infty| + |\mathfrak{a} - \mathfrak{a}^\infty|) \to 0.$$

*Moreover, if $A^\infty + i\mathfrak{a}^\infty$ defines a regular flat $SL(2;\mathbb{C})$ connection, then the second bullet's identification can be chosen so that*

$$\lim_{t \to \infty} \sup_{\{t\} \times Y} e^{\delta t}(|A - A^\infty| + |\mathfrak{a} - \mathfrak{a}^\infty|) \to 0$$

*with $\delta$ being a positive number that depends only on the limit pair $(A^\infty, \mathfrak{a}^\infty)$*

The proof of the first assertion about pointwise convergence is in Sections IV.1-IV.4. The proof of the second assertion about exponential convergence is in Section IV.5.

By way of a remark concerning the hypothesis: If the values of the $\{t\} \times Y$ integrals of $|\mathfrak{a}|^2$ diverge as $t \to \infty$, then there is still something to be said when the function on $[1, \infty)$ given by the various $\{t\} \times Y$ integrals of $\langle \mathfrak{a} \wedge *B_A - \frac{1}{3} \mathfrak{a} \wedge \mathfrak{a} \wedge \mathfrak{a} \rangle$ is bounded. Assuming such is the case, then the main theorem in [T1] can be invoked for the restriction of $(A, \mathfrak{a})$ to the constant t slices of $[1, \infty) \times Y$ to obtain a $\mathbb{Z}/2$ harmonic 1-form data set limit for the $t \to \infty$ limit of $(A, \mathfrak{a})$. More is said about this in Section IV.6.

## 1. A priori bounds and sequential convergence

The purpose of this section is to first establish a priori bounds for solutions to (I.2.1) on $[1, \infty) \times Y$ when the $\{t\} \times Y$ integrals of $|\mathfrak{a}|^2$ are bounded. These are used to prove that there are sequential limits of $(A, \mathfrak{a})$ as $t \to \infty$ after the application of suitable automorphisms of P.

### a) A priori bounds for solutions to (I.2.1)

The assumption in what follows is that $(A, \mathfrak{a})$ is a pair of smooth connection on P over $[1, \infty) \times Y$ and smooth $\mathrm{ad}(P)$-valued section of $T^*Y$ over this same domain. Also: Take this pair to obey the equations in (I.2.1). The following proposition gives a priori bounds on the norms of $E_A$, $B_A$ and $\nabla \mathfrak{a}$ and their derivatives to any order assuming a relatively weak, integral bound for $|\mathfrak{a}|^2$.

**Proposition IV.1.1**: *Define a function on $[1, \infty)$ by the rule*

$$t \to \int_{\{t\} \times Y} |\mathfrak{a}|^2 .$$

*If this function is bounded, then*

- $|\mathfrak{a}|$ *is bounded*.
- *The norms of* $E_A$. $B_A$ *and* $\nabla\mathfrak{a}$ *are bounded as are those of their A-covariant derivatives to any given order*.
- $|B_A - *(\mathfrak{a}\wedge\mathfrak{a})|^2$ *and* $|d_A\mathfrak{a}|^2$ *have finite integral on* $[2,\infty)\times Y$; *thus so do* $|\nabla_t\mathfrak{a}|^2$ *and* $|E_A|^2$.

*Proof of Proposition IV.1.1*: The proof of the assertions in the first two bullets won't be given in detail because it is a fairly straightforward example of what are called 'elliptic bootstrap' arguments. See [T2] for more about how the elliptic bootstrap arguments for the solutions to (I.2.1) work. What follows is a brief outline of these arguments: By virtue of (I.2.1), the section $\mathfrak{a}$ obeys the second order equation

$$\nabla^\dagger\nabla\mathfrak{a} + [\mathfrak{a}_i,[\mathfrak{a},\mathfrak{a}_i]] + \text{Ric}(\mathfrak{a}) = 0 ,$$

(IV.1.1)

where $\text{Ric}(\cdot)$ denotes the endomorphism of $T^*Y$ that is defined by the Ricci tensor and the metric. This last identity implies turn the identity

$$\nabla^\dagger\nabla|\mathfrak{a}|^2 + |\nabla\mathfrak{a}|^2 + |[\mathfrak{a},\mathfrak{a}]|^2 + \langle\mathfrak{a}, \text{Ric}(\mathfrak{a})\rangle = 0 .$$

(IV.1.2)

Integration of both sides of (IV.1.2) against suitably chosen Green's functions for $\nabla^\dagger\nabla$ leads to an a priori bound for $|\mathfrak{a}|$. Integrating both sides of (IV.1.2) against bump functions supported in small radius balls with centers where $t \geq \frac{3}{2}$ leads to a prior bounds for the integrals of $|\nabla\mathfrak{a}|^2$ and $|[\mathfrak{a},\mathfrak{a}]|^2$ over these balls: If a ball has radius r, then the bound has the form $cr^2$. These bounds also apply to the integrals of $|E_A|^2$ and $|B_A|^2$. Granted all of those bounds, take the $\nabla$-covariant derivative of (IV.1.1) to get an elliptic equation for $\nabla\mathfrak{a}$ which leads (using Green's functions and bump functions) to pointwise bounds for $\nabla\mathfrak{a}$ and $E_A$ and $B_A$, and integral bounds for the square of the norm of their covariant derivatives. Then, one takes a second covariant derivatives of (IV.1.1) and repeats the process to get pointwise bounds for $\nabla\nabla\mathfrak{a}$ and $\nabla E_A$ and $\nabla B_A$; and so on… .

The third bullet follows from the first two by integrating the identity

$$\tfrac{d}{dt}\int_{\{\cdot\}\times Y}\langle\mathfrak{a}\wedge *B_A - \tfrac{1}{3}\mathfrak{a}\wedge\mathfrak{a}\wedge\mathfrak{a}\rangle = \int_{\{\cdot\}\times Y}(|B_A - *(\mathfrak{a}\wedge\mathfrak{a})|^2 + |d_A\mathfrak{a}|^2)$$

(IV.1.3)

on intervals of the form [1, n] for n = 2, 3, … . (The top two bullets of the proposition imply that the integral of the right hand side of (10.3) on [1, n] has an n-independent upper bound.)

### b) Sequential convergence to flat connections

Proposition IV.1.1 has an immediate almost corollary which is stated below as Proposition IV.1.2. To set the stage for the lemma, suppose for the moment that $P_\infty$ denotes a principal SU(2) bundle over Y. An SL(2;$\mathbb{C}$) connection on $P_\infty \times_{SU(2)} SL(2;\mathbb{C})$ is said to be *stable* if, when writing it as $A^\infty + i\mathfrak{a}^\infty$, the pair $(A^\infty, \mathfrak{a}^\infty)$ is such that $d_{A^\infty} *\mathfrak{a}^\infty = 0$. (As has been the case previoulsy, $A^\infty$ signifies a connection on P and $\mathfrak{a}^\infty$ signifies an ad(P) valued section of T*Y.)

***Proposition IV.1.2***: *Suppose that the function depicted in Proposition IV.1.1 is bounded. Let $\{s_n\}_{n \in \mathbb{N}}$ denote any unbounded, increasing sequence in $[1, \infty)$. There exists*
- *A subsequence of $\mathbb{N}$ to be denoted by $\Lambda$,*
- *A principal SU(2) bundle over Y to be denoted by $P^\infty$,*
- *A stable flat SL(2;$\mathbb{C}$) connection on $P^\infty \times_{SU(2)} SL(2;\mathbb{C})$ to be denoted by $\mathbb{A}$.*

*These have the following significance: For each $n \in \Lambda$, pull-back P and the Sl(2;$\mathbb{C}$) connection $A + i\mathfrak{a}$ to the domain $(1 - s_n, \infty) \times Y$ using the map from $[1, \infty) \times Y$ given by the rule $(t, y) \to (t + s_n, y)$. For each such n, there is an isomorphism between this pull-back of P over $(1 - s_n, \infty) \times Y$ and the pull-back of $P^\infty$ via the projection to Y that identifies the pull-back of $A + i\mathfrak{a}$ with an SL(2;$\mathbb{C}$) connection such that the resulting $\Lambda$-labeled sequence of SL(2;$\mathbb{C}$) connections converges to $\mathbb{A}$ in the $C^\infty$-topology on bounded subsets of $\mathbb{R} \times Y$.*

***Proof of Proposition IV.1.2***: For each $n \in \mathbb{Z}$, let $P_n$ denote the pull-back of P to the domain $[1 - s_n, s_n] \in Y$; and let $A^n + i\mathfrak{a}^n$ denote the corresponding pull-back of $A + i\mathfrak{a}$. As done in [T2] for example, the top two bullets of Proposition IV.1.1 with Karen Uhlenbeck's bounded curvature theorem in [U] can be used to obtain the items in the first two bullets, and also:

- *A principal SU(2) bundle over $\mathbb{R} \times Y$ and a stable, SL(2;$\mathbb{C}$) connection on the associated, principal SL(2;$\mathbb{C}$) bundle.*
- *For each $n \in \Lambda$, a corresponding isomorphism over $(1 - s_n, \infty) \times Y$ from the principal SU(2) bundle in the preceding bullet to $P_n$ (the isomorphism is denoted by $g_n$).*

(IV.1.4)

Their significance is as follows: The $\Lambda$-labeled sequence whose n'th term is the $g_n$ pull-back of $A^n + i\mathfrak{a}^n$ converges in the $C^\infty$-topology on compact subsets of $\mathbb{R} \times Y$ to the SL(2;$\mathbb{C}$) connection. (Given Uhlenbeck's theorem, the convergence is due to the Arzoli-Ascoli theorem.) The SL(2;$\mathbb{C}$) connection in the first bullet of (IV.1.4) has to be stable because $d_{A^n} *\mathfrak{a}^n \equiv 0$ (because $d_A *\mathfrak{a} \equiv 0$).

To see that it is flat and that it is isomorphic to one pulled back from Y, fix some large (>> 1) positive number T. If n is such that $s_n > T + 1$, then the $[-T, T] \times Y$ integral of $|B_{A^n} - *(\mathfrak{a}^n \wedge \mathfrak{a}^n)|^2 + |d_{A^n}\mathfrak{a}^n|^2$ is equal to the $[s_n - T, T] \times Y$ integral of $|B_A - \mathfrak{a} \wedge \mathfrak{a}|^2 + |d_A\mathfrak{a}|^2$. The sequence of these $n \in \Lambda$ labeled integrals has to limit to zero as $n \to \infty$ for any fixed T because of what is said by the third bullet of Proposition IV.1.1. In addition, because any given $n \in \Lambda$ labeled integral is equal to the integral of $|\nabla_{A^n t}\mathfrak{a}^n|^2 + |E_{A^n}|^2$ over the domain $[-T, T] \times Y$, the corresponding $n \in \Lambda$ labeled sequence of the latter integrals also limits to zero as $n \to \infty$ (for fixed T). The fact that these two limits are zero for any fixed T implies that the limit $SL(2;\mathbb{C})$ connection from the first bullet in (IV.1.4) is flat and that it and its principal SU(2) bundle are isomorphic to the pull-backs of a flat $SL(2;\mathbb{C})$ connection on the associated $SL(2;\mathbb{C})$ bundle to a principal SU(2) bundle over Y. (Write the limit $SL(2;\mathbb{C})$ connection from the first bullet in (IV.1.4) as $A^\infty + i\mathfrak{a}^\infty$. Then use the $A^\infty$ parallel transport along the $\mathbb{R}$ factor in $\mathbb{R} \times Y$ to construct the desired isomorphism with the principal SU(2) bundle on Y being the restriction of the principal SU(2) bundle from the first bullet of (IV.1.4) to $\{0\} \times Y$.)

## 2. Leon Simon's Lojasiewicz inequality

The proof of the unique $t \to \infty$ limit assertion in Theorem IV.1 is a version of a theorem of Leon Simon (this is Theorem 3 in [S]) and worked out in the related context of the anti-self dual Yang-Mills equation by Morgan, Mrowka and Ruberman [MMR]. The account here closely follows [MMR].

### a) The Lojasiewicz inequality

To set the stage for the story, suppose for the moment that P denotes a principal SU(2) bundle over $[1, \infty) \times Y$ and that $(A, \mathfrak{a})$ is a pair of connection on P and section of $ad(P) \otimes T^*Y$ that obeys the equations in (I.2.1) and is such that the function depicted of t depicted in Proposition IV.1.1 is bounded. Define a second function on $[1, \infty)$, to be denoted by $\mathfrak{cs}$, by the rule

$$t \to \mathfrak{cs}(t) = \int_{\{\cdot\} \times Y} \langle \mathfrak{a} \wedge B_A - \tfrac{1}{3}\mathfrak{a} \wedge \mathfrak{a} \wedge \mathfrak{a}\rangle \ .$$

(IV.2.1)

The t-derivative of $\mathfrak{cs}$ is depicted in (IV.1.3). It can also be written using (I.2.1) as

$$\tfrac{d}{dt}\mathfrak{cs} = \int_{\{\cdot\} \times Y} (|E_A|^2 + |\nabla_t\mathfrak{a}|^2)$$

(IV.2.2)

Since $(A, \mathfrak{a})$ obeys the conditions for Proposition IV.1.1, the $t \to \infty$ limit of $\mathfrak{cs}$ converges to some fixed number to be denoted by $\mathfrak{cs}_\infty$.

The crucial proposition for Theorem IX.1 is called a Lojasiewicz inequality:

**Proposition IV.2.1**: *Supposing that $(A, \mathfrak{a})$ is as described above, then there exists numbers $C > 0$ and $\mu \in (0, \frac{1}{2}]$ such that* $|\mathfrak{cs} - \mathfrak{cs}_\infty|^{1-\mu} \leq C ( \int_{\{\cdot\} \times Y} (|E_A|^2 + |\nabla_t \mathfrak{a}|^2))^{1/2}$.

This proposition is proved in Section IV.3. Assume it is true for the moment.

Two observations are needed to exploit this proposition. The first is a direct consequence of (IV.2.2) which is that $\mathfrak{cs}$ is an increasing function of $t$. As a consequence of this $\mathfrak{cs}(t) < \mathfrak{cs}_\infty$ for all $t \geq 1$ unless $E_A$ and $\nabla_{At} \mathfrak{a}$ are identically zero in which case the conclusions of Theorem IV.1 hold automatically. The second observation is this:

$$\tfrac{d}{dt} (\mathfrak{cs}_\infty - \mathfrak{cs})^\mu = -\mu (\mathfrak{cs}_\infty - \mathfrak{cs})^{\mu-1} \int_{\{\cdot\} \times Y} (|E_A|^2 + |\nabla_{At} \mathfrak{a}|^2) \ .$$

(IV.2.3)

Now input Proposition IV.2.1 on the right hand side of (IV.2.3) to obtain the inequality

$$\tfrac{d}{dt} (\mathfrak{cs} - \mathfrak{cs}_\infty)^\mu \leq -\mu C^{-1} ( \int_{\{\cdot\} \times Y} (|E_A|^2 + |\nabla_{At} \mathfrak{a}|^2))^{1/2} \ .$$

(IV.2.4)

Then integrate this last inequality to see that

$$\int_1^\infty ( \int_{\{\cdot\} \times Y} (|E_A|^2 + |\nabla_{At} \mathfrak{a}|^2))^{1/2} < \infty \ .$$

(IV.2.5)

As explained in the next subsection, the fact that the integral on the left side of (IV.2.5) is finite leads to the desired unique convergence. What follows directly is parenthetical.

Another inequality coming from the proposition and (IV.2.6) is this:

$$(\mathfrak{cs}_\infty - \mathfrak{cs})^{2(1-\mu)} \leq C^2 \tfrac{d}{dt} \mathfrak{cs} \ .$$

(IV.2.6)

Since the right hand side of this is $-C^2 \tfrac{d}{dt} (\mathfrak{cs}_\infty - \mathfrak{cs})$, this inequality can be integrated. To do that, first divide both sides by $(\mathfrak{cs}_\infty - \mathfrak{cs})^{2(1-\mu)}$. If $\mu$ is strictly less than $\frac{1}{2}$, then the right hand side is $\tfrac{1}{1-2\mu} C^2$ times the $t$-derivative of $(\tfrac{1}{\mathfrak{cs}_\infty - \mathfrak{cs}})^{1-2\mu}$. Then, integrate to see that

$$\mathfrak{cs}_\infty - \mathfrak{cs} < c \, \frac{1}{t^{1/(1-2\mu)}} \ .$$

(IV.2.7)

when $t \geq 2$. (Note that $c$ depends on $\mu$ and $C$.) If $\mu = \frac{1}{2}$, the inequality in (IV.2.6) says that $\frac{d}{dt}(e^{t/C^2}(\mathfrak{cs}_\infty - \mathfrak{cs})) \leq 0$ which gives an exponential decay rate: $\mathfrak{cs}_\infty - \mathfrak{cs} < c \, e^{-t/C^2}$.

**b) Proof of unique convergence**

This subsection uses (IV.2.5) to prove $(A, \mathfrak{a})$ has a unique $t \to \infty$ limit after the application of a suitable automorphism of P. The proof has five parts.

*Part 1*: Define an isomorphism between $P|_{\{1\} \times Y}$ and $P|_{[1,\infty) \times Y}$ by using the parallel transport by A along the line segments $[1, t] \times \{y\}$ for $y \in Y$. Use this isomorphism and pull-back via the projection from $[1, \infty) \times Y$ to $\{1\} \times Y$ to view $A|_{t=1}$ as a connection on P over $[1, \infty) \times Y$ and likewise view $\mathfrak{a}|_{t=1}$ as a section of $\mathrm{ad}(P) \otimes T^*Y$ over this same domain. Let $\partial_t$ denote the $A|_{t=1}$-covariant derivative along the $[1, \infty)$ factor of $[1, \infty) \times Y$. Now write A as $A|_{t=1} + \mathrm{A}$ and $\mathfrak{a}$ as $\mathfrak{a}|_{t=1} + \mathfrak{A}$. (Both $\mathrm{A}$ and $\mathfrak{A}$ are sections over $[1, \infty) \times Y$ of $\mathrm{ad}(P) \otimes T^*Y$.) The integral in (IV.2.5) is the following

$$\int_1^\infty ( \int_{\{\cdot\} \times Y} (|\partial_t \mathrm{A}|^2 + |\partial_t \mathfrak{A}|^2))^{1/2}$$

(IV.2.8)

because $E_i = \partial_t A_i$ and $\nabla_t \mathfrak{a}_i = \partial_t \mathfrak{A}_i$. Since the preceding integral is finite, the families $\{A|_{\{t\} \times Y}\}_{t \in [1,\infty)}$ and $\{\mathfrak{A}|_{\{t\} \times Y}\}_{t \in [1,\infty)}$ converge as $t \to \infty$ in the $L^2$-topology on the space of $\mathrm{ad}(P)$ valued 1-forms on Y with the limit giving a pair of sections in the $L^2$ completion of $C^\infty(Y; \mathrm{ad}(P) \otimes T^*Y)$. (The convergence can be viewed as taking place in the Hilbert space $\mathbb{L}_Y$ defined previously in Section III.6b.) Now it remains to show that $\mathrm{A}$ and $\mathfrak{A}$ converge is the $C^\infty$ topology.

*Part 2*: What follows describes the strategy for obtaining a priori $C^k$ bounds on $(\mathrm{A}, \mathfrak{A})$ for positive integers k: Consider first the case of $C^0$ bounds for $\mathrm{A}$ and $\mathfrak{A}$. These are bounded on $[1, \infty) \times Y$ if every function on $[1, \infty)$ from the set $\{|E|(\cdot, y) + |\nabla_{A_t} \mathfrak{a}|(\cdot, y)\}_{y \in Y}$ has finite integral on $[1, \infty)$ and if the resulting set of integrals is itself bounded. Meanwhile, it follows from (IV.2.8) that these integrals are finite and that the set of such integrals is bounded if there exists $r > 1$ such that for any $y \in Y$ and $t \geq 2$, one has

$$|E|(t, y) + |\nabla_t \mathfrak{a}|(t, y) \leq r \int_{t-2}^{t+2} ( \int_{\{\cdot\} \times Y} (|E_A|^2 + |\nabla_t \mathfrak{a}|^2))^{1/2} \ .$$

(IV.2.9)

Given $C^0$ bounds for $(A, \mathfrak{A})$, one can then obtain their $C^1$ bounds if there exists some positive $r_1$ such that the inequality

$$|\nabla E|(t,y) + |\nabla \nabla_t \mathfrak{a}|(t,y) \leq r_1 \int_{t-2}^{t+2} ( \int_{\{\cdot\} \times Y} (|E_A|^2 + |\nabla_t \mathfrak{a}|^2))^{1/2} .$$

(IV.2.10)

holds for all $t \in [3, \infty)$ and all $y \in Y$. (In this regard: Keep in mind that the norm of $\nabla_A E$ is bounded by the sum of the norms of $\nabla_{A|_{t=1}} E$ and $[A, E]$. Likewise, the norm of $\nabla \nabla_t \mathfrak{a}$ is bounded by the sum of those of $\nabla_{A|_{t=1}} (\nabla_t \mathfrak{a})$ and $[A, \nabla_t \mathfrak{a}]$.) Continuing in this vein, one can obtain $C^2$ bounds for $(A, \mathfrak{A})$ if the respective norms of $\nabla^{\otimes 2} E$ and $\nabla^{\otimes 2} \nabla_t \mathfrak{a}$ at any $(t,y)$ with $t > 4$ are bounded by what is written on the right hand side of (IV.2.10) for a suitable choice of $r_1$. Having obtained $C^2$ bounds, move next to $C^3$ bounds and so on. In each case, the key step is to bound the norms $\nabla^{\otimes k} E$ and $\nabla^{\otimes k} \nabla_t \mathfrak{a}$ where $t \geq k+2$ by a (possibly k-dependent) constant multiple of the integral on the right hand side of (IV.2.9).

*Part 3*: This part of the proof explains how to derive (IV.2.9), (IV.2.10) and their $\{(\nabla^{\otimes k} E, \nabla^{\otimes k} \nabla_t \mathfrak{a})\}_{k \geq 2}$ analogs. To do this, introduce $\psi$ to denote the pair of ad(P)-valued sections of $T^*Y$ given by $(B = E_A, C = \nabla_{At} \mathfrak{a})$. This $\psi$ will be viewed as a section of the bundle $\oplus_2 (\text{ad}(P) \oplus (\text{ad}(P) \otimes T^*Y))$ over $[1, \infty) \times Y$ with zero for its two ad(P) components. (This is the bundle $\mathbb{W}$ from Section I.1b.) When viewed in this light, the $(A, \mathfrak{a})$ version of $\mathcal{D}$ can act on $\psi$, and in doing so, it annihilates $\psi$: $\mathcal{D}\psi = 0$. It follows as a consequence that $\psi$ is also annihilated by $\mathcal{D}^\dagger \mathcal{D}$. Then, the Bochner-Weitzenboch formula for the latter leads to an equation that has the schematic form

$$\nabla^\dagger \nabla \psi + \mathfrak{F} \psi = 0 ,$$

(IV.2.11)

with $\mathfrak{F}$ denoting an endomorphism of $\mathbb{W}$ whose norm is bounded on $[2, \infty)$ as are the norms of its covariant derivatives to any given finite order. In this regard: Proposition IV.1.1 guarantees that the part of the endomorphism $\mathfrak{F}$ that involve commutators with components of $F_A$ and $\nabla \mathfrak{a}$ and $\mathfrak{a}$ are all uniformly bounded on $[2, \infty) \times Y$ and likewise so are all of their covariant derivatives to any order. Of course, this is also the case for the Riemannian curvature contributions to $\mathfrak{F}$.

Given what Proposition IV.1.1 says, a straightforward use of bump functions with (IV.2.11) and with its consecutive $\{\nabla^{\otimes k}\}_{k=1,2,\ldots}$ covariant derivatives leads to bounds of the following sort: Supposing that $t \geq k+1$, then

$$\int_{[t-\frac{1}{2},t+\frac{1}{2}]\times Y} |\nabla^{\otimes k}\psi|^2 \leq r_k \int_{[t-1,t+1]\times Y} (|E_A|^2 + |\nabla_t \mathfrak{a}|^2)$$

(IV.2.12)

where $r_k$ depends on k but not on t. These bounds lead to sup norm bounds using the Sobolev inequalities: For example, the bounds for k = 1, 2, 3 leads to an $L^4$ bound for $|\nabla^{\otimes 2}\psi|$ on a slightly smaller domain whose square is bounded by a t-independent multiple of the right hand side of (IV.2.12); that in turn leads to an analogous $L^8$ bound for $|\nabla\psi|$ and thus $|d|\psi||$. The latter bound then leads to a $C^0$ bound for $\psi$ on $[t-\frac{1}{4}, t+\frac{1}{4}]\times Y$ via a dimension four Sobolev inequality.

To summarize: Supposing that $m \geq 0$, then the various versions of (IV.2.12) for integers $k \leq m+3$ lead to bounds of this sort:

$$\sup_{\{t\}\times Y} |\nabla^{\otimes m}\psi| \leq r_{m*} \left( \int_{[t-1,t+1]\times Y} (|E_A|^2 + |\nabla_t \mathfrak{a}|^2) \right)^{1/2},$$

(IV.2.13)

with $r_{m*}$ being independent of t if t > m+2.

*Part 4*: This last part of the proof explains why the integral on the right hand side of (IV.2.13) is bounded by a t-independent multiple of the integral on the right hand side of (IV.2.9). (The opposite inequality is a direct consequence of an instance of the Cauchy-Schwarz inequality which bounds the integral on the right hand side of (IV.2.7) by $\sqrt{2}$ times the integral on the right hand side of (IV.2.13).)

To start the explanation, take the inner product of (IV.2.11) with $\psi$ and then integrate the resulting identity over the slices $\{t\}\times Y$ to obtain this inequality:

$$-\frac{d^2}{dt^2} \int_{\{t\}\times Y} |\psi|^2 \leq c_\diamond \int_{\{t\}\times Y} |\psi|^2$$

(IV.2.14)

with $c_\diamond$ being independent of t. Now fix s > 0 and some $t_0 > 1$ and let M denote the maximum of the function $t \to \int_{\{t\}\times Y} |\psi|^2$ on the interval $[t_0-s, t_0+s]$. Now set

$$\phi = \int_{\{t\}\times Y} |\psi|^2 - \tfrac{1}{2} c_\diamond M(t+s-t_0)(t_0+s-t).$$

(IV.2.15)

This function obeys $-\frac{d^2}{dt^2}\phi \leq 0$. As a consequence, its maximum is taken at $t_0-s$ or $t_0+s$:

$$\int_{\{t\}\times Y} |\psi|^2 - \tfrac{1}{2} c_\lozenge M (t+s-t_0)(t_0+s-t) \le \max\{ \int_{\{t_0-s\}\times Y} |\psi|^2 , \int_{\{t_0+s\}\times Y} |\psi|^2 \} .$$

(IV.2.16)

Taking this inequality at the time t where M is realized gives a bound for M if $\tfrac{1}{2} c_\lozenge M s^2 < 1$. In particular, if $s < c_\lozenge^{-1} M^{-1}$, then M must be less than twice the right hand side of (IV.2.16). Supposing this upper bound for s, then it follows that

$$\int_{\{t\}\times Y} |\psi|^2 \le 2( \int_{\{t_0-s\}\times Y} |\psi|^2 + \int_{\{t_0+s\}\times Y} |\psi|^2 )$$

(IV.2.17)

for all $t \in [t_0-s, t_0+s]$. With regards to choosing s: It follows from Proposition IV.1.1 that a positive integer $n_* > 100$ so that (IV.2.17) holds for any $s < \tfrac{1}{n}$ for an integer $n \ge n_*$ and for any choice of $t_0$ from the interval $[2, \infty)$.

With the preceding understood, return now to the integral on the right hand side of (IV.2.13). Having fixed t, then fix $s \in [\tfrac{1}{2n}, \tfrac{1}{n}]$ and for any integer $k \ge 0$, let $I_k$ denote the interval $[t-1+(k-1)s, t-1+(k+1)s]$. Since the set $\{I_k\}_{k=0,\ldots,4n}$ covers $[t-1,t+1]$

$$( \int_{[t-1,t+1]\times Y} (|E_A|^2 + |\nabla_t \mathfrak{a}|^2) )^{1/2} \le \sum_{k=1,\ldots,4n} ( \int_{I_k \times Y} (|E_A|^2 + |\nabla_t \mathfrak{a}|^2) )^{1/2} .$$

(IV.2.18)

With this undersood, use (IV.2.17) to bound each of the $I_k$ intervals. This is the result:

$$( \int_{[t-1,t+1]\times Y} (|E_A|^2 + |\nabla_t \mathfrak{a}|^2) )^{1/2} \le \tfrac{16}{\sqrt{n}} \sum_{k=1,\ldots,4n} ( \int_{\{t-1+ks\}\times Y} (|E_A|^2 + |\nabla_t \mathfrak{a}|^2) )^{1/2} .$$

(IV.2.19)

Now average both sides over the allowed values of s from $[\tfrac{1}{2n}, \tfrac{1}{n}]$ to see that

$$( \int_{[t-1,t+1]\times Y} (|E_A|^2 + |\nabla_t \mathfrak{a}|^2) )^{1/2} \le 32\sqrt{n} \sum_{k=1,\ldots,4n} ( \int_{\tfrac{1}{2}n_*}^{n_*} ( \int_{\{t-1+k(\cdot)\}\times Y} (|E_A|^2 + |\nabla_t \mathfrak{a}|^2))^{1/2} ds )^{1/2} .$$

(IV.2.20)

To finish, note that the right hand side of this is at most $128 n^{3/2} \int_{t-2}^{t+2} ( \int_{\{\cdot\}\times Y} (|E_A|^2 + |\nabla_t \mathfrak{a}|^2))^{1/2}$.

### 3. Proof of Proposition IV.2.1

A proof of Proposition 10.3 can be had by invoking a version of Theorem 3 in [S] along essentially the same lines as was done by Morgan, Mrowka and Ruberman in Section 4.2 of [MMR]. (They state and prove an analog of Proposition IV.2.1 for the anti-self dual Yang-Mills equations, their Proposition 4.2.1) The presentation below paraphrases what is basically the [MMR] story.

To set the notation for what is to come: Let $\mathcal{A}$ denote the set of pairs of the form $(A, \mathfrak{a})$ consisting of a connection on P over Y and ad(P) valued 1-form on Y that obey the constraint $d_A *\mathfrak{a} = 0$. Let $\mathfrak{f}$ denote the function on $\mathcal{A}$ that sends $(A, \mathfrak{a})$ to

$$\int_Y \langle \mathfrak{a} \wedge * B_A - \tfrac{1}{3} \mathfrak{a} \wedge \mathfrak{a} \wedge \mathfrak{a} \rangle .$$

(IV.3.1)

Of particular note is that $\mathfrak{f}$ is an Aut(P)-invariant function on $\mathcal{A}$. The formal gradient of this function $\mathfrak{f}$ is denoted by $\nabla \mathfrak{f}$:

$$\nabla \mathfrak{f} \equiv (*d_A \mathfrak{a}, B_A - *(\mathfrak{a} \wedge \mathfrak{a}))$$

(IV.3.2)

It is a gradient in following sense: If $\mathfrak{b}$ and $\mathfrak{c}$ are any pair of ad(P)-valued 1-forms, then the derivative at the origin of the function on (-1, 1) given by the rule $s \to \mathfrak{f}(A + s\mathfrak{b}, \mathfrak{a} + s\mathfrak{c})$ is equal to

$$\int_Y (\langle \mathfrak{b} \wedge d_A \mathfrak{a} \rangle + \langle \mathfrak{c} \wedge *(B_A - \mathfrak{a} \wedge \mathfrak{a}) \rangle)$$

(IV.3.3)

which is the $\mathbb{L}_Y$-inner product between $(\mathfrak{b}, \mathfrak{c})$ and $\nabla \mathfrak{f}$ when these are viewed as sections of $\mathbb{W}_Y = \oplus_2 (\text{ad}(P) \oplus (\text{ad}(P) \otimes T^*Y))$ with vanishing ad(P) components. (By way of a reminder from Section III.6b: The $\mathbb{L}_Y$-inner product between pairs $(\mathfrak{b}, \mathfrak{c})$ and $(\mathfrak{b}', \mathfrak{c}')$ is the integral over Y of $\langle \mathfrak{b} \wedge * \mathfrak{b}' \rangle + \langle \mathfrak{c} \wedge * \mathfrak{c}' \rangle$. The $\mathbb{L}_Y$-norm of $(\mathfrak{b}, \mathfrak{c})$ is the square root of the Y-integral of $|\mathfrak{b}|^2 + |\mathfrak{c}|^2$.)

The preceding definitions are relevant with regards to Proposition IV.2.1 because:

- *The pull-back of Proposition IV.2.1's pair $(A, \mathfrak{a})$ to any slice $\{t\} \times Y$ is in $\mathcal{A}$*
- *The value of $\mathfrak{cs}(t)$ is the value of $\mathfrak{f}$ on the pull-back of $(A, \mathfrak{a})$ to any $\{t\} \times Y$.*
- *$(E_A, \nabla_{At} \mathfrak{a})$ is $(*d_A \mathfrak{a}, B_A - *(\mathfrak{a} \wedge \mathfrak{a}))$ which is $\nabla \mathfrak{f}$ at any such pull-back of $(A, \mathfrak{a})$.*

(IV.3.4)

As explained momentarily, Proposition IV.2.1 is a consequence of (IV.3.4) and the next proposition which says in effect that $\mathfrak{f}$ near a compact, connected set of Aut(P) orbits of pairs in $\mathcal{A}$ defining flat $SL(2; \mathbb{C})$ connections obeys a bound of the form $|\mathfrak{f} - \mathfrak{f}_*|^{1-\mu} \leq \|\nabla \mathfrak{f}\|_{\mathbb{L}}$ with $\mu \in (0, \tfrac{1}{2}]$ and $\mathfrak{f}_*$ both being constant. (This upcoming proposition is the analog of Theorem 3 from [S] in this context, thus the analog of Proposition 4.2.1 in [MMR].)

By way of terminology, the promised proposition refers to the $\mathbb{H}$-distance between pairs in $\mathcal{A}$ and Aut(P) orbits in $\mathcal{A}$. This $\mathbb{H}$-distance is denoted by $\mathfrak{d}$ and it is defined as follows: Supposing that (A, 𝔞) and (A´, 𝔞´) are from $\mathcal{A}$, introduce ς to denote (A´- A, 𝔞´- 𝔞). The distance ∂ is the square root of the minimum of the Y-integrals of $|\nabla_A ς|^2 + |ς|^2$ and $|\nabla_{A´} ς|^2 + |ς|^2$. The ℍ-distance between (A, 𝔞) and the Aut(P) orbit of the pair (A´, 𝔞´) is defined to be the infimum of the ℍ-distances between (A, 𝔞) and pairs on that Aut(P) orbit. This version of ℍ-distance is also denoted by ∂. The latter is also defined to be the ℍ-distance between their respective Aut(P) orbits.

***Proposition IV.3.1*** *Let $\mathcal{R}$ denote a sequentially compact, connected set of* Aut(P) *equivalence classes of pairs from $\mathcal{A}$ that define flat* SL(2; ℂ) *connections. Given $\mathcal{R}$, there exists $\mu \in (0, \frac{1}{2}]$ and $C > 0$ and $\delta > 0$ and $\mathfrak{f}_* \in \mathbb{R}$ such that if (A, 𝔞) is in $\mathcal{A}$ and has ℍ-distance less than $\delta$ from an orbit in $\mathcal{R}$, then $\mathfrak{f}$ at (A, 𝔞) and the $\mathbb{L}$-norm of $\nabla\mathfrak{f}$ at (A, 𝔞) obey the inequality $|\mathfrak{f} - \mathfrak{f}_*|^{1-\mu} \leq \|\nabla\mathfrak{f}\|_\mathbb{L}$.*

Accept this proposition as true for a moment to derive the assertion of Proposition IV.2.1.

***Proof of Proposition IV.2.1***: The assertion of Proposition IV.2.1 is a direct consequence of this new proposition given what is said in Proposition IV.2.2 about sequential compactness. To elaborate about the role of Proposition IV.2.2: It plays two roles. In the first, it implies this: Given $\varepsilon > 0$, there exists $t_\varepsilon$ such that if $t > t_\varepsilon$, then the pull-back to $\{t\} \times Y$ of the solution (A, 𝔞) to (I.2.1) has ℍ-distance less than $\varepsilon$ from some a pair from $\mathcal{A}$ that defines a flat SL(2; ℂ) connection. (Infact, given a positive integer k, the number $t_\varepsilon$ can be chosen so that (A, 𝔞) on $\{t\} \times Y$ is close in the $C^k$ to a pair from $\mathcal{A}$ that defines a flat SL(2; ℂ) connection.) In the second role, Propostion IV.2.2 implies that the set of Aut(P) orbits of such $t \to \infty$ limit pairs from $\mathcal{A}$ is connected and sequentially compact with respect to ℍ-distance (and with respect to any given integer k version of the $C^k$-distance between pairs in $\mathcal{A}$).

The rest of this section makes various observations (some in the form of lemmas) that are used as input for Section IV.4's proof of Proposition IV.3.1.

**a) On the role of $\mathcal{R}$ in Proposition IV.3.1**

Because the set $\mathcal{R}$ from Proposition IV.3.1 is assumed to be compact, it is sufficient to prove that the assertion of Proposition IV.3.1 holds when $\mathcal{R}$ is replaced by the Aut(P) orbit of any given pair in $\mathcal{A}$ that defines a flat SL(2; ℂ) connection. Indeed, suppose that the assertion does hold in these instances. Then each Aut(P) orbit in $\mathcal{R}$ has an associated set $(\mu, C, \delta)$ which can depend on the particular Aut(P) orbit (but $\mathfrak{f}_*$ would be the same for each such orbit because $\mathcal{R}$ is assumed to be connected). Given an orbit $\mathcal{O} \subset \mathcal{R}$, let $(\mu_\mathcal{O}, C_\mathcal{O}, \delta_\mathcal{O})$ denote the corresponding version of $\mu$, $\delta$ and $C$. Also, let $\mathcal{B}_\mathcal{O}$ denote the set of pairs in $\mathcal{A}$ with $\mathbb{H}$-distance less than $\delta_\mathcal{O}$ from $\mathcal{O}$. The collection of these $\mathcal{O} \subset \mathcal{R}$ versions of $\mathcal{B}_\mathcal{O}$ when mapped to $\mathcal{A}/\text{Aut}(P)$ intersect $\mathcal{R}$ so as to define an open cover of $\mathcal{R}$. Since $\mathcal{R}$ is compact, a finite subset of these $\mathcal{B}_\mathcal{O}$ will do that also. Take this finite subset and define $\mu$ and $\delta$ for Proposition IV.3.1 to be the smallest of the corresponding $\mu_\mathcal{O}$ and $\delta_\mathcal{O}$ with $\mathcal{O}$ taken from the finite subset. Meanwhile, define $C$ for Proposition IV.3.1 to be the largest of the corresponding $C_\mathcal{O}$ with again $\mathcal{O}$ taken from the finite subset.

Section IV.4 proves that the assertion of Proposition IV.3.1 holds when $\mathcal{R}$ is the Aut(P) orbit of any given pair from $\mathcal{A}$ that defines a flat $SL(2;\mathbb{C})$ connection. (This is basically the argument used by Morgan, Mrowka and Ruberman in [MMR] with a proof of Simon's Lojasiewicz inequality in [S] written out for these specific circumstances.)

**b) An Aut(P) slice lemma for Proposition IV.3.1**

The central lemma in this subsection states what is often called a slice lemma for the Aut(P) orbits of pairs in $\mathcal{A}$ near any chosen pair. It identifies these Aut(P) orbits with the quotient of a ball in $C^\infty(Y; \otimes_2(\text{ad}(P) \otimes T^*Y))$ by the action of the stabilizer in Aut(P) of the chosen pair. With regards to this stabilizer: It is the set of automorphisms of P that fix the chosen pair in $\mathcal{A}$. It is isomorphic to a finite subgroup of SU(2). By way of an example: The chosen pair is irreducible if this subgroup is $\{\pm 1\}$.

**Lemma IV.3.2**: *Having fix a pair of connection on P and ad(P) valued 1-form (to be denoted here by $(A, \mathfrak{a})$), there exists a positive number, $\delta$, with the following significance: Suppose that $\mathcal{O}$ is an Aut(P) orbit of a pair of connection on P and ad(P) valued 1-form with $\mathbb{H}$-distance less that $\delta$ from the Aut(P) orbit of $(A, \mathfrak{a})$. Modulo the action of the stabilizer of $(A, \mathfrak{a})$ in Aut(P), there is a unique pair on $\mathcal{O}$ (to be denoted by $(A + \mathfrak{b}, \mathfrak{a} + \mathfrak{c})$) that minimizes the $\mathbb{H}$-distance from $(A, \mathfrak{a})$ to the pairs on $\mathcal{O}$. Moreover, this distance minimizing pair $(\mathfrak{b}, \mathfrak{c})$ obeys the slice identity $(\nabla_{A_i} \mathfrak{b})_i + [\mathfrak{a}_i, \mathfrak{c}_i] = 0$.*

*Proof of Lemma IV.3.2*: The proof that there is a $\delta$ as described above mimics the proof of a similar slice condition for SU(2) connections (the $\mathfrak{a} \equiv 0$ case) which dates from the dawn of time (see, e.g. [D], [FU].) What follows is an outline of a proof for people who have yet to see one.

Having fixed $\delta$ for the moment, choose a pair on $\mathcal{O}'$ which can be written as $(A + \mathfrak{b}', \mathfrak{a} + \mathfrak{c}')$ with $\varsigma = (\mathfrak{b}', \mathfrak{c}')$ obeying

$$\int_Y (|\nabla_A \varsigma|^2 + |\varsigma|^2) < \delta^2 \, .$$

(IV.3.5)

Any other pair on $\mathcal{O}'$ has the form $(A + g^{-1}\mathfrak{b}'g + g^{-1}\nabla_A g, g^{-1}\mathfrak{a}g + g^{-1}\mathfrak{c}'g)$ with g being a section of Aut(P). This being the case, the strategy is to look for $g \in$ Aut(P) near the identity so that $\mathfrak{b} = g^{-1}\mathfrak{b}'g + g^{-1}\nabla_A g$ and $\mathfrak{c} = g^{-1}\mathfrak{a}g - \mathfrak{a} + g^{-1}\mathfrak{c}'g$ obeys the slice identity. In particular, the strategy looks for g of the form $\exp(\sigma)$ with $\sigma$ being a section of ad(P) whose norm is everywhere much less than 1. Then, the slice identity is obeyed if $\sigma$ obeys

$$\nabla_{Ai}(\nabla_A \sigma)_i + [\mathfrak{a}_i, [\mathfrak{a}_i, \sigma]] + \nabla_{Ai}\mathfrak{b}'_i + [\mathfrak{a}_i, \mathfrak{c}'_i] + \mathfrak{k} = 0$$

(IV.3.6)

where $\mathfrak{k}$ is a section of ad(P) obeying:

- $|\mathfrak{k}| \leq c(|\sigma| + |\nabla_A \sigma|)(|\nabla_A \sigma| + |\mathfrak{b}'| + |\mathfrak{c}'|)$ .
- $\int_Y \langle \phi \mathfrak{k} \rangle = 0$ *whenever* $\phi$ *obeys* $\nabla_A \phi \equiv 0$ *and* $[\mathfrak{a}, \phi] \equiv 0$.

(IV.3.7)

The lemma now follows from the following claim with regards to (IV.3.6): If $\delta$ is small, then (IV.3.6) has a unique point-wise small solution obeying the following:

- $\int_Y (|\nabla_A \nabla_A \sigma|^2 + |\nabla_A \sigma|^2 + |[\mathfrak{a}, \sigma]|^2) + \sup_Y |\sigma|^2 \leq c \int_Y (|\nabla_A \varsigma|^2 + |[\mathfrak{a}, \varsigma]|^2)$ .
- $\int_Y \langle \phi \sigma \rangle = 0$ *whenever* $\phi$ *obeys* $\nabla_A \phi \equiv 0$ *and* $[\mathfrak{a}, \phi] \equiv 0$.

(IV.3.8)

The preceding claim can be proved by rewriting (IV.3.6) as an identity that is obeyed by a fixed point of a map from a small radius ball in a certain Hilbert space. This space is the completion of the vector subspace in $C^\infty(Y; \text{ad}(P))$ that obeys the second bullet in (IV.3.8) using the inner product whose associated norm has square given

$$\sigma \to \int_Y (|\nabla_A \nabla_A \sigma|^2 + |\nabla_A \sigma|^2 + |[\mathfrak{a}, \sigma]|^2) \, .$$

(IV.3.9)

Denote this space by $\mathbb{H}^\perp_2$. Let $\mathbb{L}^\perp$ denote the completion of this same vector subspace of $C^\infty(Y; \text{ad}(P))$ using the norm whose square assigns to $\sigma$ the Y-integral of $|\sigma|^2$. The operator $\nabla_{Ai}\nabla_{Ai} + [\mathfrak{a}_i, [\mathfrak{a}_i, \cdot]]$ can be shown to define a bounded, invertible map from $\mathbb{H}^\perp_2$ to $\mathbb{L}^\perp$. Let $\mathcal{K}$ denote its inverse. Meanwhile, the function $\mathfrak{k}$ can be shown to define a smooth map from $\mathbb{H}^\perp_2$ to $\mathbb{L}^\perp$ whose $\mathbb{L}^\perp$ norm is bounded by $c$ times the square of the $\mathbb{H}^\perp_2$ norm. (Note in this regard that the Sobolev inequalities can be used to prove that tautological map from $C^\infty(Y; \text{ad}(P))$ to $C^0(Y; \text{ad}(P))$ extends to $\mathbb{H}^\perp_2$ as a bounded, linear map.) With the preceding understood, it then follows that a solution to (IV.3.6) in $\mathbb{H}^\perp_2$ with small $\mathbb{H}^\perp_2$-norm is a fixed point of the following map (and vice-versa without the small norm constraint):

$$\sigma \to \mathbb{T}(\sigma) = -\mathcal{K}(\nabla_{Ai}\mathfrak{b}'_i + [\mathfrak{a}_i, \mathfrak{c}'_i] + \mathfrak{k}) \:.$$

(IV.3.10)

The fact that $\mathfrak{k}$ is quadratic in $\sigma$ can be used to prove that $\mathbb{T}$ defines a contraction mapping on a ball in $\mathbb{H}^\perp_2$ about the origin of radius $c^{-1}\delta$ if $\delta$ is sufficently small (but positive) and if (IV.3.5) is obeyed by $\varsigma = (\mathfrak{b}', \mathfrak{c}')$.

### c) The operator $\mathfrak{L}$

Any given pair $(A, \mathfrak{a})$ from $\mathcal{A}$ defines a corresponding $\mathbb{H}_Y$-norm on the space of sections of $\mathbb{W}_Y = \oplus_2(\mathrm{ad}(P) \oplus (\mathrm{ad}(P) \otimes T^*Y))$ as follows: The norm is the square root of the functional on the space of sections that sends any given section $\psi$ to the Y-integral of $|\nabla_A\psi|^2 + |\psi|^2$. The completion of $\mathbb{W}_Y$ using the $\mathbb{H}_Y$-norm is a Hilbert space; it is denoted by $\mathbb{H}_Y$. (The norm and inner product depend on the choice of A but the space $\mathbb{H}$ does not.) Also: Remember that the square root of the Y-integral of $|\psi|^2$ defines the $\mathbb{L}_Y$-norm on $\mathbb{W}_Y$; and remember that the completion of the latter space using the $\mathbb{L}_Y$-norm is a Hilbert space that is denoted by $\mathbb{L}_Y$.

With regards to conventions: An unwritten convention in what follows is to view a pair $(\mathfrak{b}, \mathfrak{c})$ from $C^\infty(Y; \oplus_2 (\mathrm{ad}(P) \otimes T^*Y))$ as an element in $\mathbb{H}_Y$ with both $\mathrm{ad}(P)$ components being zero.

More convention/notation: At the risk of introducing confusion with regards to the notation used by Proposition IV.3.1, a chosen pair from $\mathcal{A}$ that defines a flat $SL(2;\mathbb{C})$ connection is denoted subsequently by $(A, \mathfrak{a})$. (Just to be sure, the $SL(2;\mathbb{C})$ connection $\mathbb{A} = A + i\mathfrak{a}$ is henceforth assumed to be flat.) Nearby pairs in $\mathcal{A}$ will be denoted by $(A', \mathfrak{a}')$. These need not define flat $SL(2;\mathbb{C})$ connections. This notation change is introduced because most of the subsequent discussion centers around the pair that defines the flat $SL(2;\mathbb{C})$ connection.

Let $(A, \mathfrak{a})$ now denote a pair from $\mathcal{A}$ with $A + i\mathfrak{a}$ being a flat $SL(2;\mathbb{C})$ connection. Having specified $(A, \mathfrak{a})$, return to the $\gamma$ and $\rho$ matrix notation from Section I.1a and let $\mathfrak{L}$ denote the operator $\gamma^i\nabla_{Ai} + \rho_i[\mathfrak{a}_i, \cdot]$. This operator defines a bounded, symmetric (essentially self-adjoint) operator mapping $\mathbb{H}_Y$ to $\mathbb{L}_Y$. As such, it is Fredholm with index zero. Viewed as an unbounded operator on $\mathbb{L}_Y$, it is self-adjoint with dense domain $\mathbb{H}_Y$. It has a complete orthonormal set of eigensections in this incarnation (which are a priori smooth) with the corresponding set of eigenvalues being a set of real numbers with no accumulation points. Also, each eigenvalue has finite multiplicity.

Let $\Pi^0$ denote the $\mathbb{L}_Y$-orthogonal projection to the kernel of $\mathfrak{L}$ (which is finite dimensional). The image of $\mathfrak{L}$ is $\mathbb{L}_Y$-orthogonal to kernel($\mathfrak{L}$) by virtue of $\mathfrak{L}$ being symmetric. It therefore has bounded inverse as a map from $(1-\Pi^0)\mathbb{H}_Y$ to $(1-\Pi^0)\mathbb{L}_Y$. With regards to the kernel of $\mathfrak{L}$: This vector space will be written as the $\mathbb{L}_Y$-orthogonal direct sum $H^1 \oplus H^0$ where $H^1$ consists of elements in $\mathbb{H}_Y$ with both ad(P) components zero, and where $H^0$ consists of the elements with only ad(P) components that are annihilated by $\nabla_A(\cdot) = 0$ and $[\mathfrak{a}, \cdot]$. (The existence of this decomposition follows from (III.3.15)–(III.3.17) and what is said about them in Part 7 of Section III.3b.)

**d) The map $\mathfrak{w}$**

This subsection uses $\mathfrak{L}$ to construct a *Kuranishi* picture of a neighborhood of a pair from $\mathcal{A}$ that defines a flat $Sl(2;\mathbb{C})$ connections. To this end, suppose that $(A, \mathfrak{a})$ is such a pair just as in the previous section. Letting $\psi = ((\mathfrak{b}_t, \mathfrak{b}), (\mathfrak{c}_t, \mathfrak{c}))$ denote an element from $\mathbb{H}_Y$, set

- $\mathfrak{p} = -d_A \mathfrak{b}_t - *d_A \mathfrak{c} + *(\mathfrak{b} \wedge \mathfrak{a} + \mathfrak{a} \wedge \mathfrak{b}) + [\mathfrak{a}_i, \mathfrak{c}_t] - [\mathfrak{b}, \mathfrak{b}_t] - *(\mathfrak{b} \wedge \mathfrak{c} + \mathfrak{c} \wedge \mathfrak{b}) + [\mathfrak{c}, \mathfrak{c}_t]$ .
- $\mathfrak{p}_t = (\nabla_i \mathfrak{b})_i + [\mathfrak{a}_i, \mathfrak{c}_i]$ .
- $\mathfrak{q} = -d_A \mathfrak{c}_t - *d_A \mathfrak{b} + *(\mathfrak{c} \wedge \mathfrak{a} + \mathfrak{a} \wedge \mathfrak{c}) - [\mathfrak{a}_i, \mathfrak{b}_t] - [\mathfrak{b}, \mathfrak{c}_t] - *(\mathfrak{b} \wedge \mathfrak{b} - \mathfrak{c} \wedge \mathfrak{c}) - [\mathfrak{c}, \mathfrak{b}_t]$ .
- $\mathfrak{q}_t = (\nabla_i \mathfrak{c}_i) - [\mathfrak{a}_i, \mathfrak{c}_i] + [\mathfrak{b}_t, \mathfrak{c}_t] + [\mathfrak{b}_i, \mathfrak{c}_i]$ .

(IV.3.11)

Of particular note is when $\mathfrak{b}_t$ and $\mathfrak{c}_t$ are zero, then $\mathfrak{p}$ and $\mathfrak{q}$ are the components of $-\nabla \mathfrak{f}$ at $(A' = A + \mathfrak{b}, \mathfrak{a}' = \mathfrak{a} + \mathfrak{c})$; and $\mathfrak{q}_t$ is $*d_{A'}*\mathfrak{a}'$; and $\mathfrak{p}_t$ is what has to vanish if $(\mathfrak{b}, \mathfrak{c})$ obeys Lemma IV.3.2's slice condition.

With regards to the operator $\mathfrak{L}$: The right hand side of (IV.3.11) depicts a smooth map from $\mathbb{H}_Y$ to $\mathbb{L}_Y$ to be denoted by $\mathfrak{F}$ that has the schematic form

$$\psi \to \mathfrak{F}(\psi) = \mathfrak{L}\psi + \psi \# \psi$$

(IV.3.12)

where $\psi \# \psi$ signifies the element in $\mathbb{L}_Y$ with its $\mathbb{W}_Y$ components given by the commutators between the various components of $\psi$ that appear in (IV.3.11). (By way of a look ahead: The map $\mathfrak{F}$ is relevant by virtue of its relation to $-\nabla \mathfrak{f}$ and $d_{A'}*\mathfrak{a}'$ and Lemma IV.3.2 when $\mathfrak{b}_t$ and $\mathfrak{c}_t$ are absent.)

The next lemma describes a construction that will be used to exploit the fact that the map $\mathfrak{F}$ is the exactly the operator $\mathfrak{L}$ to leading order and thus linear to leading order on the $\mathbb{L}_Y$-orthogonal complement of $\mathfrak{L}$'s kernel.

**Lemma IV.3.3**: *Let* $(A, \mathfrak{a})$ *denote a pair from* $\mathcal{A}$ *such that* $A + i\mathfrak{a}$ *is flat. There exists* $\delta > 0$ *and* $\kappa > 1$ *and a real analytic map (denoted by* $\mathfrak{w}$*) from the radius* $\delta$ *ball about the origin in* $\mathbb{H}^1$ *to* $(1 - \Pi^0)C^\infty(Y; \mathbb{W}_Y)$ *with the following significance: Let* B *denote this radius* $\delta$*-ball in* $\mathbb{H}^1$. *If* $\phi \in B$, *then*

- $(1 - \Pi^0)\mathfrak{F}(\phi + \mathfrak{w}(\phi)) = 0$.
- $\|\mathfrak{w}(\phi)\|_\mathbb{H} \leq \kappa \|\phi\|_\mathbb{H}^2$.
- $\|\partial_\phi \mathfrak{w}\|_\mathbb{H} \leq \kappa \|\phi\|_\mathbb{H}$ *where* $\partial_\phi$ *denotes the exterior derivative on* B.

*Proof of Lemma IV.3.3*: The existence and analyticity of $\mathfrak{w}$ can be proved using a contraction mapping argument. To set up the argument, introduce by way of notation $\mathfrak{L}^{-1}$ to denote the inverse of the restriction of $\mathfrak{L}$ as a map from $(1 - \Pi^0)\mathbb{H}_Y$ to $(1 - \Pi^0)\mathbb{L}_Y$. Having fixed $\phi$ in the $\mathbb{H}^1$ summand of the kernel of $\mathfrak{L}$, use $\phi$ with $\mathfrak{L}^{-1}$ to define a map to be denoted by $G_\phi$ from $(1 - \Pi^0)\mathbb{H}_Y$ to itself by the rule

$$w \to G_\phi(w) = -\mathfrak{L}^{-1}(1 - \Pi^0)((\phi + w) \# (\phi + w))$$

(IV.3.13)

By construction: A fixed point of this non-linear map $G_\phi$ is a point $w \in (1 - \Pi^0)\mathbb{H}_Y$ that obeys $w = G_\phi(w)$. A fixed point is therefore a solution to the equation in the top bullet of the lemma, and vice-versa.

To find a fixed point: Fix $\rho > 0$ and let $\mathbb{B}_\rho$ denote the ball of radius $\rho$ about the orgin in $(1 - \Pi^0)\mathbb{H}_Y$. The map $G_\phi$ is said to be a contraction mapping on $\mathbb{B}_\rho$ if it maps $\mathbb{B}_\rho$ to to some slightly smaller radius concentric ball in $(1 - \Pi^0)\mathbb{H}_Y$ and if there exists some positive number $\varepsilon$ such that

$$\|G_\phi(w) - G_\phi(w')\|_\mathbb{H} \leq (1 - \varepsilon)\|w - w'\|_\mathbb{H}$$

(IV.3.14)

when every $w$ and $w'$ are both from $\mathbb{B}_\rho$. If it is a contraction mapping on $\mathbb{B}_\rho$, then there is a unique fixed point. This is because the sequence

$$\{w_0 = 0, w_1 = G_\phi(w_0), \ldots, w_k = G_\phi(w_{k-1}), \ldots\}$$

(IV.3.15)

is guaranteed by (IV.3.14) to be a Cauchy sequence in $\mathbb{B}_\rho$ whose limit is a fixed point of $G_\phi$. Meanwhile, if $w$ and $w'$ are fixed points of $G_\phi$ in $\mathbb{B}_\rho$, then (IV.3.14) guarantees that they are equal.

To prove that (IV.3.14) holds when $\phi$ has small $\mathbb{H}_Y$-norm: A Sobolev inequality bounds the Y-integral of $|\eta|^4$ for $\eta \in \mathbb{H}_Y$ by a $\eta$-independent multiple of $\|\psi\|_{\mathbb{H}}^4$. As a consequence, there exists $c_* \geq 1$ such that

$$\| \mathfrak{L}^{-1}(1-\Pi^0)(\eta \# \xi) \|_{\mathbb{H}} \leq c_* \|\eta\|_{\mathbb{H}} \|\xi\|_{\mathbb{H}}$$
(IV.3.16)

for any pair $\eta, \xi$ from $\mathbb{H}_Y$. This implies that the map $\psi \to \mathfrak{L}^{-1}(1-\Pi^0)(\psi \# \psi)$ sends the radius $\rho$ ball in $\mathbb{H}_Y$ into the radius $c_* \rho^2$ ball. In particular, if $\rho < \frac{1}{128} c_*$, and if $\|\phi\|_{\mathbb{H}} \leq \rho$, then $G_\phi$ maps $\mathbb{B}_\rho$ into $\mathbb{B}_{\rho/64}$. The fact that $G_\phi$ is a contraction on these small $\rho$ balls follows from (IV.3.16) also because

$$\mathfrak{L}^{-1}(1-\Pi^0)(\eta \# \eta) - \mathfrak{L}^{-1}(1-\Pi^0)(\xi \# \xi)$$
(IV.3.17)

can be written as

$$\mathfrak{L}^{-1}(1-\Pi^0)(\eta \# (\eta - \xi)) - \mathfrak{L}^{-1}(1-\Pi^0)((\xi - \eta) \# \xi) \,.$$
(IV.3.18)

Supposing that $\|\phi\|_{\mathbb{H}} < \rho$ (with $\rho < \frac{1}{128} c_*^{-1}$), let $\mathfrak{w}(\phi)$ denote the unique small normed fixed point of $G_\phi$. By virtue of (IV.3.16), the norm of $\mathfrak{w}$ obeys

$$\|\mathfrak{w}(\phi)\|_{\mathbb{H}} \leq c_* \|\phi\|_{\mathbb{H}}^2 + c c_* \rho \|\mathfrak{w}(\phi)\| \,.$$
(IV.3.19)

This leads to the bound in the second bullet of the lemma when $\rho < \frac{1}{128} c^{-1} c_*^{-1}$.

To see about the third bullet's bound, suppose for the moment that $\phi$ and $\phi'$ are elements in the radius $\rho < \frac{1}{128} c^{-1} c_*^{-1}$ ball centered on the origin in $H^1$. Set $\psi = \phi + \mathfrak{w}(\phi)$ and set $\psi' = \phi' + \mathfrak{w}(\phi')$. By virtue of the fixed point condition,

$$\mathfrak{w}(\phi') - \mathfrak{w}(\phi) = -\mathfrak{L}^{-1}(1-\Pi^0)((\psi' - \psi) \# \psi') - \mathfrak{L}^{-1}(1-\Pi^0)(\psi \# (\psi' - \psi)) \,,$$
(IV.3.20)

which implies that

$$\|\mathfrak{w}(\phi') - \mathfrak{w}(\phi)\|_{\mathbb{H}} \leq c \|\psi - \psi'\|_{\mathbb{H}} (\|\psi\|_{\mathbb{H}} + \|\psi'\|_{\mathbb{H}}) \,,$$
(IV.3.21)

Supposing that $\rho < c^{-1} c_*^{-1}$ so that $\|\psi\|_{\mathbb{H}} + \|\psi'\|_{\mathbb{H}} \leq c^{-1}$, then this leads directly to the bound

$$\|\mathfrak{w}(\phi') - \mathfrak{w}(\phi)\|_{\mathbb{H}} \leq c \|\phi - \phi'\|_{\mathbb{H}} ((\|\phi\|_{\mathbb{H}} + \|\phi'\|_{\mathbb{H}})$$
(IV.3.22)

which leads in turn to the inequality in the lemma's third bullet.

To see why the map $\phi \to \mathfrak{w}(\phi)$ is real analytic on a positive radius ball in $H^1$: The proceding analysis can be repeated almost verbatim with $H^1$ and $\mathbb{H}_Y$ replaced by their respective complexifications (these are denoted by $H^1_{\mathbb{C}}$ and $\mathbb{H}_{Y\mathbb{C}}$ respectively). The inequality in (IV.3.16) and (IV.3.19) and (IV.3.22) still hold because $G_\phi$ with $\phi$ from a small radius ball in $H^1_{\mathbb{C}}$ defines a quadratic polynomial mapping from the radius $\rho$ ball in $(1-\Pi^0)\mathbb{H}_{Y\mathbb{C}}$ to itself. (It is a sum of a term that is independent of w, a term that is $\mathbb{C}$-linear in w, and a term that is quadratic over $\mathbb{C}$.) It then follows directly that this mapping is a contraction mapping on a small radius ball about the origin in $(1-\Pi^0)\mathbb{H}_{Y\mathbb{C}}$. As a consequence, it has a unique fixed point in this ball for each choice of $\phi$ from a small radius ball about the origin in $H^1_{\mathbb{C}}$. The complex versions of (IV.3.17) and (IV.3.18) can be used to prove that the corresponding complex fixed point depends holomorphically on the choice of $\phi$ when the latter is from a small radius ball in $H^1_{\mathbb{C}}$. This implies that the map $\mathfrak{w}$ is real analytic on the real slice of $H^1_{\mathbb{C}}$.

### e) The $\mathbb{L}$-norm of $\mathfrak{F}$

Suppose once again that $(A, \mathfrak{a}) \in \mathcal{A}$ is such that $A + i\mathfrak{a}$ is a flat $SL(2;\mathbb{C})$ connection. Fix a pair of ad(P)-valued 1-forms on Y to be denoted by $\mathfrak{b}$ and $\mathfrak{c}$. Assume that this pair obeys the slice constraint $\nabla_{Ai}\mathfrak{b}_i + [\mathfrak{a}_i, \mathfrak{c}_i] = 0$. Let $\psi$ denote the pair $(\mathfrak{b}, \mathfrak{c})$, henceforth viewed as an element in $\mathbb{H}_Y$ with zero ad(P) components. Let $\phi$ denote $\Pi^0\psi$. Noting that $\|\phi\|_{\mathbb{H}} \leq c\|\psi\|_{\mathbb{H}}$: If the $\mathbb{H}_Y$ norm of $\psi$ is bounded by $c^{-1}\delta$ with $\delta$ from the $(A, \mathfrak{a})$ version of Lemma IV.3.3, then the $\mathbb{H}_Y$ norm of $\phi$ will be bounded by $\frac{1}{100}\delta$. Assuming this, then Lemma IV.3.3 can be invoked to define $\mathfrak{w}(\phi)$. Then, writing $\psi$ as $\psi \equiv \phi + \mathfrak{w}(\phi) + \eta$ with $\Pi_0\eta = 0$, it follows from the top bullet of Lemma IV.3.3 that

$$(1-\Pi^0)\mathfrak{F} = \mathfrak{L}\eta + (1-\Pi^0)(\eta \# \eta + \eta \#(\phi + \mathfrak{w}) + (\phi + \mathfrak{w}) \# \eta).$$

(IV.3.23)

(Here, $\mathfrak{w}$ is shorthand for $\mathfrak{w}(\phi)$.) This implies in particular that there exists $c_\diamond > 1$ such that if $\|\psi\|_{\mathbb{H}} \leq c_\diamond^{-1}$, then

$$\|(1-\Pi^0)\mathfrak{F}\|_{\mathbb{L}} \geq c_\diamond^{-1}\|\eta\|_{\mathbb{H}}.$$

(IV.3.24)

As for the $\mathbb{L}$-norm of $\Pi^0\mathfrak{F}$, well it obeys the bound

$$\|\Pi^0\mathfrak{F}\|_{\mathbb{L}} \geq \tfrac{1}{2}\|\Pi^0((\phi+\mathfrak{w})\#(\phi+\mathfrak{w}))\|_{\mathbb{L}} - c\|\psi\|_{\mathbb{H}}\|\eta\|_{\mathbb{H}}.$$

(IV.3.25)

Thus, if $\|\psi\|_{\mathbb{H}} \leq c^{-1} c_\diamond^{-1}$, then (IV.3.24) and (IV.3.25) together imply that

$$\|\mathfrak{F}\|_{\mathbb{L}} \geq \tfrac{1}{2} \|\Pi^0((\phi+\mathfrak{w}) \# (\phi+\mathfrak{w}))\|_{\mathbb{L}} + \tfrac{1}{2} c_\diamond^{-1} \|\eta\|_{\mathbb{H}} \ .$$

(IV.3.26)

This last inequality plays a crucial role in the proof of Proposition IV.3.1.

## 4. Proof of Proposition IV.3.1

To start the proof: Let $(A, \mathfrak{a})$ denote a pair from $\mathcal{A}$ with $A + i\mathfrak{a}$ being a flat $SL(2; \mathbb{C})$ connection. Define a function to be denoted by $\mathfrak{f}_\diamond$ on $C^\infty(Y; \oplus_2(ad(P) \otimes T^*Y))$ by the rule whereby $\mathfrak{f}_\diamond(\mathfrak{b}, \mathfrak{c}) = \mathfrak{f}(A+\mathfrak{b}, \mathfrak{a}+\mathfrak{c}) - \mathfrak{f}(A, \mathfrak{a})$ with $\mathfrak{f}$ being Proposition IV.3.1's function (see (IV.3.1)). To be explicit, the value of $\mathfrak{f}_\diamond(\mathfrak{b}, \mathfrak{c})$ is given by this integral:

$$\int_Y \langle \mathfrak{c} \wedge d_A \mathfrak{b} + \mathfrak{a} \wedge \mathfrak{b} \wedge \mathfrak{b} - \mathfrak{a} \wedge \mathfrak{c} \wedge \mathfrak{c} + \mathfrak{c} \wedge \mathfrak{b} \wedge \mathfrak{b} - \tfrac{1}{3}\mathfrak{c} \wedge \mathfrak{c} \wedge \mathfrak{c} \rangle \ .$$

(IV.4.1)

Suppose now that $\psi$ is a given element from $\mathbb{H}_Y$ with vanishing ad(P) components whose norm is small enough so that there is the decomposition $\psi = \eta + \phi + \mathfrak{w}(\phi)$ with $\phi$ being the $H^1$ part (with respect to $\mathbb{L}_Y$-orthogonal projection) and with both $\eta$ and $\mathfrak{w}$ being $\mathbb{L}_Y$-orthogonal to the kernel of $\mathfrak{L}$. (As noted above, this is always doable if $\|\psi\|_{\mathbb{H}} \leq c^{-1}$.) With the preceding understood, write $\mathfrak{f}_\diamond(\psi)$ using Taylor's theorem with remainder as

$$\mathfrak{f}_\diamond(\psi) = \mathfrak{f}_\diamond(t(\phi+\mathfrak{w})) + \nabla \mathfrak{f}_\diamond|_{t(\phi+\mathfrak{w})} \bullet (t\eta) + \mathfrak{r},$$

(IV.4.2)

where $t(\cdot)$ indicates the $\oplus_2(adP \otimes T^*Y)$ part of the indicated element in $\mathbb{H}_Y$ and where $\mathfrak{r}$ is a term with norm bound $|\mathfrak{r}| \leq c \|t\eta\|_{\mathbb{H}}^2$ when $\psi$ has small $\mathbb{H}_Y$-norm (which again means that $\|\psi\|_{\mathbb{H}} \leq c^{-1}$). Meanwhile, Taylors theorem also allows for writing

$$\nabla \mathfrak{f}_\diamond|_{t(\phi+\mathfrak{w})} = \nabla \mathfrak{f}_\diamond|_\psi + \mathfrak{e},$$

(IV.4.3)

where $\mathfrak{e}$ obeys the bound $|\mathfrak{e}| \leq c \|t\eta\|_{\mathbb{H}}$ when $\|\psi\|_{\mathbb{H}} \leq c^{-1}$. Then, (IV.4.2) and (IV.4.3) together lead to this:

$$|\mathfrak{f}_\diamond(\psi) - \mathfrak{f}_\diamond(t(\phi+\mathfrak{w}))| \leq c (\|\nabla \mathfrak{f}_\diamond|_\psi\|_{\mathbb{L}} \|\eta\|_{\mathbb{H}} + \|\eta\|_{\mathbb{H}}^2)$$

(IV.4.4)

when $\|\psi\|_{\mathbb{H}} \leq c^{-1}$ (which is assumed henceforth).

To continue the analysis, let $\wp$ denote the function on the radius $c^{-1}$ ball in $\mathbb{H}^1$ that is given by the rule $\phi \to \wp(\phi) \equiv \mathfrak{f}_\diamond(\mathfrak{r}(\phi + \mathfrak{w}(\phi)))$. This is a real analytic function (see Lemma IV.3.2) on a finite dimensional ball. As such, the classical Lojasiewicz inequality can be invoked (see [L1], [L2]) to find numbers $c_\ddagger > 1$ and $\mu \in (0, \tfrac{1}{2}]$ such that

$$|\wp|^{1-\mu} \leq c_\ddagger |\nabla \wp| \qquad (IV.4.5)$$

on the radius $c_\ddagger^{-1}$ ball about the origin in $\mathbb{H}^1$. Meanwhile, $\nabla \wp$ can be computed in terms of $\nabla \mathfrak{f}_\diamond$ using the chain rule; and what with the third bullet of Lemma IV.3.2, doing that leads to the following inequality

$$|\mathfrak{f}_\diamond(\phi + \mathfrak{w}(\phi))|^{1-\mu} \leq c\, c_\ddagger \|\nabla \mathfrak{f}_\diamond|_{\mathfrak{r}(\phi+\mathfrak{w})}\|_{\mathbb{L}}, \qquad (IV.4.6)$$

which leads via (IV.4.3) to the inequality

$$|\mathfrak{f}_\diamond(\phi + \mathfrak{w}(\phi))|^{1-\mu} \leq c\, (\|\nabla \mathfrak{f}_\diamond|_\psi\|_{\mathbb{L}} + \|\eta\|_{\mathbb{H}}). \qquad (IV.4.7)$$

Putting the preceding bound into (IV.4.4):

$$|\mathfrak{f}_\diamond(\psi)| \leq c\, (\|\nabla \mathfrak{f}_\diamond|_\psi\|_{\mathbb{L}} + \|\eta\|_{\mathbb{H}})^{1/(1-\mu)} + c\, (\|\nabla \mathfrak{f}_\diamond|_\psi\|_{\mathbb{L}} \|\eta\|_{\mathbb{H}} + \|\eta\|_{\mathbb{H}}^2), \qquad (IV.4.8)$$

assuming that $\|\psi\|_{\mathbb{H}} \leq c^{-1}$.

With (IV.4.8) in hand, suppose now that $\psi$ has no ad(P) components so it can be written as $(\mathfrak{b}, \mathfrak{c})$ with both $\mathfrak{b}$ and $\mathfrak{c}$ denoting ad(P)-valued 1-forms on Y. With regards to the latter: Assume that the pair $(A + \mathfrak{b}, \mathfrak{a} + \mathfrak{c})$ is in $\mathcal{A}$ (so that $d_{A+\mathfrak{b}}*(\mathfrak{a} + \mathfrak{c}) = 0$) and that Lemma IV.3.1's slice condition $\nabla_{A_i} \mathfrak{b}_i + [\mathfrak{a}_i, \mathfrak{c}_i] = 0$. In this event, $\nabla \mathfrak{f}_\diamond|_\psi$ is the map $\mathfrak{F}$ that is depicted in (4.3.12). (See (IV.3.11) with $\mathfrak{b}_t \equiv \mathfrak{c}_t \equiv 0$.) For $\psi$ of this sort, the inequality in (IV.3.26) can be employed with (IV.4.8)'s inequality to see that

$$|\mathfrak{f}_\diamond(\psi)| \leq c\, (\|\nabla \mathfrak{f}_\diamond|_\psi\|_{\mathbb{L}}^{1/(1-\mu)} + \|\nabla \mathfrak{f}_\diamond|_\psi\|_{\mathbb{L}}^2) \qquad (IV.4.9)$$

when $\|\psi\|_{\mathbb{H}} \leq c^{-1}$ which is Proposition IV.3.1's claim because $\tfrac{1}{1-\mu} \leq \tfrac{1}{2}$.

## 5. Exponential decay as t → ∞

To set the stage for what is to come, let $(A, \mathfrak{a})$ denote a pair of connection on P and T*Y valued 1-form, both defined over $[1, \infty) \times Y$. Assume that this pair obeys (I.2.1) on the interior of this domain. With regards to the $t \to \infty$ asymptotics: Identify P with the pull-back via the projection map to Y of a principal SU(2) bundle over Y (also denoted by P); and let $\mathbb{A} = A_\infty + i\mathfrak{a}_\infty$ denote a flat, $SL(2;\mathbb{C})$ connection on $P|_Y \times_{SU(2)} SL(2;\mathbb{C})$. Use the projection map from $[1,\infty) \times Y$ to view $(A_\infty, \mathfrak{a}_\infty)$ simultaneously as a pair of connection on P over $[1,\infty) \times Y$ and section over this same domain of $ad(P) \otimes T^*Y$.

The lemma that follows restates the assertion in Theorem IV.1 with regards to the $t \to 0$ exponential convergence.

**Lemma IV.5.1**: *With $(A, \mathfrak{a})$ and $\mathbb{A}$ as just described, suppose the following: Given $\varepsilon > 0$, there exists $t_\varepsilon$ and and an automorphism of P on $[t_\varepsilon, \infty) \times Y$ that pulls $(A,\mathfrak{a})$ back to give a pair that differs pointwise by at most $\varepsilon$ from $(A_\infty, \mathfrak{a}_\infty)$. Assuming this, there exists $\delta > 0$ that depends only on the $\mathrm{Aut}(P)$ orbit of $(A_\infty, \mathfrak{a}_\infty)$, and there exists $\kappa > 1$ and an automorphism of P over $[t_\varepsilon, \infty) \times Y$ with the following significance: The pull-back of $(A, \mathfrak{a})$ by the automorphism obeys $|A^\infty - A| + |\mathfrak{a}^\infty - \mathfrak{a}| < \kappa e^{-\delta t}$. Moreover there are pointwise bounds proportional to $e^{-\delta t}$ for the $A_\infty$-covariant derviatives to any given order of the both $A^\infty - A$ and $\mathfrak{a}^\infty - \mathfrak{a}$.*

*Proof of Lemma IV.5.1*: Write the curvature 2-form of A as $dt \wedge E_A + *B_A$ with $E_A$ and $B_A$ being $ad(P)$ valued sections of $T^*Y$ over $[1,\infty) \times Y$. Granted this notation, let $\psi$ denote the pair $(\mathfrak{b} = E_A, \mathfrak{c} = \nabla_t \mathfrak{a})$. This element $\psi$ obeys $\mathcal{D}\psi = 0$ with $\mathcal{D}$ defined by $(A, \mathfrak{a})$ and as depicted in Section 1a.

Having fixed $\varepsilon > 0$, use it to define the time $t_\varepsilon$ and the automorphism of P that makes $(A, \mathfrak{a})$ differ by at most $\varepsilon$ from $(A_\infty, \mathfrak{a}_\infty)$ on $[t_\varepsilon, \infty) \times Y$. This automorphism writes the operator $\mathcal{D}$ as $\mathcal{D}^\infty + \mathfrak{k}$ with $\mathcal{D}^\infty$ denoting the $(A_\infty, \mathfrak{a}_\infty)$ version of (I.1.2) and with $|\mathfrak{k}| \leq c\varepsilon$. Thus, the pull-back of $\psi$ by this automorphism (still denoted by $\psi$) obeys

$$\mathcal{D}^\infty \psi + \mathfrak{k}\psi = 0 \ .$$

(IV.5.1)

Let $\Pi^+$ and $\Pi^-$ denote the repective $\mathbb{L}_Y$-orthogonal projections on Y to the span of the eigenvectors of $\mathcal{D}^\infty$ with respective positive and negative eigenvalues. It is important to note that $\psi$ is equal to the sum of $\Pi^+\psi$ and $\Pi^-\psi$ because $\psi$ when written as in (I.1.5) has no $\mathfrak{c}_t$ or $\mathfrak{b}_t$ components whereas the zero eigenvalue of $\mathcal{D}^\infty$ has only these components. (There is no zero eigenvalue at all if $\mathbb{A}$ is irreducible.) Let $f_+$ and $f_-$ denote the respective respective functions on $[t_\varepsilon, \infty)$ whose values are the $\mathbb{L}_Y$-norms on $\{t\} \times Y$ of $\Pi^+\psi$ and $\Pi^-\psi$. These functions obey (by virtue of (IV.5.1)) the analogs of (III.6.5) that are written below when $\varepsilon < c^{-1}$ (as in (III.6.5), the inequalities use $\lambda$ to denote half of the smallest of the absolute values of the non-zero eigenvalues of the operator $\gamma_i \nabla_{A^\infty_i} + \rho_i[\mathfrak{a}^\infty_i, \cdot\,])$.

- $\frac{d}{dt} f_+ + \lambda f_+ \leq c\varepsilon f_-$.
- $\frac{d}{dt} f_- - \lambda f_- \geq -c\varepsilon f_+$.

(IV.5.2)

Now fix $z > 1$ for the moment. The preceding inequalities imply that

$$\tfrac{d}{dt}(f_+ - z\varepsilon f_-) + (\lambda - cz\varepsilon^2)f_+ + (z\lambda - c)\varepsilon f_- < 0 .$$

(IV.5.3)

Thus, if $z > c$ and $\varepsilon < c^{-1}z^{-2}$, then

$$\tfrac{d}{dt}(f_+ - z\varepsilon f_-) + \tfrac{1}{2}\lambda(f_+ - z\varepsilon f_-) < 0 .$$

(IV.5.4)

This implies in turn that either $f_+ < -z\varepsilon f_-$ on $[t_\varepsilon, \infty)$ or that $(f_+ - z\varepsilon f_-) \geq e^{\lambda t/4}$ where t is large. Since the latter event runs afoul of the a priori bounds on curvatures in Proposition IV.1.1, it must be that that $f_+ < -z\varepsilon f_-$. Granted this, and granted that $\varepsilon < c^{-1}z^{-2}$, then the second bullet in (IV.5.2) can be integrated to see that $f_- \leq e^{-\lambda t/4}$ where t is large; and thus $f_+ \leq z\varepsilon\, e^{-\lambda t/4}$ also where t is large.

These $\mathbb{L}_Y$-norms for $\nabla_t \mathfrak{a}$ and $E_t$ can now be used as input for the analysis in Section IV.2b to complete the proof of Lemma IV.5.1.

## 6. When the $\{t\} \times Y$ integrals of $|\mathfrak{a}|^2$ diverge

This short section states and then proves a proposition about instances when the assumptions in Theorem IV.1 are not met with regards to the $\{t\} \times Y$ integrals of $|\mathfrak{a}|^2$.

To set the stage for this proposition, suppose that $(A, \mathfrak{a})$ is a solution to (I.2.1) on $[1, \infty) \times Y$ with the two properties listed below. (The second bullet refers to the function $\mathfrak{cs}$ on $[1, \infty)$ that is defined in (IV.2.1).)

- *Given $x > 1$, there exists $t_x \geq 1$ such that $\int_{\{t\} \times Y} |\mathfrak{a}|^2 \geq x$ when $t \geq t_x$.*

- *The function $\mathfrak{cs}$ is bounded.*

(IV.6.1)

By way of notation for the upcoming proposition, let M denote the positive function on $[1, \infty)$ whose value at time t is the square root of the $\{t\} \times Y$ integral of $|\mathfrak{a}|^2$.

By way of more notation for this proposition: A $\mathbb{Z}/2$-harmonic 1-form data set on Y consists of a triple $(Z, \mathcal{I}, \nu)$ with these being as follows:

- *Z is a closed subset in Y with finite 1-dimensional Hausdorff dimension,*

- $\mathcal{I}$ *is a real line bundle over* Y−Z *with fiber metric,*
- $\nu$ *is an $\mathcal{I}$-valued harmonic 1-form on* Y−Z *whose norm* $|\nu|$ *extends across* Z *as a Hölder continuous function on* Y *and with* $|\nabla\nu|^2$ *having finite integral on* Y.

(IV.6.2)

(Note that a theorem of Zhang [Z?] asserts in part that Z has finite 1-dimensional Hausdorff measure.)

As a final bit of notation: Given $s \in [1, \infty)$, let $\mathfrak{s}: [1,\infty)\times Y \to \mathbb{R}\times Y$ diffeomorphism given by the rule $(t, y) \to (t-s, y)$. (This is a constant translation along the $\mathbb{R}$ factor.)

**Proposition IV.6.1**: *Let* $(A,\mathfrak{a})$ *denote a solution to (I.2.1) on* $[1,\infty)\times Y$ *that obeys the conditions set forth in (IV.6.1). Let* $\{s_n\}_{n\in\mathbb{N}} \subset [1, \infty)\times Y$ *denote an unbounded, increasing sequence. There exists the following data:*

a) *A $\mathbb{Z}/2$ harmonic 1-form data set* $(\mathcal{I}, Z, \nu)$ *on* Y.
b) *A principal $SU(2)$-bundle* $P^\infty \to Y-Z$ *with a flat connection* $A^\infty$ *and an $A^\infty$−covariantly constant isometry* $\sigma: \mathcal{I} \to \mathrm{ad}(P^\infty)$.
c) *A subsequence* $\Lambda \subset \mathbb{N}$.
d) *For each* $n \in \Lambda$, *an isomorphism over* $[1, \infty)\times(Y-Z)$ *from* $P^\infty$ *(pulled-back via the projection map) and* P. *This isomorphism is denoted by* $g_n$.

*These are such that*

- *The pull-back sequence* $\{\mathfrak{s}_n^*g_n^*A\}_{n\in\Lambda}$ *converges weakly to* $A^\infty$ *in the $L^2_1$ topology on compact subsets of* $\mathbb{R}\times(Y-Z)$.
- *The pull-back sequence* $\{\mathfrak{s}_n^*g_n^*\frac{1}{M}\mathfrak{a})\}_{n\in\Lambda}$ *converges to* $\nu\sigma$ *weakly in the $L^2_1$-topology on compact subsets of* $\mathbb{R}\times(Y-Z)$.
- *The sequence* $\{\mathfrak{s}_n^*\frac{1}{M}|\mathfrak{a}|)\}_{n\in\Lambda}$ *converges to* $|\nu|$ *in the $C^0$ topology on compact subsets of* $\mathbb{R}\times Y$ *and weakly in the $L^2_1$-topology on compact subsets of* $\mathbb{R}\times Y$.

*Proof of Proposition IV.6.1*: Because $\mathfrak{cs}$ is bounded, this identity in (IV.2.2) implies this: Given $\varepsilon > 0$, there exists $t_\varepsilon > 1$ such that

$$\int_{[t_\varepsilon,\infty)\times Y} (|E_A|^2 + |\nabla_t\mathfrak{a}|^2) < \varepsilon .$$

(IV.6.3)

This implies in turn (via (I.2.1)) that

$$\int_{[t_\varepsilon,\infty)\times Y} (|B_A - *(\mathfrak{a}\wedge\mathfrak{a})|^2 + |d_A\mathfrak{a}|^2) < \varepsilon$$

(IV.6.4)

also. The inequality (IV.6.3) also implies this: Given $\varepsilon > 0$ and times $t, t´ \in [t_\varepsilon, \infty)$, then

$$|M(t) - M(t´)| < \varepsilon |t-t´|^{1/2}.$$

(IV.6.5)

Now let $(A^n, \mathfrak{a}^n)$ for $n \in \mathbb{Z}$ denote the translated pair $s_n^*(A,\mathfrak{a})$. Since $d_A *\mathfrak{a} = 0$, it is also the case that $d_{A^n} *\mathfrak{a}^n = 0$. With that observation and (IV.6.4) and (IV.6.5), the assumptions in the second bullet of Theorem 1.2 of [T1] are met for any fixed bounded interval $\mathbb{I}$ in $(-\infty, \infty)$ by the sequence $\{(A^{n+N}, \mathfrak{a}^{n+N})\}_{n\geq N}$ with N fixed by the choice of $\mathbb{I}$. Therefore, the conclusions of that second bullet can be invoked. Except for the three points, the conclusions of the second bullet in Theorem 1.2 are the assertions of Proposition IV.6.1. The first point concerns the renormalization $\mathfrak{a} \to \frac{1}{M}\mathfrak{a}$ that is used here. Because of (IV.6.5) and because $\mathbb{I}$ is bounded, this renormalization is compatible with the renormalization used in Theorem 1.2 of [T1] which renormalizes $\mathfrak{a}$ by dividing by its $L^2$ norm over $\mathbb{I}\times Y$. The different renormalizations lead to versions of $\nu$ that are proportional in the sense that one is a constant, non-zero multiple of the other.

The second point is this: Theorem 1.2 in [T1] does not say that its limit $\nu$ is the pull-back via the projection map to Y. That this is so follows from the IV.6.3 given the manner of convergence dictated by Theorem 1.2 in [T1] and by virtue of (IV.6.5). Indeed, (IV.6.5) implies that the version of $\nu$ from Theorem 1.2 in [T1] is bounded on $\mathbb{R}\times Y$. Since it is also closed, and coclosed, an integration by parts can be employed (taking care near Z) to see that the $\partial_t \nu = 0$. Meanwhile, the dt component of $\nu$ is zero because the dt component of $\mathfrak{a}$ is zero. A t-independent 1-form with no dt component is necessarily the pull-back of a form from Y.

The third point is the assertion that the limit connection $A^\infty$ is flat. This follows from (IV.6.4) if the following is true: Let $B \subset Y$ denote a small radius ball with compact closure in the complement of Z. Then for any fixed t,

$$\lim_{n\in \Lambda} \int_{[t+s_n-1, t+s_n+1)\times B} |\mathfrak{a}\wedge\mathfrak{a}|^2 = 0.$$

(IV.6.6)

The proof that (IV.6.6) holds has three parts.

*Part 1*: Let $(A_n, \mathfrak{a}_n)$ denote the pull-back of $(A, \mathfrak{a})$ by the translation $s_n$. Let B´ denote a ball in Y concentric to B with slightly larger radius whose closure is also disjoint from Z. Then set X to denote $[t-2, t+2] \times B´$. Having fixed $\varepsilon > 0$, then the integral of all sufficiently large n versions of $|E_{A_n}|^2 + |B_{A_n} - *(\mathfrak{a}_n \wedge \mathfrak{a}_n)|^2$ over X will be less than $\varepsilon$. Keeping this in mind: Lemmas 6.2–6.4 in [T1] refer to a set $\Theta_c$ that is associated to each $(A_n, \mathfrak{a}_n)$ on the 4-manifold X. Because of the preceding observation, this set will be empty for all sufficiently large n. As a consequence of that, the conclusions of Lemma 6.4 in [T1] hold for $(A_n, \mathfrak{a}_n)$ on the subset $U \equiv [t - \frac{3}{2}, t + \frac{3}{2}] \times B$ in X for all sufficiently large n. This lemma gives a positive lower bound for the $(A_n, \mathfrak{a}_n)$ and $p \in U$ and $c > 1$ versions of numbers $r_{c\wedge}$ and $r_{cF}$ that are defined via (3.1) and (3.2) at the very beginning of Section 3 of [T1]. This lower bound depends on $c$ but it can be chosen so that it is independent of n and p. Let $r_{c*}$ denote such a lower bound with the property that any ball of raduis $r_{c*}$ with center in U has compact closure in X.

*Part 2*: Supposing that $p \in U$ and having fixed a large positive integer (call it n), let $\kappa_n$ and $N_n$ denote the $(A_n, \mathfrak{a}_n)$ versions of the functions that are defined in Section 3a of [T1] with these viewed as functions on $[0, r_{c*}]$. If n is sufficiently large (independent of p), then $\kappa_n$ will be very close at each $r \in [0, r_{c*}]$ to the version of $\kappa$ that is defined as in Section 3a of [T1] using $\nu\sigma$. (That version is denoted by $\kappa_\infty$.) Because of this and because $|\nu|$ is uniformly bounded away from zero on X, one has this: Given $\varepsilon > 0$, all sufficiently large n versions of $\kappa_n$ obey

$$|\kappa_n(r) - \kappa_\infty(r)| \leq (1+\varepsilon)\kappa_\infty(r)$$

(IV.6.7)

for all $r \in (0, r_{c*}]$.

*Part 3*: Having fixed $p \in U$, there is also an $(A^\infty, \nu\sigma)$ version of the function $N$ that is defined in Section 3a of [T1]. It is $r\frac{d}{dr}\ln\kappa_\infty$. This version of $N$ is denoted by $N_\infty$. It is bounded by a constant multiple of r at each point in U with the constant being independent of the given point. This can be proved using Taylor's theorem with remainder applied to $\nu$ because $\nu$ is nothing but a smooth, harmonic 1-form on X.

Granted that $N_\infty(r)$ is bounded by a constant multiple of r, and granted that $N_n$ is defined by writing it as $r\frac{d}{dr}\ln\kappa_n$ (see (3.6) in [T1]) then (IV.6.7) implies the following: Given $\varepsilon > 0$, there exists $r_\varepsilon > 0$ (which is independent of $p \in U$) such that if n is sufficiently large so that (IV.6.7) holds on $(0, r_{c*}]$, and if $r < r_\varepsilon$, then

$$N_n(\cdot) < 10^4 \varepsilon.$$

(IV.6.8)

on a subset of $[r, 2r]$ with measure greater than $\frac{9}{10}r$.

*Part 4*: The preceding observations about $N_n$ for large n implies this: There exists $r_{**} \in (0, r_{c*}]$ such that if n is sufficiently large and if $p \in U$, then the conclusions of Propositions 3.1-3.3 in [T1] hold for $(A_n, \mathfrak{a}_n)$ and the point p if $r_c$ is chosen judiciously from any interval of the form $[r, 2r]$ in $(0, r_{**}]$. In this instance, what is said by the $(A_n, \mathfrak{a}_n)$ version of Lemma 4.3 in [T1] applies for all $p \in U$ if n is sufficiently large if $\mu = \frac{1}{2}$. That lemma with the $(A_n, \mathfrak{a}_n)$ version of (4.5) in [T1] implies in turn that the integral over the radius $\frac{1}{2} r_c$ ball centered at p of $|\mathfrak{a} \wedge \mathfrak{a}|^4$ is bounded by an n and p independent multiple of $K_n^{-2}(t + s_n)$. The latter bound leads directly to (IV.6.6).

## A. APPENDIX

The following are versions of the γ and ρ matrices that are introduced in [I.1.7].

$$\gamma_1 \equiv \begin{pmatrix} 0 & 0 & 0 & -1 & 0 & 0 & 0 & 0 \\ 0 & 0 & 0 & 0 & 0 & 0 & 1 & 0 \\ 0 & 0 & 0 & 0 & 0 & -1 & 0 & 0 \\ 1 & 0 & 0 & 0 & 0 & 0 & 0 & 0 \\ 0 & 0 & 0 & 0 & 0 & 0 & 0 & -1 \\ 0 & 0 & 1 & 0 & 0 & 0 & 0 & 0 \\ 0 & -1 & 0 & 0 & 0 & 0 & 0 & 0 \\ 0 & 0 & 0 & 0 & 1 & 0 & 0 & 0 \end{pmatrix} \quad \gamma_2 \equiv \begin{pmatrix} 0 & 0 & 0 & 0 & 0 & 0 & -1 & 0 \\ 0 & 0 & 0 & -1 & 0 & 0 & 0 & 0 \\ 0 & 0 & 0 & 0 & 1 & 0 & 0 & 0 \\ 0 & 1 & 0 & 0 & 0 & 0 & 0 & 0 \\ 0 & 0 & -1 & 0 & 0 & 0 & 0 & 0 \\ 0 & 0 & 0 & 0 & 0 & 0 & 0 & -1 \\ 1 & 0 & 0 & 0 & 0 & 0 & 0 & 0 \\ 0 & 0 & 0 & 0 & 0 & 1 & 0 & 0 \end{pmatrix} \quad \gamma_3 \equiv \begin{pmatrix} 0 & 0 & 0 & 0 & 0 & 1 & 0 & 0 \\ 0 & 0 & 0 & 0 & -1 & 0 & 0 & 0 \\ 0 & 0 & 0 & -1 & 0 & 0 & 0 & 0 \\ 0 & 0 & 1 & 0 & 0 & 0 & 0 & 0 \\ 0 & 1 & 0 & 0 & 0 & 0 & 0 & 0 \\ -1 & 0 & 0 & 0 & 0 & 0 & 0 & 0 \\ 0 & 0 & 0 & 0 & 0 & 0 & 0 & -1 \\ 0 & 0 & 0 & 0 & 0 & 0 & 1 & 0 \end{pmatrix}$$

$$\rho_1 \equiv \begin{pmatrix} 0 & 0 & 0 & 0 & 0 & 0 & 0 & 1 \\ 0 & 0 & 1 & 0 & 0 & 0 & 0 & 0 \\ 0 & -1 & 0 & 0 & 0 & 0 & 0 & 0 \\ 0 & 0 & 0 & 0 & 1 & 0 & 0 & 0 \\ 0 & 0 & 0 & -1 & 0 & 0 & 0 & 0 \\ 0 & 0 & 0 & 0 & 0 & 0 & -1 & 0 \\ 0 & 0 & 0 & 0 & 0 & 1 & 0 & 0 \\ -1 & 0 & 0 & 0 & 0 & 0 & 0 & 0 \end{pmatrix} \quad \rho_2 \equiv \begin{pmatrix} 0 & 0 & -1 & 0 & 0 & 0 & 0 & 0 \\ 0 & 0 & 0 & 0 & 0 & 0 & 0 & 1 \\ 1 & 0 & 0 & 0 & 0 & 0 & 0 & 0 \\ 0 & 0 & 0 & 0 & 0 & 1 & 0 & 0 \\ 0 & 0 & 0 & 0 & 0 & 0 & 1 & 0 \\ 0 & 0 & 0 & -1 & 0 & 0 & 0 & 0 \\ 0 & 0 & 0 & 0 & -1 & 0 & 0 & 0 \\ 0 & -1 & 0 & 0 & 0 & 0 & 0 & 0 \end{pmatrix} \quad \rho_3 \equiv \begin{pmatrix} 0 & 1 & 0 & 0 & 0 & 0 & 0 & 0 \\ -1 & 0 & 0 & 0 & 0 & 0 & 0 & 0 \\ 0 & 0 & 0 & 0 & 0 & 0 & 0 & 1 \\ 0 & 0 & 0 & 0 & 0 & 0 & 1 & 0 \\ 0 & 0 & 0 & 0 & 0 & -1 & 0 & 0 \\ 0 & 0 & 0 & 0 & 1 & 0 & 0 & 0 \\ 0 & 0 & 0 & -1 & 0 & 0 & 0 & 0 \\ 0 & 0 & -1 & 0 & 0 & 0 & 0 & 0 \end{pmatrix}$$